\documentclass[12pt]{amsart}
\usepackage{amsmath,amssymb,amscd}
\usepackage{times,color, graphicx}
\newcommand{\beqn}{\begin{eqnarray}}
\newcommand{\eeqn}{\end{eqnarray}}
\newcommand{\W}{\pmb{W}}
\newcommand{\vecu}{\pmb{u}}
\newcommand{\vecv}{\pmb{v}}
\newcommand{\R}{\pmb{R}}
\newcommand{\C}{\mathbb{C}}
\newcommand{\N}{\mathbb{N}}

\newcommand{\Z}{\mathbb{Z}}
\newcommand{\Q}{\mathbb{Q}}

\newcommand{\s}{\mathfrak S}
\newcommand{\F}{\mathcal{F}}
\newcommand{\I}{\pmb{I}}
\newcommand{\invlim}{\underleftarrow{\lim}}
\newcommand{\eqcoh}{H^*_{T_n}}
\newcommand{\cS}{\mathcal{S}}

\newcommand{\e}{\pmb{e}}
\newcommand{\diag}{{\mathrm{diag}}}

\theoremstyle{plain}
\newtheorem{thm}{Theorem}[section]
\newtheorem{Def}[thm]{Definition}
\newtheorem{lem}[thm]{Lemma}
\newtheorem{sublem}[thm]{Sublemma}
\newtheorem{prop}[thm]{Proposition}
\newtheorem{cor}[thm]{Corollary}

\newtheorem{remark}[thm]{Remark}

\def\1{\bar{1}}
\def\2{\bar{2}}
\def\3{\bar{3}}

\begin{document}\title{Double Schubert polynomials 
for the classical groups}\author{Takeshi Ikeda} \address{Department of Applied Mathematics,                                            
Okayama University of Science, Okayama 700-0005, Japan} \email{ike@xmath.ous.ac.jp}

\author{Leonardo~C.~Mihalcea}
\address{Department of Mathematics, Duke University, P.O. Box 90320 Durham, NC 27708-0320 USA}
\email{lmihalce@math.duke.edu}

\author{Hiroshi Naruse}\address{Graduate School of Education, Okayama University, Okayama 700-8530, Japan}\email{rdcv1654@cc.okayama-u.ac.jp}


\date{October 5, 2008}

\maketitle

\begin{abstract}
For each infinite series of the classical Lie groups of type $\mathrm{B}$, $\mathrm{C}$ or $\mathrm{D}$,
we introduce a family of polynomials parametrized by the elements of the
corresponding Weyl group of infinite rank.
These polynomials represent the Schubert classes in the equivariant 
cohomology of the appropriate flag variety. They satisfy a stability property,
and are a natural extension of the (single) Schubert polynomials of Billey 
and Haiman, which represent non-equivariant 
Schubert classes. They are also positive in a certain sense, and when indexed by maximal Grassmannian elements, or by the longest element in a finite Weyl group, these polynomials can be expressed in terms of the factorial analogues of
Schur's $Q$- or $P$-functions defined earlier
by Ivanov.  
\end{abstract}


\section{Introduction}


A classical result of Borel states that
the (rational) cohomology ring of 
the flag variety of any 
simple complex Lie group $G$
is isomorphic, as a graded ring, to the {\it coinvariant\/} algebra
of the corresponding Weyl group $W$,
i.e. to the quotient of a polynomial
ring modulo the ideal generated by 
the $W$-invariant polynomials
of positive degree.
The Schubert classes form a 
distinguished additive basis of the cohomology ring
indexed by the elements in the Weyl group.
Bernstein-Gelfand-Gelfand \cite{BGG} (see also Demazure \cite{D}) 
showed that if one starts with a polynomial that represents the 
cohomology class of
highest codimension (the Schubert 
class of a point), one obtains all the other Schubert
classes by applying a succession of 
{\it divided difference\/} operators
corresponding to simple roots. This construction depends on 
the choice of a polynomial representative for the  ``top'' cohomology class.
For $SL(n,\C)$, Lascoux and Sch\"utzenberger \cite{LS} considered one particular choice, which yielded polynomials - the Schubert polynomials -
with particularly good combinatorial and geometric
properties.

It is a natural problem to extend 
the construction in \cite{LS} to the other
classical Lie groups. 
To this end, Fomin and Kirillov
\cite{FK} listed up five properties that characterize
the Schubert polynomials in type $\mathrm{A},$
but they showed that it is impossible to 
construct a theory of ``Schubert polynomials''
in type $\mathrm{B}$ satisfying the same properties.
For type $\mathrm{B}_{n}$, they constructed several families
of polynomials which satisfy all but one of these properties.

There is another approach to this 
problem due to Billey and Haiman
\cite{BH}. Consider one of the series of Lie types
$\mathrm{B}_{n},\mathrm{C}_{n}$, and $\mathrm{D}_{n}$
and denote by $G_{n},T_{n},$ and $B_{n}$
the corresponding classical Lie group,
a maximal torus, and a Borel subgroup
containing the maximal torus. The associated flag variety is $\mathcal{F}_{n}=G_{n}/B_{n}$, and the cohomology Schubert classes $\sigma_{w}^{(n)}$ are labeled by elements $w$ in the  Weyl group
$W_{n}$ of $G_n$. There is a natural embedding of groups $G_{n}\hookrightarrow G_{n+1}$, and this induces an embedding
of the flag varieties
$\mathcal{F}_{n}\hookrightarrow \mathcal{F}_{n+1}$ 
and reverse maps in cohomology $H^*(\F_{n+1},\Q) \to H^{*}(\F_n,\Q)$, compatible with the Schubert classes.
For each element $w$ in the infinite Weyl group $W_\infty = \bigcup_{n \ge 1} W_n$, there is 
a {\it stable\/} Schubert class 
$\sigma_{w}^{(\infty)}=\invlim \,
\sigma_{w}^{(n)}$
in the inverse system
$\invlim \,H^{*}(\mathcal{F}_{n},\Q)$
of the cohomology rings. 
A priori, the class $\sigma_{w}^{(\infty)}$
is represented by a homogeneous element in the ring of 
power series $\Q[[z_1, z_2, \dots ]]$, but
Billey and Haiman showed in \cite{BH} that 
it is represented by
a {\em unique} element $\mathfrak{S}_{w}$
in the subring\footnote{
The elements $\{z_{1},z_{2},\ldots;p_{1}(z),p_{3}(z),p_{5}(z),\ldots\}$ are algebraically independent, so
the ring (\ref{eq:BHring}) can also be 
regarded simply as the polynomial
ring in $z_{i}$ and $p_{k}$ ($k=1,3,5,\ldots)$. }
\begin{equation}
\mathbb{Q}[z_1,z_2, \ldots;
p_{1}(z),p_{3}(z),p_{5}(z),\ldots],\label{eq:BHring}
\end{equation}
where $p_{k}(z)=\sum_{i=1}^{\infty}z_{i}^{k}$
denotes the power-sum
 symmetric function. 
Note that 
the images of the even power-sums $p_{2i}(z)$ in the limit of the coinvariant rings
vanish for each of the types $\mathrm{B,C}$ and $\mathrm{D}$.
The elements $\{\mathfrak{S}_{w}\}$ are obtained as
the {\em unique} solutions of a system of equations involving infinitely many
divided difference operators.
These polynomials will satisfy the main combinatorial properties of the type $\mathrm{A}$ Schubert polynomials, if interpreted appropriately.
In particular, the polynomial 
$\mathfrak{S}_{w}$ 
is stable, i.e. 
it represents the Schubert classes
$\sigma_{w}^{(n)}$ 
in $H^*(\F_n,\Q)$ simultaneously for all positive integers $n.$

The flag varieties admit an action of
the maximal torus $T_{n}$, and the inclusion $\F_n \hookrightarrow \F_{n+1}$ is equivariant
with respect to $T_{n}\hookrightarrow T_{n+1}
$. Therefore one can define an {\em equivariant} version 
of the stable Schubert classes, and one can ask whether 
we can `lift' the polynomials of Billey and Haiman
to the equivariant setting.
These will be the ``double Schubert polynomials", which we will define and study next.
The terminology comes from type $\mathrm{A}$, where Lascoux and Sch\"utzenberger defined a double version of their Schubert polynomials in \cite{LS} (see also \cite{MacSchubert}).

As shown by Fulton in \cite{F:deg}, the type $\mathrm{A}$ double Schubert polynomials can also be constructed as polynomials which represent the cohomology classes of some degeneracy loci. More recently, two related constructions connecting double Schubert polynomials to equivariant cohomology of flag manifolds, using either Thom polynomials or Gr\"obner degenerations were obtained independently by Feh\'er and Rim\'anyi \cite{FR} and by Knutson and Miller \cite{KM}. The degeneracy locus construction was extended to other types by Fulton \cite{F}, Pragacz-Ratajski \cite{PR} and Kresch-Tamvakis \cite{KT}. The resulting polynomials are expressed in terms of Chern classes canonically associated to the geometric situation at hand.
Their construction depends again on the choice of a 
polynomial to represent the ``top class" - the diagonal class 
in the cohomology of a flag bundle.
Unfortunately 
different choices lead to polynomials having some desirable combinatorial properties - but not all. In particular, the polynomials in \cite{KT,F} do not satisfy the
stability property.

In this paper we will work in the equivariant cohomology of flag varieties.
As in \cite{BH}, there is a unique family of stable polynomials,
which is the unique solution of a system of divided difference operators. In our study we will make full use of {\it localization\/} techniques in equivariant cohomology. In the process, we will reprove, and put on a more solid geometric foundation, the results from \cite{BH}.

\subsection{Infinite hyperoctahedral groups}\label{ssec:infinite_weyl} To fix notations, let $W_\infty$ be the infinite hyperoctahedral group, i.e.~the Weyl group of type $\mathrm{C}_\infty$ (or $\mathrm{B}_\infty$). It is generated by elements $s_0,s_1, \ldots$ subject to the braid relations described in (\ref{eq:Wrel}) below. For each nonnegative integer $n$, the subgroup $W_n$ of $W_\infty$ generated by $s_0, \ldots , s_{n-1}$ is the Weyl group of type $\mathrm{C}_n$. $W_\infty$ contains a distinguished subgroup $W_\infty'$ of index $2$ - the Weyl group of $\mathrm{D}_\infty$ -  which is generated by $s_{\hat{1}},s_{1}, s_2, \ldots$, where $s_{\hat{1}} = s_0s_1s_0$. The corresponding finite subgroup $W_n'=W_\infty' \cap W_n$ is the type $\mathrm{D}_n$ Weyl group. To be able to make statements which are uniform across all classical types, we use $\W_{\infty}$ to denote $W_{\infty}$ when we consider types $\mathrm{C}$ or $\mathrm{B}$ and
$W_{\infty}'$ for type $\mathrm{D}$; similar notation is used for $\W_{n}\subset \W_{\infty}$. Finally,
set $\I_{\infty}$ to be the indexing set $\{0,1,2,\ldots\}$ 
for types $\mathrm{B,C}$ and $\{\hat{1},1,2,\ldots\}$ for type $\mathrm{D}.$



\subsection{Stable Schubert classes}\label{ssec:SSch} To each element $w \in \W_n$ there corresponds a torus-fixed point $e_{w}$ in the flag variety $\mathcal{F}_{n}=G_{n}/B_{n}$ (see \S \ref{ssec:NotationFlag} below).
The {\em Schubert variety}
$X_{w}$ is defined to be the closure of the {\em Schubert cell}
$B^{-}_{n}e_{w}$ in $\mathcal{F}_{n},$
where $B^{-}_{n}$ is the Borel subgroup opposite
to $B_{n}$. The fundamental class of $X_w$ determines a {\em Schubert class}
$\sigma_{w}^{(n)}$ in the
equivariant cohomology ring\begin{footnote}{
Unless otherwise stated, from now on we work over cohomology with coefficients over $\Z$.}\end{footnote} $H_{T_{n}}^{2 \ell(w)}(\mathcal{F}_{n})$, where $\ell(w)$ is the length of $w$. The classes $\{\sigma_{w}^{(n)}\}_{w\in \W_{n}}$ 
form an $H_{T_{n}}^{*}(pt)$-basis
of the equivariant cohomology ring
$H_{T_{n}}^{*}(\mathcal{F}_{n}).$
Note that $H_{T_{n}}^{*}(pt)$ is canonically
identified with the polynomial ring
$\Z[t]^{(n)}:=\Z[t_{1},\ldots,t_{n}]$ generated by 
a standard basis $\{t_{1},\ldots,t_{n}\}$ of the character group of $T_{n}.$

Since the torus $T_n$ acts trivially on $e_v$, the inclusion map $\iota_v:e_v \to \F_n$ is equivariant, and induces the {\em localization map} $\iota_{v}^{*}: H_{T_{n}}^{*}(\mathcal{F}_{n})
\rightarrow H_{T_{n}}^{*}(e_{v})$.
It is well-known (cf. e.g. \cite{A}) that 
the product map
$$
\iota^{*}=(\iota_{v}^{*})_{v}: H_{T_{n}}^{*}(\mathcal{F}_{n})
\longrightarrow
\prod_{v\in \W_{n}}
H_{T_{n}}^{*}(e_{v})
=\prod_{v\in \W_{n}}
\Z[t]^{(n)}$$
is injective, so we will often 
identify $\sigma_{w}^{(n)}$ with an
element in $\prod_{v\in \W_{n}}
\Z[t]^{(n)}$ via $\iota^{*}.$
It turns out that the localizations of Schubert classes stabilize, in the sense that for $v,w \in \W_n \subset \W_m$, the polynomial $\iota_v^*(\sigma_w^{(m)})$ in $H^*_{T_{m}}(e_v)=\Z[t]^{(m)}$ remains constant as $m$ varies.
Therefore, we can pass to the limit to define the {\em stable (equivariant) Schubert class} $\sigma_{w}^{(\infty)}
$ 
in $\prod_{v\in \W_{\infty}}\Z[t]$, 
where $\Z[t]$ is the polynomial ring
in the variables $t_i$ ($i \ge 1$).
Denote by $H_{\infty}$ the $\Z[t]$-submodule of $\prod_{v\in \W_{\infty}}\Z[t]$ spanned by the stable Schubert classes. We will show that $H_{\infty}$ is actually a subalgebra of $\prod_{v\in \W_{\infty}}\Z[t]$ which has a $\Z[t]$-basis consisting of stable Schubert classes $\{\sigma_{w}^{(\infty)}\}.$

 
A crucial part in the theory of Schubert polynomials of classical types is played by the $P$ and $Q$ Schur functions \cite{Schur}. These are symmetric functions $P_\lambda(x)$ and $Q_\lambda(x)$ in a new set of variables $x=(x_1, x_2, \ldots)$, and are indexed by strict partitions $\lambda$ (see \S \ref{sec:Schur} below for details). The $P$ or $Q$ Schur function corresponding to $\lambda$ with one part of length $i$ is denoted respectively by $P_i(x)$ and $Q_i(x)$. 
Define $\Gamma=\Z[Q_{1}(x),Q_{2}(x),\ldots]$ and 
$\Gamma'=\Z[P_{1}(x),P_{2}(x),\ldots]$. Note that $\Gamma$ and $\Gamma'$ are not polynomial rings, since $Q_i(x)$ respectively $P_i(x)$ are not algebraically independent (see \S \ref{sec:Schur} for the relations among them), but they have canonical $\Z$-bases consisting of the $Q$-Schur functions $Q_\lambda(x)$ (respectively $P$-Schur functions $P_\lambda(x)$). 
 We define next the $\Z[t]$-algebras of {\em Schubert polynomials} $$R_{\infty}=\Z[t]\otimes_{\Z}\Gamma\otimes_{\Z}\Z[z],\quad
R'_{\infty}=\Z[t]\otimes_{\Z}\Gamma'\otimes_{\Z}\Z[z],
$$
where $\Z[z]=\Z[z_{1},z_{2},\ldots]$ is the polynomial
ring in $z=(z_{1},z_{2},\ldots).$ We will justify the terminology in the next paragraph.
Again, in order to state results uniformly in all types, we use the bold letter $\R_{\infty}$ 
to denote $R_{\infty}$ for type $\mathrm{C}$ and $R'_{\infty}$
for types $\mathrm{B}$ and $\mathrm{D}$.

There exists a homomorphism 
$$
\Phi=(\Phi_{v})_{v\in \W_{\infty}} : \R_{\infty}\longrightarrow 
\prod_{v\in \W_{\infty}}\Z[t]
$$
of graded $\Z[t]$-algebras,
which we call {\em universal localization map}. Its precise (algebraic) definition is given in \S \ref{sec:UnivLoc} and it has a natural geometrical interpretation explained in \S \ref{sec:geometry}. One of the main results in this paper is that $\Phi$ is an isomorphism
from $\R_{\infty}$ onto $H_{\infty}$
(cf. Theorem \ref{thm:isom} below). While injectivity is easily proved algebraically, surjectivity is more subtle. It uses the {\em transition equations}, which are recursive formulas which allow writing any stable Schubert class $\sigma_w^{(\infty)}$ in terms of 
the (Lagrangian or Orthogonal) Grassmannian Schubert classes. Once reduced to the Grassmannian case, earlier results of the first and third author \cite{Ik,IN}, which show that the classes in question are represented by Ivanov's {\em factorial $Q$ (or $P$)} Schur functions \cite{Iv}, finish the proof of surjectivity.



By pulling back - via $\Phi$ - the stable Schubert classes,  
we introduce polynomials $\mathfrak{S}_{w}
=\mathfrak{S}_{w}(z,t;x)$ in $\R_{\infty}$, which are uniquely determined by
the "localization equations"
\begin{equation}
\Phi_{v}\left(
\mathfrak{S}_{w}(z,t;x)\right)
=\iota_{v}^{*}(\sigma_{w}^{(\infty)}),
\label{eq:Phisigma}
\end{equation}
where 
$\iota_{v}^{*}(\sigma_{w}^{(\infty)})
\in \Z[t]$ is the stable limit
of $\iota_{v}^{*}(\sigma_{w}^{(n)})$.

\subsection{Divided difference operators}
Alternatively, $\{\mathfrak{S}_{w}(z,t;x)\}$ can be characterized
in a purely algebraic manner by using the divided difference operators.
There are two families of operators $\partial_i, \delta_i$ ($i \in \I_\infty$) on $\R_\infty$, such that operators from one family commute with those from the other (see \S \ref{ssec:divdiff}
for the definition). Then:
\begin{thm}
\label{existC}
There exists a unique family of elements $\mathfrak{\s}_{w}= \s_w(z,t;x) $ in $\R_\infty$, where $w\in \W_{\infty}$,
satisfying the equations
\begin{equation}\label{E:dd} \partial_{i}\mathfrak{S}_{w}=\begin{cases}
\mathfrak{S}_{ws_{i}}\quad\mbox{if}\quad \ell(ws_{i})<\ell(w)\\
0 \quad \mbox{otherwise}
\end{cases},\quad
\delta_{i}\mathfrak{S}_{w}=\begin{cases}
\mathfrak{S}_{s_{i}w}\quad\mbox{if}\quad \ell(s_{i}w)<\ell(w)\\
0 \quad \mbox{otherwise}
\end{cases},
\end{equation}
for all $i\in \I_{\infty}$,
and such that $\mathfrak{S}_{w}$ has no
constant term except for 
$\mathfrak{S}_{e}=1$. 
\end{thm}
The operators $\partial_{i},\delta_{i}$
are the limits of the same operators on the equivariant
cohomology $\eqcoh (\F_n)$, since the latter are compatible with the projections $H^*_{T_{n+1}}(\F_{n+1}) \to \eqcoh(\F_n)$. In this context, the operator $\partial_i$ is an equivariant generalization the operator defined in \cite{BGG,D}, and it can shown that it is induced by the right action of the Weyl group on the equivariant cohomology (cf. \cite{KK,Kn}). The operator $\delta_i$ exists only in equivariant cohomology, and it was used in \cite
{Kn,Ty} to study equivariant Schubert classes. It turns out that it corresponds to a left Weyl group action on $\eqcoh (\F_n)$.

\subsection{Billey and Haiman's polynomials and a change of variables}
The polynomials from Theorem \ref{existC} are lifts of Billey-Haiman polynomials from the non-equivariant
to equivariant cohomology.
Concretely, if we forget the torus action,
then
$$
\mathfrak{S}_{w}(z,0;x)=\mathfrak{S}_{w}^{BH}(z;x)
$$
where $\mathfrak{S}_{w}^{BH}(z;x)$ denotes
the Billey-Haiman's polynomial. This follows immediately from the Theorem \ref{existC}, since $\{\mathfrak{S}_{w}^{BH}(z;x)\}$
are the unique elements in
$\Gamma\otimes_{\Z}\Z[z]$ (for type $\mathrm{C}$)
or $\Gamma'\otimes_{\Z}\Z[z]$ (for type $\mathrm{B,D}$)
which satisfy the equation involving the right divided difference operators
$\partial_{i}.$

The variables $z_i$ in $\mathfrak{S}_w(z,t;x)$ correspond geometrically to the limits of Chern classes of the tautological line bundles, while the variables $t_i$ to the equivariant parameters. To understand why the variables $x_i$ are needed - both algebraically and geometrically - we comment next on  the ``change of variables'' which relates the Billey-Haiman polynomials to those defined by Fomin and Kirrilov in \cite[Thm. 7.3]{FK}. The general formula in our situation - along with its geometrical explanation - will be given in section \ref{sec:geometry} below. The relation between
$x$ and $z$ is given by
$$
\prod_{i=1}^{\infty}\frac{1+x_{i}u}{1-x_{i}u}
=\prod_{i=1}^{\infty}(1-z_{i}u),
$$ 
or equivalently
$p_{2i}(z)=0$ and $-2p_{2i+1}(x)=p_{2i+1}(z).$ In type $\mathrm{C}$, it is known that $p_{2i}(z)$ generates the (limit of the) ideal of relations in cohomology, therefore such a variable change eliminates the ambiguity of representatives
coming from $p_{2i}(z)=0.$ Note that the change of variables can be also
expressed as $(-1)^{i}e_{i}(z)=Q_{i}(x)=2P_{i}(x)\;(i\geq 1).$
It follows that after extending the scalars from $\Z$ to $\Q$, the ring $\Gamma\otimes_{\Z}\Z[z]$
(or $\Gamma'\otimes_{\Z}\Z[z]$) is identified with 
the ring (\ref{eq:BHring}).
Since both $\Gamma$ and $\Gamma'$ have
distinguished $\Z$-bases,
the polynomials $\mathfrak{S}^{BH}_{w}(z;x)$ will expand uniquely as a combination of $Q$-Schur (or $P$-Schur) functions with coefficients in $\Q[z]$.

In type $\mathrm{C}$, the change of variables corresponds to making $Q_i(x)=c_i(\mathcal{S}^*)$ - where $c_i(\mathcal{S}^*)$ is the limit of the Chern classes of the duals of the tautological subbundles of the Lagrangian Grassmannians, regarded as elements in $\invlim H^*(\F_n)$. This is the identification which was used by Pragacz (see e.g. \cite[pag. 32]{FP}) to study the cohomology of the Lagrangian Grassmannian.



\subsection{Combinatorial properties of the double Schubert polynomials} We state next
the combinatorial properties of the double Schubert polynomials $\s_w(z,t;x)$:
\begin{itemize}
\item(Basis) $\s_w(z,t;x)$ form a $\Z[t]$-basis of $\R_\infty$;
\item(Symmetry) $\s_w(z,t;x) = \s_{w^{-1}}(-t,-z;x)$;
\item(Positivity) The double Schubert polynomial $\s_w(z,t;x)$ can be uniquely written as \[ \s_w(z,t;x) = \sum_{\lambda} f_\lambda(z,t) F_\lambda(x) \/, \] where the sum is over strict partitions $\lambda=(\lambda_1, \dots , \lambda_r)$ such that $\lambda_1+ \dots +\lambda_r \le \ell(w)$, $f_\lambda(z,t)$ is a homogeneous polynomial in $\N[z,-t]$, 
and $F_\lambda(x)$ is the $Q$-Schur function $Q_\lambda(x)$ in type $\mathrm{C}$, respectively the $P$-Schur function $P_\lambda(x)$ in types $\mathrm{B,D}$.
For a precise combinatorial formula for the 
coefficients $f_\lambda(z,t)$ see Cor. \ref{cor:typeAexpand} and Lem. \ref{lem:S00} below.
\end{itemize} 

The basis property implies that we can define the structure constants 
$c_{uv}^{w}(t)\in \Z[t]$ by the
expansion 
$$
\mathfrak{S}_{u}(z,t;x)
\mathfrak{S}_{v}(z,t;x)
=\sum_{w \in \W_{\infty}}c_{uv}^{w}(t)\mathfrak{S}_{w}(z,t;x).
$$ These coincide with the structure constants in equivariant cohomology of $\F_n$, written in a stable form. The same phenomenon happens in \cite{LS,BH}.

\subsection{Grassmannian Schubert classes}
To each strict partition $\lambda$ one can associate a {\em Grassmannian element} $w_\lambda \in \W_\infty$. 
Geometrically these arise as the elements 
in $\W_{\infty}$ which index the pull-backs of the Schubert classes from the appropriate Lagrangian or Orthogonal Grassmannian, via the natural projection from the flag variety. For the Lagrangian Grassmannian,
the first author \cite{Ik} identified the equivariant Schubert
classes with the factorial analogues
of Schur $Q$-function defined by Ivanov \cite{Iv}.
This result was extended 
to the maximal isotropic
Grassmannians of orthogonal types $\mathrm{B}$ and $\mathrm{D}$
by Ikeda and Naruse \cite{IN}.
See \S \ref{sec:Schur} for the definition 
of Ivanov's functions 
$Q_{\lambda}(x|t),P_{\lambda}(x|t).$
We only mention here that if
all $t_{i}=0$ they coincide with the ordinary $Q$ or $P$
functions;
in that case, these results recover Pragacz's results from \cite{Pr} (see also \cite{Jo}).
We will show in Theorem \ref{PhiFacQ}
that the polynomial $\s_{w_\lambda}(z,t;x)$ coincides with $Q_\lambda(x|t)$
or $P_\lambda(x|t)$, depending on the type at hand.
In particular, the double Schubert polynomials
for the Grassmannian elements
are Pfaffians - this is a {\em Giambelli formula} in this case.


\subsection{Longest element formulas}
Next we present the combinatorial formula for the double Schubert polynomial indexed by $w_0^{(n)}$, the longest element in $\W_n$ (regarded as a subgroup of $\W_\infty$). This formula has a particular significance since this is the top class mentioned in the first section.
We denote by $\mathfrak{B}_{w},\mathfrak{C}_{w},
\mathfrak{D}_{w}$ the double Schubert polynomial
$\mathfrak{S}_{w}$ for types $\mathrm{B,C},$ and $\mathrm{D}$ respectively.
Note that $\mathfrak{B}_{w}=2^{-s(w)}
\mathfrak{C}_{w},$ where $s(w)$ is the number of
signs changed by $w$ (cf. \S \ref{sec:signed.permut} below).
\begin{thm}[Top classes]\label{thm:Top} The double Schubert polynomial 
associated with the longest element $w_{0}^{(n)}$ in $\W_{n}$
is equal to:
\begin{enumerate}
\item 
$\mathfrak{C}_{w_{0}^{(n)}}(z,t;x) =
Q_{\rho_{n}+\rho_{n-1}}
(x|t_{1},-z_{1},t_{2},-z_{2},\ldots,t_{n-1},-z_{n-1}),\label{longC}$
\item $\mathfrak{D}_{w_{0}^{(n)}}(z,t;x) =P_{2\rho_{n-1}}
(x|t_{1},-z_{1},t_{2},-z_{2},\ldots,t_{n-1},-z_{n-1}),\label{longD}$
\end{enumerate} 
where $\rho_{k}=(k,k-1,\ldots,1).$
\end{thm}

\subsection{Comparison with degeneracy loci formulas}
One motivation for the present paper was to give a 
geometric interpretation 
to the factorial Schur $Q$-function by means of degeneracy loci formulas.
In type $\mathrm{A}$,
this problem was treated by the second author in \cite{Mi},
where the Kempf-Laksov formula for degeneracy loci 
is identified with the Jacobi-Trudi type
formula for the factorial (ordinary) Schur function. To this end, 
we will reprove a multi-Pfaffian expression for $\sigma_{w_\lambda}$ (see \S \ref{sec:Kazarian} below) obtained by Kazarian \cite{Ka} while studying Lagrangian degeneracy loci.

\subsection{Organization}
Section \ref{sec:EqSch} is devoted to some general
facts about the equivariant cohomology of the 
flag variety.
In section \ref{sec:Classical} we fix notation concerning
root systems and Weyl groups, while in section \ref{sec:Schur} we give the definitions and some 
properties of 
$Q$- and $P$-Schur functions, and of their factorial analogues. 
The stable (equivariant) Schubert classes
$\{\sigma_{w}^{(\infty)}\}$
and the ring $H_{\infty}$ spanned by these classes are introduced in section \ref{sec:StableSchubert}.
In section \ref{sec:UnivLoc} we define the ring of Schubert polynomials $\R_{\infty}$
and establish the isomorphism
$\Phi: \R_{\infty}\rightarrow H_{\infty}.$
In the course of the proof, we recall 
the previous results on isotropic Grassmannians
(Theorem \ref{PhiFacQ}). 
In section \ref{sec:WactsR} we define the left and right action of the infinite Weyl group on ring
$\R_{\infty},$ and then use them to define the divided
difference operators. 
We also discuss the compatibility
of the actions on both $\R_{\infty}$ and $H_{\infty}$
under the isomorphism $\Phi.$
We will prove the existence and 
uniqueness theorem for the double Schubert polynomials
in section \ref{sec:DSP}, along with some basic combinatorial properties 
of them. The formula for the Schubert polynomials indexed by the longest Weyl group element is proved in
section \ref{sec:Long}.
Finally, in section \ref{sec:geometry} we give an alternative geometric
construction of our universal localization map $\Phi$, and
in section \ref{sec:Kazarian}, we prove the formula for $Q_{\lambda}(x|t)$
in terms of a multi-Pfaffian.

\subsection{Note}
After the present work was completed 
we were informed that 
A. Kirillov \cite{K} had introduced
double Schubert polynomials
of type $\mathrm{B}$ (and $\mathrm{C}$) in 1994
by using Yang-Baxter operators (cf. \cite{FK}), independently to us,
although no connection with (equivariant) cohomology had been established.
His approach is quite different from ours, nevertheless  
the polynomials are the same, after a suitable identification
of variables. Details will be given elsewhere.

\bigskip

This is the full paper version of 
`extended abstract' \cite{fpsac} for the FPSAC 2008
conference held in Vi\~na del Mar, Chile, June 2008.
Some results in this paper 
were announced without proof in \cite{rims}.

\subsection{Acknowledgements}
We would like to thank S.~Billey and H.~Tamvakis for stimulating conversations 
that motivated this work, and to S. Kumar, K.~Kuribayashi, M.~Mimura, M.~Nakagawa, T. Ohmoto, N.~Yagita and M. Yasuo for helpful comments. This work was facilitated by the "Workshop on Contemporary Schubert Calculus and Schubert Geometry" organized at Banff in March 2007. We are grateful to the organizers J. Carrell and F. Sottile for inviting all the authors there.
\section{Equivariant Schubert classes of the flag variety}
\label{sec:EqSch}
\setcounter{equation}{0}
In this section we will recall some basic facts about the equivariant cohomology of the flag variety $\F=G/B$. 
The main references are \cite{A} and \cite{KK} (see also \cite{Ku}).
\subsection{Schubert varieties and equivariant cohomology}\label{ssec:NotationFlag}
Let $G$ be a complex connected semisimple Lie group,
$T$ a maximal torus, $W=N_{G}(T)/T$ its Weyl group,
and $B$ a Borel subgroup such that $T\subset B.$ The flag variety is the variety $\F=G/B$ of translates of the Borel subgroup $G$, and it admits a $T$-action, induced by the left $G$-action. Each Weyl group element determines a $T$-fixed point $e_w$ in the flag variety (by taking a representative of $w$), and these are all the torus-fixed points. Let $B^{-}$ denote the opposite Borel subgroup. The {\em Schubert variety} $X_w$ is the closure of $B^- e_w$ in the flag variety; it has codimension $\ell(w)$ - the length of $w$ in the Weyl group $W$.




In general, if $X$ is a topological space with a left $T$-action, the equivariant cohomology of $X$ is the ordinary cohomology of a ``mixed space'' $(X)_T$, whose definition (see e.g.
\cite{GKM} and references therein) we recall. Let $ET
\longrightarrow BT$ be the universal $T-$bundle. The $T-$action on
$X$ induces an action on the product $ET \times X$ by $t\cdot
(e,x)=(et^{-1},tx)$. The quotient space $(X)_T=(ET \times X)/T$ is
the ``homotopic quotient'' of $X$ and the ($T-$)equivariant
cohomology of $X$ is by definition
\[ H^{i}_T(X)=H^i(X_T). \] In particular, the equivariant
cohomology of a point, denoted by $\mathcal{S}$, is equal to the
ordinary cohomology of the classifying space $BT$. If $\chi$ is a
character in $\hat{T}=Hom(T, \C^*)$ it determines a line
bundle $L_\chi: ET\times_T \C_\chi \longrightarrow BT$ where
$\C_\chi$ is the $1-$dimensional $T-$module determined by $\chi$.
It turns out that the morphism $\hat{T} \longrightarrow
H^2_T(pt)$ taking the character $\chi$ to the first Chern
class $c_1(L_\chi)$ extends to an isomorphism from the symmetric
algebra of $\hat{T}$ to $H^*_T(pt)$. Therefore, if one chooses a basis $t_1, \ldots , t_n$ for $\hat{T}$, then $\mathcal{S}$ is the polynomial ring $\Z[t_1, \ldots , t_n]$.



Returning to the situation when $X=\F$, note that $X_{w}$ is a $T$-stable,
therefore its fundamental class determines
the {\em (equivariant) Schubert class} $\sigma_{w}=[X_{w}]_T$ in
$H_{T}^{2\ell(w)}(\mathcal{F}).$
It is well-known that the Schubert classes form an $H_{T}^{*}(pt)$-basis of $H_{T}^{*}(\mathcal{F}).$

\subsection{Localization map}
Denote by $\F^T= \{e_v| v \in W\}$ the set of $T$-fixed points in $\F$; the inclusion
$\iota: \mathcal{F}^{T}\hookrightarrow \mathcal{F}$
is $T$-equivariant and 
induces 
a homomorphism
$\iota^{*}: H_{T}^{*}(\mathcal{F})\longrightarrow 
H_{T}^{*}(\mathcal{F}^{T})=\prod_{v\in W}H_{T}^{*}(e_{v}).$
We identify each $H_{T}^{*}(e_{v})$ with $\cS$ and for $\eta \in H^*_T (\F)$ we denote its localization in $H_{T}^{*}(e_{v})$ by $\eta|_v$.
Let $R^{+}$ denote the set of positive roots corresponding 
to $B$ and set $R^{-}=-R^{+},
R=R^{+}\cup R^{-}.$
Each root $\alpha$ in $R$ can be 
regarded as a linear 
form in $\cS.$ Let $s_{\alpha}$ denote the reflection
corresponding to the root $\alpha.$ Remarkably, the localization map $\iota^*$ in injective,
and the elements $\eta=(\eta|_{v})_{v}$ in $\prod_{v\in W}\cS$ in the image of $\iota^*$
are characterized by the {\em GKM conditions} (see e.g. \cite{A}):
$$
\eta|_{v}-\eta|_{s_{\alpha}v}\;\mbox{is a multiple
of}\;\alpha
$$
for all $v$ in $W$ and $\alpha\in R^{+},$ where $s_\alpha \in W$ is the reflection associated to $\alpha$.


\subsection{Schubert classes}
We recall a characterization of the Schubert class $\sigma_{w}.$  
Let $\leq$ denote the Bruhat-Chevalley ordering
on $W$; then $e_{v}\in X_{w}$ if and only if
$w\leq v.$
\begin{prop}\cite{A},\cite{KK} \label{xi} 
The Schubert class $\sigma_{w}$ is characterized by the
following conditions:
\begin{enumerate}
\item $\sigma_w|_{v}$ vanishes unless $w\leq v,$
\item If $w\leq v$ then $\sigma_{w}|_{v}$ is homogeneous of degree $\ell(w)$,
\item $\sigma_w|_{w}=\prod_{\alpha\in R^{+}\cap wR^{-}}\alpha.$
\end{enumerate}
\end{prop}
\begin{prop}\label{prop:basis}
Any cohomology class $\eta$ in $H_{T}^{*}(\mathcal{F})$
can be written uniquely as an 
an $H_{T}^{*}(pt)$-linear 
combination of $\sigma_{w}$ using only 
those $w$ such that $w\geq u$ for some 
$u$ with $\eta|_{u}\ne 0.$
\end{prop}
\begin{proof}
The corresponding fact
for the Grassmann variety
is proved in \cite{KnTa}.
The same proof works for the general flag variety also.
\end{proof}

\subsection{Actions of Weyl group}\label{ssec:WeylAction}
There are two actions of the Weyl group on the equivariant cohomology ring $H_{T}^{*}(\F)$, which are used to define corresponding divided-difference operators. 
In this section we will follow the approach presented in \cite{Kn}. Identify $\eta \in H_{T}^{*}(\F)$ with the sequence of polynomials $(\eta|_v)_{v \in W}$ arising from the localization map. For $w \in W$ define
\[ (w^{R}\eta)|_v =  \eta|_{vw} \quad (w^{L}\eta)|_v = w \cdot (\eta|_{w^{-1}v}) \/.\] 
It is proved in \cite{Kn} that these are well defined actions on $\eqcoh(\F_n),$ and that $w^{R}$ is $H^*_{T}(pt)$-linear, while $w^{L}$ it is not (precisely because it acts on the polynomials' coefficients).

\subsection{Divided difference operators}\label{ssec:divdiff}
For each simple root $\alpha_{i}$, we define
the {\it divided difference operators\/} 
$\partial_{i}$ and $\delta_{i}$
on 
$H_{T}^{*}(\mathcal{F})$
by 
$$
(\partial_{i}\eta)|_{v}
=\frac{\eta|_{v}-(s_{i}^{R}\eta)|_{v}}{-v(\alpha_{i})},\quad
(\delta_{i}\eta)_{v}=\frac{\eta|_{v}-(s_{i}^{L}\eta)|_{v}}{\alpha_{i}}\quad (v\in W).
$$
These rational functions are proved to be actually polynomials. They satisfy the GKM conditions, and thus
give elements in $H_{T}^{*}(\mathcal{F})$
(see \cite{Kn}).
We call $\partial_{i}$'s (resp. $\delta_{i}$'s) {\it right\/} 
(resp. {\it left\/}) divided difference operators.
The operator $\partial_{i}$ was introduced 
in \cite{KK}.
On the ordinary cohomology, analogous operators 
to $\partial_{i}$'s are introduced 
independently by Bernstein et al. \cite{BGG}
and Demazure \cite{D}.
The left divided difference operators
$\delta_{i}$ was studied by Knutson in \cite{Kn}
(see also \cite{Ty}). Note that $\partial_{i}$ is $H_{T}^{*}(pt)$-linear 
whereas $\delta_{i}$ is not. The next proposition was stated \cite[Prop.2]{Kn} (see also \cite{KK,Ty}).

\begin{prop}
\label{prop:propertiesdiv}
\begin{enumerate}
\item \label{welldef} Operators $\partial_{i}$ and $\delta_{i}$
are well-defined on the ring $H_{T}^{*}(\mathcal{F})$;
\item The left and right divided difference 
operators commute with each other;
\item \label{eq:partial_sigma} We have
\begin{equation}
\partial_{i}\sigma_{w}=\begin{cases}
\sigma_{ws_{i}} &\mbox{if}\;\ell(ws_{i})=\ell(w)-1\\
0 &\mbox{if}\;\ell(ws_{i})=\ell(w)+1
\end{cases},\quad
\delta_{i}\sigma_{w}=\begin{cases}
\sigma_{s_{i}w} &\mbox{if}\;\ell(s_{i}w)=\ell(w)-1\\
0 &\mbox{if}\;\ell(s_{i}w)=\ell(w)+1
\end{cases}.\label{eq:divdiffsigma}
\end{equation}
\end{enumerate}
\end{prop}
\begin{proof} We only prove (\ref{eq:divdiffsigma})
for $\delta_{i}$ here, as the rest is proved in \cite{Kn}.
By Prop. \ref{xi} 
$(\delta_{i}{\sigma}_{w})|_{v}$ is nonzero
only for $\{v|\;v\geq w\;\mbox{or}\; s_{i}v\geq w\}.$
This implies that
the element $\delta_{i}{\sigma}_{w}$
is a $H_{T}^{*}(pt)$-linear combination of 
$\{{\sigma}_{v}\,|\,v\geq w\;\mbox{or}\; s_{i}v\geq w\}$
by Prop. \ref{prop:basis}.
Moreover $\sigma_{v}$ appearing in the linear combination
have degree at most $\ell(w)-1.$
Thus if $\ell(s_{i}w)=\ell(w)+1$ then  
$\delta_{i}{\sigma}_{w}$ must vanish. 
If $\ell(s_{i}w)=\ell(w)-1$
the only possible
term is a multiple of ${\sigma}_{s_{i}w}.$
In this case we calculate
$$
\left(\delta_{i}{\sigma}_{w}\right)|_{s_{i}w}
=\frac{{\sigma}_{w}|_{s_{i}w}
-s_{i}({\sigma}_{w}|_{w})}{\alpha_{i}}
=-\frac{s_{i}({\sigma}_{w}|_{w})}{\alpha_{i}},
$$
where we used
${\sigma}_{w}|_{s_{i}w}=0$
since $s_{i}w<w.$
Here we recall the following well-known fact that 
$$w>s_{i}w\Longrightarrow 
s_{i}(R^{+}\cap w R^{-})=
(R^{+}\cap s_{i}w R^{-})\sqcup \{-\alpha_{i}\}.$$
So we have $$s_{i}({\sigma}_{w}|_{w})
=\prod_{\beta\in s_{i}(R^{+}\cap w R^{-})}\beta
=(-\alpha_{i})\prod_{\beta\in (R^{+}\cap s_{i}w R^{-})}\beta
=(-\alpha_{i})\cdot{\sigma}_{s_{i}w}|_{s_{i}w}.$$
By the characterization (Prop. \ref{xi}),
we have $\delta_{i}{\sigma}_{w}=\sigma_{s_{i}w}.$
\end{proof} 

\section{Classical groups}
\label{sec:Classical}
\setcounter{equation}{0} 
In this section, we fix the notations 
for the root systems, Weyl groups,
for the classical groups
used throughout the paper. 
\subsection{Root systems}\label{ssection:root_systems} Let $G_n$ be the classical Lie group of one of the types $\mathrm{B}_n,\mathrm{C}_n$ or $\mathrm{D}_n$, i.e. the symplectic group $Sp(2n,\C)$ in type $\mathrm{C}_n$, the odd orthogonal group $SO(2n+1,\C)$ 
in types $\mathrm{B}_n$ and $SO(2n,\C)$ in type $\mathrm{D}_n$. 
Correspondingly we have the set $R_{n}$ of 
roots, and the set  of simple roots. These are subsets of
the character group 
$\hat{T}_{n}=\bigoplus_{i=1}^{n}\Z t_{i}$
of $T_{n},$ the maximal torus of $G_{n}.$

The positive roots $R_{n}^{+}$ 
(set $R_{n}^{-}:=-R_{n}^{+}$ the negative roots) 
are given by 
\begin{eqnarray*}
\mbox{Type}\; \mathrm{B}_{n}:&
R_{n}^{+}&=\{t_{i}\;|\; 1\leq i\leq n\}
\cup\{t_{j}\pm t_{i}\;|\;1\leq i<j\leq n\},\\
\mbox{Type}\; \mathrm{C}_{n}:&
R_{n}^{+}&=\{2t_{i}\;|\; 1\leq i\leq n\}
\cup\{t_{j}\pm t_{i}\;|\;1\leq i<j\leq n\},\\
\mbox{Type}\; \mathrm{D}_{n}:&
R_{n}^{+}&=\{t_{j}\pm t_{i}\;|\;1\leq i<j\leq n\}.
\end{eqnarray*}

The following are the simple roots:
\begin{eqnarray*}
\mbox{Type}\; \mathrm{B}_{n}:&
\alpha_{0}&=t_{1},
\quad \alpha_{i}=t_{i+1}-t_{i}\quad(1\leq i\leq n-1) ,\\
\mbox{Type}\; \mathrm{C}_{n}:&
\alpha_{0}&=2t_{1},
\quad
\alpha_{i}=t_{i+1}-t_{i}\quad(1\leq i\leq n-1),\\
\mbox{Type}\; \mathrm{D}_{n}:&
\alpha_{\hat{1}}&=t_{1}+t_{2},
\quad \alpha_{i}=t_{i+1}-t_{i}\quad(1\leq i\leq n-1).
\end{eqnarray*}

We introduce a symmetric
bilinear form $(\cdot,\cdot)$ on 
$\hat{T}_{n}\otimes_{\Z}\Q$
by $(t_{i},t_{j})=\delta_{i,j}.$
The simple {\it coroots\/} $\alpha_{i}^{\vee}$
are defined to be $\alpha_{i}^{\vee}=2\alpha_{i}/(\alpha_{i},\alpha_{i}).$
Let $\omega_{i}$ denote the {\it fundamental weights}, i.e. those elements in $\hat{T}_{n}\otimes_{\Z}\Q$ such that $(\omega_{i},\alpha_{i}^{\vee})=\delta_{i,j}$.
They are explicitly given as follows:
\begin{eqnarray*}
\mbox{Type}\; \mathrm{B}_{n}:&
\omega_{0}&={\textstyle\frac{1}{2}}\left(t_{1}+t_{2}+\cdots+t_{n}\right),
\quad \omega_{i}=t_{i+1}+\cdots+t_{n}\quad(1\leq i\leq n-1),\\
\mbox{Type}\; \mathrm{C}_{n}:&
\omega_{i}&=t_{i+1}+\cdots+t_{n}\quad(0\leq i\leq n-1),\\
\mbox{Type}\; \mathrm{D}_{n}:&
\omega_{\hat{1}}&={\textstyle\frac{1}{2}}\left(t_{1}+t_{2}+\cdots+t_{n}\right),\quad
\omega_{1}={\textstyle\frac{1}{2}}\left(-t_1+t_2+\cdots+t_n\right),\\
&\omega_{i}&=t_{i+1}+\cdots+t_{n}\quad(2\leq i\leq n-1).
\end{eqnarray*}

\subsection{Weyl groups}\label{ssec:Weyl}
Set $I_{\infty}=\{0,1,2,\ldots\}$ and $I_{\infty}'
=\{\hat{1},1,2,\ldots\}.$
We define the Coxeter
group $(W_{\infty},I_{\infty})$
(resp. $(W_{\infty}',I_{\infty})$)
of infinite rank, and its finite {\it parabolic\/} subgroup
$W_{n}$
(resp. $W'_{n}$)
by the following Coxeter graphs:

\noindent
$\mathrm{C}_{n}\subset \mathrm{C}_{\infty}$ ($\mathrm{B}_{n}\subset \mathrm{B}_{\infty}$)
\setlength{\unitlength}{0.4mm}
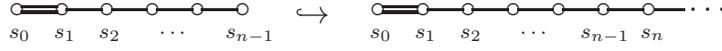
\begin{center}
  \begin{picture}(200,25)
  \thicklines
  \put(0,15){$\circ$}
  \put(4,16.5){\line(1,0){12}}
  \put(4,18.5){\line(1,0){12}}
  \multiput(15,15)(15,0){4}{
  \put(0,0){$\circ$}
  \put(4,2.4){\line(1,0){12}}}
  \put(75,15){$\circ$}
  \put(120,10)
  {\put(0,5){$\circ$}
  \put(4,6.5){\line(1,0){12}}
  \put(4,8.5){\line(1,0){12}}
  \multiput(15,5)(15,0){5}{
  \put(0,0){$\circ$}
  \put(4,2.4){\line(1,0){12}}}
  \put(90,5){$\circ$}}
  \put(0,8){\tiny{$s_0$}}
  \put(15,8){\tiny{$s_1$}}
  \put(30,8){\tiny{$s_2$}}
  \put(72,8){\tiny{$s_{n-1}$}}
  \put(120,8){\tiny{$s_0$}}
  \put(135,8){\tiny{$s_1$}}
  \put(150,8){\tiny{$s_2$}}
  \put(190,8){\tiny{$s_{n-1}$}}
  \put(210,8){\tiny{$s_{n}$}}
  \put(170,8){\tiny{$\cdots$}}
  \put(50,8){\tiny{$\cdots$}}
  \put(95,13){$\hookrightarrow$}
  \put(214.5,17.5){\line(1,0){10}}
  \put(226,15){$\cdots$}
  \end{picture}
\end{center}
$\mathrm{D}_{n}\subset \mathrm{D}_{\infty}$
\setlength{\unitlength}{0.4mm}
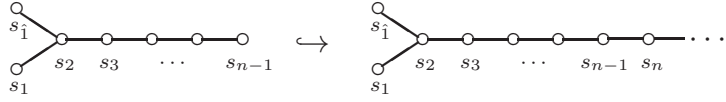
\begin{center}
  \begin{picture}(200,20)
  \thicklines
  \put(0,25){$\circ$}
  \put(0,5){$\circ$}
  \put(4,26){\line(3,-2){12}}
  \put(4,8.5){\line(3,2){12}}
  \multiput(15,15)(15,0){4}{
  \put(0,0){$\circ$}
  \put(4,2.4){\line(1,0){12}}}
  \put(75,15){$\circ$}
  \put(120,0)
  {\put(0,25){$\circ$}
  \put(0,5){$\circ$}
  \put(4,26){\line(3,-2){12}}
  \put(4,8.5){\line(3,2){12}}
  \multiput(15,15)(15,0){4}{
  \put(0,0){$\circ$}
  \put(4,2.4){\line(1,0){12}}}
  \put(75,15){$\circ$}
  \put(120,10)}
  \put(0,20){\tiny{$s_{\hat{1}}$}}
  \put(0,0){\tiny{$s_1$}}
  \put(15,8){\tiny{$s_2$}}
  \put(30,8){\tiny{$s_3$}}
  \put(72,8){\tiny{$s_{n-1}$}}
  \put(120,20){\tiny{$s_{\hat{1}}$}}
  \put(120,0){\tiny{$s_1$}}
  \put(135,8){\tiny{$s_2$}}
  \put(150,8){\tiny{$s_3$}}
  \put(190,8){\tiny{$s_{n-1}$}}
  \put(210,8){\tiny{$s_{n}$}}
  \put(170,8){\tiny{$\cdots$}}
  \put(50,8){\tiny{$\cdots$}}
  \put(95,13){$\hookrightarrow$}
  \put(199,17.5){\line(1,0){12}}
  \put(210,15){$\circ$}
    \put(214.5,17.5){\line(1,0){10}}
  \put(226,15){$\cdots$}
  \end{picture}
\end{center}
More explicitly,
the group $W_{\infty}$ (resp. $W_{\infty}'$ ) is
generated by the simple reflections
$s_{i}\,(i\in I_{\infty})$ (resp. $s_{i}\,(i\in I_{\infty}')$)
subject to the relations: 
\begin{equation}
\begin{cases}s_{i}^{2}=e\;(i\in I_{\infty})\\
s_{0}s_{1}s_{0}s_{1}=s_{1}s_{0}s_{1}s_{0}\\
s_{i}s_{i+1}s_{i}=s_{i+1}s_{i}s_{i+1} 
\; (i\in I_{\infty}\setminus \{0\})\\
s_{i}s_{j}=s_{j}s_{i}\;(|i-j|\geq 2)
\end{cases},\quad
\begin{cases}
s_{i}^{2}=e\;( i\in I_{\infty}')\\
s_{\hat{1}}s_{2}s_{\hat{1}}=s_{2}s_{\hat{1}}s_{2}\\
s_{i}s_{i+1}s_{i}=s_{i+1}s_{i}s_{i+1} 
\; (i\in I_{\infty}'\setminus \{\hat{1}\})\\
s_{\hat{1}}s_{i}=s_{i}s_{\hat{1}}\;(i\in I_{\infty}',\;i\ne 2)\\
s_{i}s_{j}=s_{j}s_{i}\;(\;i,j\in I_{\infty}'\setminus \{\hat{1}\},
\,|i-j|\geq 2)
\end{cases}.\label{eq:Wrel}
\end{equation}

For general facts on Coxeter groups, we refer to
\cite{BjBr}.
Let $\leq $ denote  the {\it Bruhat-Chevalley order\/} 
on $W_{\infty}$ or $W_{\infty}'.$ The {\it length\/} $\ell(w)$ of $w\in W_{\infty}$ (resp. $w\in W'_{\infty}$) is defined to be the least number $k$ of simple reflections
in any reduced expression of $w\in W_{\infty}$.


The subgroups $W_{n}
\subset W_{\infty}$, 
$W_{n}'\subset W_{\infty}'$ are 
the Weyl groups of the following types:
$$
\mbox{Type}\; B_{n},C_{n}:
W_{n}=\langle s_{0},s_{1},s_{2},\ldots,s_{n-1}
\rangle,\quad
\mbox{Type}\; D_{n}:
W_{n}'=\langle s_{\hat{1}},s_{1},s_{2},\ldots,s_{n-1}\rangle.
$$ 
It is known that the inclusion $W_{n}\subset W_{\infty}$
(resp. $W_{n}'\subset W_{\infty}'$) preserves the length and the Bruhat-Chevalley order, while $W_{\infty}'\subset W_{\infty},$
(resp. $W_{n}'\subset W_{n}$) is not (using terminology from \cite{BjBr} this says that $W_n$ is a {\em parabolic subgroup} of $W_\infty$, while $W'_\infty$ is not). From now on, whenever possible, we will employ the notation explained in \S \ref{ssec:infinite_weyl}, and use bold fonts $\W_\infty$ respectively $\W_n$ to make uniform statements.



\subsection{Signed permutations}\label{sec:signed.permut}
The group $W_\infty$ 
is identified with the set of all permutations $w$ of the set $\{1,2,\ldots\}\cup \{\bar{1},\bar{2},\ldots\}$ such that $w(i)\ne i$ for only finitely many $i$, and $\overline{w(i)}=w(\bar{i})$ for all $i$. These can also be considered as 
signed (or barred) permutation of $\{1,2,\ldots\}$; we often use one-line notation 
$w=(w(1),w(2),\ldots)$ to denote an element $w\in W_{\infty}$. The simple reflections are
identified with the transpositions
$s_{0}=(1,\bar{1})$ and $s_{i}=(i+1,i)(\overline{i},\overline{i+1})$ for $i\geq 1.$
The subgroup $W_n\subset W_{\infty}$
is described as 
$$
W_{n}=\{w\in W_{\infty}\;|\; w(i)=i\;\mbox{for} \; i>n\}.
$$
In one-line notation,
we often denote 
an element $w\in W_{n}\subset W_{\infty}$
by the finite sequence 
$(w(1),\ldots,w(n)).$

The group $W_{\infty}'$, as a (signed) permutation group,
can be realized
as the subgroup of $W_{\infty}$
consisting of elements in $W_\infty$ with even number of sign changes.
The simple reflection $s_{\hat{1}}$ is identified with $s_{0}s_{1}s_{0}
\in W_{\infty}$, so as a permutation
$s_{\hat{1}}=(1,\bar{2})(2,\bar{1}).$

\subsection{Grassmannian elements} 
An element $w\in W_{\infty}$ is 
a {\it Grassmannian element\/} 
if $$w(1)<w(2)<\cdots<w(i)<\cdots$$ in the ordering
$ \cdots<\bar{3}<\bar{2}<\bar{1}<1<2<3<\cdots
.$ 
Let $W_{\infty}^{0}$ denote the set of all Grassmannian elements in $W_{\infty}.$
For $w\in W_{\infty}^{0}$, let $r$ be the number such that 
\begin{equation} w(1)<\cdots<w(r)<1\quad\mbox{and}
\quad \bar{1}<w(r+1)<w(r+2)<\cdots.
\label{eq:GrassPerm}
\end{equation} 
Then we define the $r$-tuple of positive integers  
$\lambda=(\lambda_{1},\ldots,\lambda_{r})$ by
$ \lambda_{i}=\overline{w(i)}$
for $1\leq i\leq r.$
This is a {\it strict partition\/} 
i.e. a partition with distinct parts:
$\lambda_{1}>\cdots>\lambda_{r}>0.$ 
Let $\mathcal{SP}$ denote the set of all
strict partitions. 
The correspondence gives a bijection
$
W_{\infty}^{0}\longrightarrow
\mathcal{SP}.
$
We denote by $w_{\lambda}\in W_{\infty}^{0}$
the Grassmannian element corresponding
to $\lambda\in \mathcal{SP}$; then 
$\ell(w_{\lambda})=|\lambda|=\sum_{i}\lambda_{i}$.
Note that this bijection preserves the partial order 
when $\mathcal{SP}$ is considered to be a partially ordered set given by 
the inclusion $\lambda\subset \mu$
of strict partitions.

We denote by $W_{\infty}^{\hat{1}}$ 
the set of all
Grassmannian elements contained in $W_{\infty}'.$ 
For $w\in W_{\infty}^{\hat{1}},$
the number $r$ in (\ref{eq:GrassPerm}) is always even.
Define the strict partition 
$\lambda'=(\lambda'_{1}>\cdots>\lambda'_{r}\geq 0)$ by setting $ \lambda'_{i}=\overline{w(i)}-1$
for $1\leq i\leq r.$
Note that $\lambda'_{r}$ can be zero
this time.
This correspondence gives also a
bijection $
W_{\infty}^{\hat{1}}
\longrightarrow 
\mathcal{SP}.$
We denote by $w'_{\lambda}\in W_{\infty}^{\hat{1}}$
the element corresponding
to $\lambda\in \mathcal{SP}.$
As before, $\ell(w'_{\lambda})=|\lambda|$ where 
$\ell(w)$ denotes the length of $w$ in $W_\infty'$.

\bigskip

{\bf Example.} Let $\lambda=(4,2,1).$
Then the corresponding Grassmannian elements are given by
$w_{\lambda}=\bar{4}\bar{2}\bar{1}3
=s_{0}s_{1}s_{0}s_{3}s_{2}s_{1}s_{0}$
and $w'_{\lambda}=\bar{5}\bar{3}\bar{2}\bar{1}4
=s_{\hat{1}}s_{2}s_{1}s_{4}s_{3}s_{2}s_{\hat{1}}.$

\bigskip

The group $W_{\infty}$ (resp. $W_{\infty}'$) 
has a parabolic subgroup 
generated by $s_{i}\;(i\in I_{\infty}\setminus \{0\})$
(resp. $s_{i}\,(i\in I'_{\infty}\setminus \{\hat{1}\}).$
We denote these subgroups by $S_{\infty}
=\langle s_{1},s_{2},\ldots\rangle$
since it is isomorphic to the infinite Weyl group of type $\mathrm{A}.$
The product map
$$
W_{\infty}^{0}\times S_{\infty}
\longrightarrow W_{\infty}
\quad(\mbox{resp.}\;
W_{\infty}^{\hat{1}}\times S_{\infty}
\longrightarrow W'_{\infty}),
$$
given by $(u,w)\mapsto uw$
is a bijection satisfying $\ell(uw)=\ell(u)+\ell(w)$
(cf. \cite[Prop. 2.4.4]{BjBr}).
As a consequence, $w_{\lambda}$ (resp. $w_{\lambda}'$) is the unique element
of minimal length
in the left coset $w_{\lambda}S_{\infty}$
(resp. $w_{\lambda}'S_{\infty}$).

\section{Schur's $Q$-functions and its factorial
analogues}\label{sec:Schur}
\setcounter{equation}{0}
\subsection{Schur's $Q$-functions}
Our main reference for symmetric functions is 
\cite{Mac}. 
Let $x=(x_{1},x_{2},\ldots)$
be infinitely many
indeterminates.
Define $Q_i(x)$ as the coefficient of $u^i$ in
the generating function
$$f(u)=
\prod_{i=1}^{\infty}
\frac{1+x_{i}u}{1-x_{i}u}
=\sum_{k\geq 0}Q_{k}(x)u^{k}.
$$
Note that $Q_{0}=1.$
Define $\Gamma$ to be $\Z[Q_1(x),Q_2(x),\ldots].$
The identity
$f(u)f(-u)=1$ yields
\begin{equation}\label{eq:quadraticQ}
Q_{i}(x)^{2}+2\sum_{j=1}^{i}(-1)^{j}Q_{i+j}(x)Q_{i-j}(x)=0
\quad \mbox{for}\quad i\geq 1.
\end{equation}
It is known that the ideal of relations among the functions $Q_k(x)$ is 
 generated by the previous relations.  
For $i\geq j\geq 0,$ define elements
$$
Q_{i,j}(x):=Q_{i}(x)Q_{j}(x)+2\sum_{k=1}^{j}(-1)^{k}
Q_{i+k}(x)Q_{j-k}(x).
$$
 Note that $Q_{i,0}(x)=Q_{i}(x)$
and $Q_{i,i}(x)\,(i\geq 1)$ is identically zero. 
For $\lambda$ a strict partition
we write $\lambda=(\lambda_{1}>\lambda_{2}>\cdots
>\lambda_{r}\geq 0)$ with $r$ even.
Then the corresponding 
Schur's $Q$-function
$Q_{\lambda}=Q_{\lambda}(x)$ is defined by
$$
Q_{\lambda}(x)=\mathrm{Pf}(Q_{\lambda_{i},\lambda_{j}}(x))_{1\leq i<j\leq r},
$$
where $\mathrm{Pf}$ denotes the Pfaffian.
It is known then that the functions
$Q_\lambda(x)$ for $\lambda\in \mathcal{SP}$ form
a $\Z$-basis of $\Gamma.$
The {\em $P$-Schur function} is defined to be $P_{\lambda}(x)=2^{-\ell(\lambda)}Q_{\lambda}(x)$ where $\ell(\lambda)$ is the number of non-zero parts in $\lambda.$ 
The next lemma shows that the $Q$-Schur function is {\em supersymmetric\/.}
\begin{lem}\label{lem:super}
Each element 
$\varphi(x)$ in $\Gamma$
satisfies 
$$
\varphi(t,-t,x_{1},x_{2},\ldots)
=\varphi(x_{1},x_{2},\ldots)
$$
where $t$ is an indeterminate.
\end{lem}
\begin{proof}
It suffices to show this for the ring generators
$Q_{i}(x).$ This follows immediately from the generating function. 
\end{proof}

\subsection{Factorial $Q$ and $P$-Schur functions}\label{ssec:FacSchur} In this section we recall the definition and some properties of the factorial $Q$-Schur and $P$-Schur functions defined by V.N. Ivanov in \cite{Iv}. Fix $n \ge 1$ an integer, $\lambda$ a strict partition of length $r \le n$ and $a=(a_i)_{i \ge 1}$ an infinite sequence. By $(x|a)^k$ we denote the product $(x-a_1)\cdots (x-a_k)$. According to \cite[Def. 2.10]{Iv} the factorial $P$-Schur function $P_\lambda^{(n)}(x_1, \ldots x_n|a)$ is defined by: \[ P_\lambda^{(n)}(x_1, \ldots x_n|a) = \frac{1}{(n-r)!} \sum_{w \in S_n} w\cdot \bigl( \prod_{i=1}^r (x_i|a)^{\lambda_i} \prod_{i \le r, i < j \le n} \frac{x_i+x_j}{x_i-x_j} \bigr) \/, \] where $w$ acts on variables $x_i$. If $a_1=0$ this function is stable, i.e. $P_\lambda^{(n+1)}(x_1, \ldots x_n,0|a) = P_\lambda^{(n)}(x_1, \ldots x_n|a)$, therefore there is a well defined limit denoted $P_\lambda(x|a)$. It was proved in \cite[Prop. 8]{IN} that if $a_1 \neq 0$, $P_\lambda^{(n)}(x_1, \ldots x_n|a)$ is stable modulo $2$, i.e. $P_\lambda^{(n+2)}(x_1, \ldots x_n,0,0|a)=P_\lambda^{(n)}(x_1, \ldots x_n|a)$; so in this case there is a well-defined even and odd limit. From now on we will denote by $P_\lambda(x|a)$ the {\em even} limit of these functions. Define also the factorial $Q$-Schur function $Q_\lambda(x|a)$ to be \[ Q_\lambda(x|a) = 2^{\ell(\lambda)}P_\lambda(x|0,a) \/,\] where $\ell(\lambda)$ is the number of non-zero parts of $\lambda$. As explained in \cite{IN}, the situation $a_1 \neq 0$ is needed to study type $\mathrm{D}$; in types $\mathrm{B,C}$, the case $a_1=0$ will suffice. For simplicity, we will also denote $P_\lambda^{(n)}(x_1, \ldots x_n|a)$ by $P_\lambda(x_1, \ldots x_n|a)$. 
Let now $t=(t_{1},t_{2},\ldots)$ be indeterminates. Define:
$$
t_{\lambda}=
(t_{\lambda_{1}},\ldots,
t_{\lambda_{r}},0,0,\ldots),
\quad
t_{\lambda}'=\begin{cases}
(t_{\lambda_{1}+1},\ldots,
t_{\lambda_{r}+1},0,0,\ldots) &\mbox{if}\;r\;\mbox{is even}\\
(t_{\lambda_{1}+1},\ldots,
t_{\lambda_{r}+1},t_{1},0,\ldots)&\mbox{if}\;r\;\mbox{is odd}
\end{cases}.
$$
Let also $w_{\lambda}\in W_{\infty}^{0}$ (resp. $w_{\lambda}'\in W_{\infty}^{\hat{1}}$ ) be the Grassmann element corresponding 
to $\lambda\in\mathcal{SP}.$ 
We associate to $\lambda$ its
{\it shifted Young diagram} $Y_\lambda$
as the set of {\it boxes}
with coordinates $(i,j)$
such that $1\leq i\leq r$ and $i\leq j\leq i+\lambda_i-1.$
We set $\lambda_j=0$ for $j>r$ by convention.
Define
$$
H_{\lambda}(t)=\prod_{(i,j)\in Y_{\lambda}}
(t_{\overline{w_{\lambda}(i)}}
+t_{\overline{w_{\lambda}(j)}}),\quad
H_{\lambda}'(t)=\prod_{(i,j)\in Y_{\lambda}}
(t_{\overline{w_{\lambda}'(i)}}
+t_{\overline{w_{\lambda}'(j+1)}}).
$$

{\bf Example.} Let $\lambda=(3,1).$ 
Then $w_{\lambda}=\bar{3}\bar{1}2,$  
$w_{\lambda}'=\bar{4}\bar{2}13,$ and
$$H_{\lambda}(t)=4t_{1}t_{3}(t_{3}+t_{1})(t_{3}-t_{2}),\quad
H_{\lambda}'(t)=(t_{4}+t_{2})(t_{4}-t_{1})(t_{4}-t_{3})(t_{2}-t_{1}).$$
\setlength{\unitlength}{0.6mm}
\begin{center}
\begin{picture}(200,45)
    \put(-25,30){Type $\mathrm{C}$}
  \put(5,35){\line(1,0){60}}
  \put(5,25){\line(1,0){60}}
  \put(25,15){\line(1,0){20}}
 \put(5,25){\line(0,1){10}}
 \put(25,15){\line(0,1){20}}
 \put(45,15){\line(0,1){20}}
 \put(65,25){\line(0,1){10}}
 \put(11,28){\small{$2t_{3}$}}
 \put(31,18){\small{$2t_{1}$}}
 \put(27,28){\small{$t_{3}\!+\!t_{1}$}}
 \put(47,28){\small{$t_{3}\!-\!t_{2}$}}
     \put(75,30){Type $D$}
  \put(105,35){\line(1,0){60}}
  \put(105,25){\line(1,0){60}}
  \put(125,15){\line(1,0){20}}
 \put(105,25){\line(0,1){10}}
 \put(125,15){\line(0,1){20}}
 \put(145,15){\line(0,1){20}}
 \put(165,25){\line(0,1){10}}
 \put(107,28){\small{$t_{4}\!+\!t_{2}$}}
 \put(127,18){\small{$t_{2}\!-\!t_{1}$}}
 \put(127,28){\small{$t_{4}\!-\!t_{1}$}}
 \put(147,28){\small{$t_{4}\!-\!t_{3}$}}
  \put(0,28){\small{$3$}}
  \put(0,18){\small{$1$}}
  \put(13,38){\small{$3$}}
  \put(33,38){\small{$1$}}
  \put(53,38){\small{$\bar{2}$}}
  \put(100,28){\small{$4$}}
  \put(100,18){\small{$2$}}
  \put(113,38){\small{${2}$}}
  \put(133,38){\small{$\bar{1}$}}
  \put(153,38){\small{$\bar{3}$}}
  \put(5,5){$w_{\lambda}=s_{0}s_{2}s_{1}s_{0}=\bar{3}\bar{1}2$}
  \put(105,5){$w_{\lambda}'=s_{1}s_{3}s_{2}s_{\hat{1}}=\bar{4}\bar{2}13$}
  \end{picture}
  \end{center}


\begin{prop} (\cite{Iv})
For any strict partition $\lambda$,
the factorial $Q$-Schur function $Q_{\lambda}(x|t)$
(resp. $P_{\lambda}(x|t)$)
satisfies the following properties:
\begin{enumerate}
\item $Q_{\lambda}(x|t)$ (resp.  $P_{\lambda}(x|t)$) is
homogeneous of degree $|\lambda|=\sum_{i=1}^{r}\lambda_{i},$
\item $Q_{\lambda}(x|t)=Q_{\lambda}(x)+\mbox{lower order terms in}\;x$\\
(resp. $P_{\lambda}(x|t)=P_{\lambda}(x)+\mbox{lower order terms in}\;x$),
\item $Q_{\lambda}(t_{\mu}|t)=0$ (resp. $P_{\lambda}(t_{\mu}'|t)=0$) unless $\lambda\subset \mu,$
\item $Q_{\lambda}(t_{\lambda}|t)=
H_{\lambda}(t)$
(resp. $P_{\lambda}(t'_{\lambda}|t)=
H'_{\lambda}(t)$).
\end{enumerate}
Moreover $Q_{\lambda}(x|t)$ 
(resp. $P_{\lambda}(x|t)$)
($\lambda\in \mathcal{SP}$)
form a $\Z[t]$-basis of $\Z[t]\otimes_{\Z}
\Gamma$ (resp.
$\Z[t]\otimes_{\Z}
\Gamma'$).
\end{prop}
\begin{proof} In the case $t_1=0$ this was proved in \cite[Thm. 5.6]{Iv}. If $t_1 \neq 0$, the identity (3) follows from definition (cf. \cite[Prop.~9]{IN}), while (4) follows from a standard computation.\end{proof}


%
\begin{remark}{\rm The statement in the previous proposition can be strengthened by showing that the properties (1)-(4) characterize the factorial $Q-$Schur (respectively $P$-Schur) functions. For $t_1=0$ this was shown in \cite[Thm.~5.6]{Iv}. A similar proof can be given for $t_1 \neq 0$, but it also follows from Thm. \ref{thm:isom} below. The characterization statement will not be used in this paper.} \end{remark}
\begin{remark}{\rm The function
$Q_{\lambda}(x|t)$ belongs actually to
$\Gamma\otimes_{\Z}\Z[t_{1},t_{2},\ldots,t_{\lambda_{1}-1}]$ and $P_{\lambda}(x|t)$ to
$\Gamma'\otimes_{\Z}\Z[t_{1},t_{2},\ldots,t_{\lambda_{1}}]$. For example we have
$$
Q_{i}(x|t)=\sum_{j=0}^{i-1}(-1)^{j}
e_{j}(t_{1},\ldots,t_{i-1})Q_{i-j}(x),\quad
P_{i}(x|t)=\sum_{j=0}^{i-1}(-1)^{j}
e_{j}(t_{1},\ldots,t_{i})P_{i-j}(x).
$$}
\end{remark}

\begin{remark}
{\rm
An alternative formula for $Q_\lambda(x|t)$, in terms of a {\em multi-Pfaffian}, will be given below in \S \ref{sec:Kazarian}}.
\end{remark}
The following proposition will only be used within the proof of the formula for the Schubert polynomial for the longest element in each type, presented in \S \ref{sec:Long} below.  
\begin{prop} [\cite{Iv}]\label{prop:Pf}
Let $\lambda=(\lambda_{1}>\cdots>\lambda_{r}\geq 0)$ a strict partition
with $r$ even. Then
$$
Q_{\lambda}(x|t)=\mathrm{Pf}\left(
Q_{\lambda_{i},\lambda_{j}}(x|t)
\right)_{1\leq i<j\leq r},\quad
P_{\lambda}(x|t)=\mathrm{Pf}\left(
P_{\lambda_{i},\lambda_{j}}(x|t)
\right)_{1\leq i<j\leq r}.
$$ 

\end{prop}
\begin{proof} Again, for $t_1=0$, this was proved in \cite[Thm.3.2]{Iv}, using the approach described in \cite[III.8 Ex.13]{Mac}. The same approach works in general, but for completeness we briefly sketch an argument. Lemma \ref{injective} below shows that there is an injective universal localization map $\Phi: \Z[z] \otimes \Z[t] \otimes \Gamma' \to \prod_{w \in W_\infty'} \Z[t]$. The image of $P_\lambda(x|t)$ is completely determined by the images at Grassmannian Weyl group elements $w_\mu'$ and it is given by $P_\lambda(t'_\mu|t)$. But by the results from \cite[\S 10]{IN} we have that $P_\lambda(t'_\mu|t)= \mathrm{Pf}( P_{\lambda_i, \lambda_j}(t'_\mu|t))_{1 \le i< j \le r}$. The result follows by injectivity of $\Phi$.\end{proof}


We record here the following formula used later. The proof is by a standard computation (see e.g. the proof of \cite[Thm. 8.4]{Iv}).
\begin{lem}
We have
\begin{equation}
P_{k,1}(x|t)=
P_{k}(x|t)P_{1}(x|t)
-P_{k+1}(x|t)-(t_{k+1}+t_{1})P_{k}(x|t)\quad
\mbox{for}\quad k\geq 1.\label{eq:P2row}
\end{equation}
\end{lem}

\subsection{Factorization formulae}
In this section we present several
{\em factorization} formulas for the factorial $P$ and $Q$-Schur functions, 
which will be used later in \S \ref{sec:Long}. To this end, we first consider the case of ordinary factorial Schur functions.
\subsubsection{Factorial Schur polynomials}
Let $\lambda=(\lambda_{1},\ldots,\lambda_{n})$ 
be a partition.
Define the {\it factorial
Schur polynomial\/} by 
$$s_{\lambda}(x_{1},\ldots,x_{n}|t)
=\frac{\det((x_{j}|t)^{\lambda_{i}+n-i})_{1\leq i,j\leq n}}{\prod_{1\leq i<j\leq n}(x_{i}-x_{j})} \/, $$
where $(x|t)^{k}$ denotes $\prod_{i=1}^{k}(x-t_{i})$.
It turns out that $s_{\lambda}(x|t)$ is an element
in 
$\Z[x_{1},\ldots,x_{n}]^{S_{n}}
\otimes \Z[t_{1},\ldots,t_{\lambda_{1}+n-1}].$ 
For some basic properties of these polynomials,
the reader can consult \cite{MS}. 
The following formula will be used to prove
Lem. \ref{lem:piDelta} in \S \ref{sec:Long}. 
\begin{lem}\label{lem:A-long} We have
$
s_{\rho_{n-1}}(t_{1},\ldots,t_{n}|t_{1},-z_{1},t_{2},-z_{2},
\ldots,
t_{n-1},-z_{n-1})=\prod_{1\leq i<j\leq n}(t_{j}+z_{i}).
$
\end{lem}
\begin{proof} When variables $z_i,t_i$ are specialized as in this Lemma,
the numerator is 
an anti-diagonal lower triangular matrix. The entry on the $i$-th row 
on the anti-diagonal is given by 
$\prod_{j=1}^{i-1}(t_{i}-t_{j})(t_{i}+z_{j}).$
The Lemma follows immediately from this.
\end{proof}

Next formula is a version of Lem. \ref{lem:A-long}
which will be used in the proof of Lem. \ref{lem:piDeltaD}. 
\begin{lem}\label{lem:A-longOdd}
If $n$ is odd
then we have
$$
s_{\rho_{n-1}+1^{n-1}}(t_{2},\ldots,t_{n}|t_{1},-z_{1},\ldots,
t_{n-1},-z_{n-1})
=\prod_{j=2}^{n}(t_{j}-t_{1})
\prod_{1\leq i<j\leq n}(t_{j}+z_{i}).
$$
\end{lem}
\begin{proof}
Similar to the proof of Lem. \ref{lem:A-long}.
\end{proof}

\subsubsection{$P$ and $Q$-Schur functions}
We need the following factorization formula.
\begin{lem}\label{lem:factorization} Let $\lambda=(\lambda_{1},\ldots,\lambda_{n})$ be a partition.
Then 
we have
\begin{eqnarray*}
Q_{\rho_{n}+\lambda}(x_{1},\ldots,x_{n}|t)
&=&\prod_{i=1}^{n}2x_{i}\prod_{1\leq i< j\leq n}
(x_{i}+x_{j})\times 
s_{\lambda}(x_{1},\ldots,x_{n}|t),
\end{eqnarray*}
\end{lem}
\begin{proof}
By their very definition
$$Q_{\rho_{n}+\lambda}(x_{1},\ldots,x_{n}|t)
=2^{n}\sum_{w\in S_{n}}w\left[\prod_{i=1}^{n}
x_{i}(x_{i}|t)^{\lambda_{i}+n-i}
\prod_{1\leq i<j\leq n}\frac{x_{i}+x_{j}}{x_{i}-x_{j}}
\right],
$$
where $w$ acts as permutation of the variables
$x_{1},\ldots,x_{n}.$
Since the polynomial
$\prod_{i=1}^{n}x_{i}
\prod_{1\leq i<j\leq n}(x_{i}+x_{j})$ in the parenthesis 
is symmetric in $x$, the last expression factorizes into
$$
2^{n}\prod_{i=1}^{n}x_{i}
\prod_{1\leq i<j\leq n}(x_{i}+x_{j})\times
\sum_{w\in S_{n}}w\left[\prod_{i=1}^{n}
(x_{i}|t)^{\lambda_{i}+n-i}
\prod_{1\leq i<j\leq n}({x_{i}-x_{j}})^{-1}
\right]$$
Then by the definition of $s_{\lambda}(x_{1},\ldots,x_{n}|t)$
we have the lemma.
\end{proof}

The following two lemmas are proved in the same way: 
\begin{lem}\label{lem:factorizationD} Assume $n$ is even.
Let 
$\lambda=(\lambda_{1},\ldots,\lambda_{n})$
be a partition.
Then we have
$$
P_{\rho_{n-1}+\lambda}
(x_{1},\ldots,x_{n}|t)
=
\displaystyle{
\prod_{1\leq i<j\leq n}}(x_{i}+x_{j})
\cdot s_{\lambda}(x_{1},\ldots,x_{n}|t).
$$
\end{lem}

\begin{lem}\label{lem:factorizationDodd}
Assume $n$ is odd. Let 
$\lambda=(\lambda_{1},\ldots,\lambda_{n-1})$
be a partition.
Then we have
$$
P_{\rho_{n-1}+\lambda}(x_{1},\ldots,x_{n-1}|t)
=
\prod_{1\leq i<j\leq n-1}(x_{i}+x_{j})
\times s_{\lambda+1^{n-1}}
(x_{1},\ldots,x_{n-1}|t).
$$
\end{lem}
%
\section{Stable Schubert classes}\label{sec:StableSchubert}
\setcounter{equation}{0} 

The aim of this section is to introduce 
{\it stable\/} Schubert classes 
indexed by the Weyl group of 
infinite rank $\W_{\infty}$. Recall that the embeddings of Dynkin diagrams shown in \S \ref{ssec:Weyl} induce embeddings $i:\W_n \to \W_{n+1}$; then $\W_\infty = \bigcup_{n \ge 1} \W_n$.

\subsection{Stable Schubert classes}\label{ssec:SSC}
Let us denote by $\sigma_{w}^{(n)}$ the equivariant Schubert 
class on $\mathcal{F}_{n}$ labeled by $w\in \W_{n}.$
\begin{prop}\label{prop:StabSch} 
The localization of Schubert classes is stable, i.e
$$\sigma_{i(w)}^{(n+1)}|_{i(v)}=\sigma_w^{(n)}|_v
\quad\mbox{for all}\quad
w,v\in \W_{n}.$$
\end{prop}
\begin{proof}
First we claim that 
$\sigma_{i(w)}^{(n+1)}|_{i(v)}\in \Z[t]^{(n)}$
for any $w,v\in \W_{n}.$
Let $w_{0}^{(n)}$ be
the longest element in $\W_{n}$. 
By Prop. \ref{xi}, we have for $v\in \W_{n}$
$$
\sigma_{i(w_{0}^{(n)})}^{(n+1)}|_{i(v)}
=\begin{cases}
\prod_{\beta\in R^{+}_{n}}\beta
&\mbox{if} \;v= w_{0}^{(n)}
\\
0 &\mbox{if} \;v\ne w_{0}^{(n)}
\end{cases}.
$$
In particular these polynomials belong to
$\Z[t]^{(n)}.$
For arbitrary $w\in \W_{n}$,
 any reduced expression of $i(w)$ contains only simple reflections $s_0, \dots , s_{n-1}$ (in type $D$, $s_0$ is replaced by $s_{\hat{1}}$). 
Hence we obtain 
the Schubert
class $\sigma_{i(w)}^{(n+1)}$ 
by applying the divided difference
operators $\partial_{0},\ldots,\partial_{n-1}$
(in type $D$, $\partial_0$ is replaced by $\partial_{\hat{1}}$)
successively to the class $\sigma_{i(w_{0}^{(n)})}^{(n+1)}.$
In this process only the variables $t_{1},\ldots,t_{n}$
are involved to compute 
$\sigma_{i(w)}^{(n+1)}|_{i(v)}\;(v\in \W_{n}).$ 
Hence the claim is proved.

For $w\in \W_{n},$
we consider the element $\eta_{w}$
in $\prod_{v\in \W_{n}}\Z[t]^{(n)}$
given by 
$\eta_{w}|_{v}=\sigma_{i(w)}^{(n+1)}|_{i(v)}\;(v\in \W_{n}).$
We will show 
the element $\eta_{w}$ 
satisfies the conditions in Prop. \ref{xi} that
characterize $\sigma_{w}^{(n)}.$ 
In fact, the vanishing condition 
holds since $i(w)\leq i(v)$
if and only if $w\leq  v.$ 
Homogeneity and the degree condition is
satisfied because $\ell(i(w))=\ell(w).$
The normalization follows from the fact
$R_{n}^{+}\cap w R^{-}_{n}=
R_{n+1}^{+}\cap i(w) R^{-}_{n+1}.$
Thus we have $\eta_{w}=\sigma_{w}^{(n)}$
and the proposition is proved.
\end{proof}

Fix $w$ be in $\W_{\infty}.$
Then, by the previous proposition, for any $v\in \W_{\infty}$, and for any sufficiently
large $n$ such that $w,v\in \W_{n}$,
the polynomial
$\sigma_{w}^{(n)}|_{v}$ does not depend on 
the choice of $n.$
Thus we can introduce
 a unique element
$\sigma_{w}^{(\infty)}
=(\sigma_{w}^{(\infty)}|_{v})_{v}$ in $\prod_{v\in \W_{\infty}}\Z[t]$ such that
$$
\sigma_{w}^{(\infty)}|_{v}
=\sigma_{w}^{(n)}|_{v}$$
for all sufficiently large $n.$ 
We call this element the {\it stable\/} Schubert class.
\begin{Def} Let  $H_{\infty}$ be the $\Z[t]$-submodule of $\prod_{v\in \W_{\infty}}\Z[t]$ spanned by the stable Schubert classes $\sigma_{w}^{(\infty)},w\in \W_{\infty}$, where the $\Z[t]$-module structure is given by the diagonal multiplication.
\end{Def}
We will show later in Cor. \ref{cor:Hsubalg} that 
$H_{\infty}$ is actually a $\Z[t]$-subalgebra 
in the product ring $\prod_{v\in \W_{\infty}}\Z[t].$
The properties of the (finite-dimensional) Schubert classes extend immediately to the stable case. For example, the classes $\sigma_{w}^{(\infty)}\;(w\in W_{\infty})$ are
linearly independent over $\Z[t]$
(Prop. \ref{prop:basis}), and they satisfy the properties from Prop. \ref{xi}. To state the latter, define $R^{+}=\bigcup_{n\geq 1}R_{n}^{+}$, regarded 
as a subset of $\Z[t].$ Then:
\begin{prop}\label{StableClass}
The stable Schubert class satisfies the following:
\begin{enumerate} 
\item[(1)](Homogeneity) $\sigma_w^{(\infty)}|_v$ is homogeneous of degree $\ell(w)$ for each $v\geq w,$ 
\item[(2)](Normalization) \label{char:normal}$\sigma_w^{(\infty)}|_w=\prod_{\beta\in R^{+}\cap w(R^{-})}\beta,$
\item[(3)](Vanishing) $\sigma_w^{(\infty)}|_{v}$ vanishes unless $v\geq w.$
\end{enumerate}
\end{prop}

It is natural to
consider the following stable version of the GKM conditions
in the ring $\prod_{v\in\W_{\infty}}\Z[t]$:
$$\eta|_{v}-\eta|_{s_{\alpha}v}\;\;\mbox{is a multiple
of}\;\alpha\;\mbox{for all}\;
\alpha\in R^{+},\; v\in \W_{\infty}.
$$
Then the stable Schubert class $\sigma_{w}^{(\infty)}$
is the unique element in 
$\prod_{v\in W_{\infty}}\Z[t]$ that 
satisfies the GKM conditions and the
three conditions Prop. \ref{StableClass}. It follows that all the elements
from $H_\infty$ satisfy the GKM conditions. In particular, the proofs from \cite{Kn} can be retraced, and one  
can define the left and right actions of 
$\W_{\infty}$ on $H_{\infty}$
by the same formulas as in \S \ref{ssec:WeylAction}
but for $i\in \I_{\infty}.$
Using these actions, we define also the divided difference operators
$\partial_{i},\delta_{i}$ on $H_{\infty}$
(see \S \ref{ssec:divdiff}). The next result follows again from the finite dimensional case (Prop. \ref{prop:propertiesdiv}).
\begin{prop}\label{prop:divdiff} 
We have
$$
\partial_{i}\sigma_{w}^{(\infty)}=\begin{cases}
\sigma_{ws_{i}}^{(\infty)}& \ell(ws_{i})=\ell(w)-1\\
0 &\ell(ws_{i})=\ell(w)+1
\end{cases},\quad
\delta_{i}\sigma_{w}^{(\infty)}=\begin{cases}
\sigma_{s_{i}w}^{(\infty)}& \ell(s_{i}w)=\ell(w)-1\\
0 &\ell(s_{i}w)=\ell(w)+1
\end{cases}.
$$
\end{prop}

\subsection{Inverse limit of cohomology groups}\label{ssec:invlim}
Let $H_{n}$ denote the image of 
the localization map
$$
\iota_{n}^{*}: H_{T_{n}}^{*}(\mathcal{F}_{n})
\longrightarrow H_{T_{n}}^{*}(\mathcal{F}_{n}^{T_{n}})
=\prod_{v\in \W_{n}}\Z[t]^{(n)}.
$$
By the stability property for the localization of Schubert classes, the natural projections $H_\infty \to H_n \simeq \eqcoh (\F_n)$ are compatible with the homomorphisms $H^*_{T_{n+1}}(\F_{n+1}) \to \eqcoh(\F_n)$ induced by the equivariant embeddings $\F_n \to \F_{n+1}$. Therefore there is a $\Z[t]$-module homomorphism \[j: H_{\infty}\hookrightarrow \invlim\, H_{T_{n}}^{*}(\mathcal{F}_{n})\/.\] The injectivity of localization maps in the finite-dimensional setting implies that $j$ is injective as well.


\section{Universal localization map}\label{sec:UnivLoc}
In this section, we introduce a $\Z[t]$-algebra $\R_\infty$ 
and establish an explicit isomorphism from 
$\R_{\infty}$ onto 
$H_\infty,$ the $\Z[t]$-module spanned by the stable Schubert
classes. This isomorphism will be used 
in the proof of the existence of 
the double Schubert polynomials from \S \ref{sec:DSP}. 
\subsection{The ring $\R_{\infty}$ and the universal localization map}\label{ssec:UnivLoc}
Set $\Z[z]=\Z[z_{1},z_{2},z_{3},\ldots]$ and define the following rings:
$$R_{\infty}:=\Z[t]
\otimes_{\Z}\Z[z]
\otimes_{\Z} \Gamma,\quad \textrm{ and } \quad 
R_{\infty}':=\Z[t]
\otimes_{\Z}\Z[z]
\otimes_{\Z} \Gamma' \/.
$$
As usual, we will use $\R_{\infty}$ to denote $R_\infty$ for type $\mathrm{C}$ and
$R_{\infty}'$ for types $\mathrm{B}$ and $\mathrm{D}$.

We introduce next the most important algebraic tool of the paper.
Let $v$ be in $W_{\infty}.$
Set $t_{v}=(t_{v,1},t_{v,2},\ldots)$ to be 
$$t_{v,i}=\begin{cases}t_{\overline{v(i)}}&
\mbox{if}\;v(i)\;\mbox{is negative}\\
0 & \mbox{otherwise}
\end{cases},
$$
where we set $t_{\overline{i}}$ to be $-t_{i}.$
Define a homomorphism of $\Z[t]$-algebras
$$
\Phi_{v}: R_{\infty}'\longrightarrow \Z[t]
\quad
\left(x\mapsto
t_{v},\quad
z_{i}\mapsto t_{v(i)}\right).$$
Note that since $v(i)=i$ for all 
sufficiently large $i$, the substitution $x\mapsto t_{v}$
to $P_\lambda(x)$ gives a {\em polynomial} $P_{\lambda}(t_{v})$ in $\Z[t]$ (rather than a formal power series). Since $R_{\infty}$ is a subalgebra of $R_{\infty}'$, the restriction map
$\Phi_v:R_{\infty}\longrightarrow \Z[t]$ sends $Q_{\lambda}(x)$ to $Q_{\lambda}(t_{v}).$ 
\begin{Def} Define the "universal localization map" 
to be the homomorphism of $\Z[t]$-algebras given
by
$$\Phi:
\R_{\infty}\longrightarrow
\prod_{v\in \W_{\infty}}\Z[t],\quad
f\mapsto (\Phi_{v}(f))_{v\in \W_{\infty}}.
$$
\end{Def}
\begin{remark}{\rm A geometric interpretation of the map $\Phi$, in terms of the usual localization map, will be given later in \S  \ref{sec:geometry}.}\end{remark}

The main result of this section is:
\begin{thm}\label{thm:isom}
The map
$\Phi$ is an isomorphism
of graded $\Z[t]$ algebras
from $\R_\infty$ onto its image. Moreover, the image of $\Phi$ is equal to $H_\infty$. 
\end{thm}
\begin{cor}\label{cor:Hsubalg}
$H_{\infty}$ is a $\Z[t]$-subalgebra in $\prod_{v\in \W_{\infty}}\Z[t].$
\end{cor}
The proof of the theorem will be given in several lemmata and propositions, and it occupies the remaining part of section 6. The more involved part is to show surjectivity, which relies on the analysis of the "transition equations" implied by the equivariant Chevalley rule, and on study of factorial $P$ and $Q$- Schur functions. The proof of injectivity is rather short, and we present it next.
\begin{lem} \label{injective}
The map $\Phi$ is injective.
\end{lem}
\begin{proof}
We first consider type $\mathrm{B}$ case. 
Write $f\in R_\infty'$ as
$f=\sum_{\lambda}c_\lambda(t,z) P_\lambda(x).$
Suppose $\Phi(f)=0.$
There are $m,n$ such that
$$c_\lambda\in \Z[t_1,\ldots,t_m,z_1,\ldots,z_n]$$
for all $\lambda$ such that $c_\lambda\ne 0.$
Define $v\in W_\infty$ by 
$v(i)=m+i\;(1\leq i\leq n),$
$ v(n+i)=i\;(1\leq i\leq m),$ 
$v(m+n+i)=\overline{m+n+i}\;(1\leq i\leq N),$
and $
v(i)=i\;(i>N),$
where $N\geq m+n+1.$
Then we have
$\Phi_v(f)=\sum_\lambda c_\lambda(t_1,\ldots,t_m;t_{m+1},\ldots,t_{m+n})
P_\lambda(t_{m+n+1},t_{m+n+2},\ldots,t_{m+n+N})=0.$
Since this holds for all sufficiently large $N,$
we have $$\sum_\lambda c_\lambda(t_1,\ldots,t_m;t_{m+1},\ldots,t_{m+n})
P_\lambda(t_{m+n+1},t_{m+n+2},\ldots)=0.$$
Since $P_\lambda(t_{m+n+1},t_{m+n+2},\ldots)$ are linearly
independent over $\Z$ 
(see \cite{Mac}, III, (8.9)),
we have $$c_\lambda(t_1,\ldots,t_m;t_{m+1},\ldots,t_{m+n})=0$$
for all $\lambda.$
This implies $c_\lambda(t_1,\ldots,t_m;z_{1},\ldots,z_{n})=0$
for all $\lambda.$
Since $R_{\infty}\subset R_{\infty}',$
type $\mathrm{C}$ case follows immediately.
Type $\mathrm{D}$ case is proved by a minor modification.
Take $N$ to be always even, and consider 
the sufficiently large {\it even} $N.$
\end{proof}

\subsection{Factorial $Q$-(and $P$-)functions 
and Grassmannian Schubert classes}
Recall that there is a natural bijection between
$W_{\infty}^{0},W_{\infty}^{\hat{1}}$ and the set of strict partitions $\mathcal{SP}$.
The next result was proved by Ikeda in \cite{Ik} for type 
$\mathrm{C}$, and Ikeda-Naruse in \cite{IN} for types $\mathrm{B,D}$.

\begin{thm}[\cite{Ik},\cite{IN}] \label{PhiFacQ}
Let $\lambda\in \mathcal{SP}$ and $w_{\lambda}
\in W_{\infty}^{0}$ and $w'_{\lambda}
\in W_{\infty}^{\hat{1}}$ be the corresponding Grassmannian
elements.
Then we have
\begin{enumerate}
\item $\Phi\left(Q_{\lambda}(x|t)\right)=\sigma_{w_{\lambda}}^{(\infty)}\;
\mbox{for type}\;C,$
\item $\Phi\left(P_{\lambda}(x|0,t)\right)
=\sigma_{w_{\lambda}}^{(\infty)}\mbox{for type}\;B,$
\item 
$\Phi\left(P_{\lambda}(x|t)\right)
=\sigma_{w'_{\lambda}}^{(\infty)}\;
\mbox{for type}\;D.
$
\end{enumerate}
\end{thm}
\begin{proof} We consider first the type $\mathrm{C}$ case.
The map on $W_{\infty}$
given by $v\mapsto \sigma_{w_{\lambda}}^{(\infty)}|_{v}$
is constant on each left coset of $W_{\infty,0}\cong 
S_{\infty}$
and it is determined by the values at the Grassmannian elements.
Let $v\in W_{\infty}$ and $w_{\mu}$ be the minimal length representative
of the coset $vS_{\infty}$ corresponding to 
a strict partition $\mu.$
Then $t_{v}$ defined in \S \ref{ssec:UnivLoc}
is a permutation of $t_{\mu}.$
Since $Q_{\lambda}(x|t)$ is 
symmetric with respect to $x$
we have 
$\Phi_{v}(Q_{\lambda}(x|t))=Q_{\lambda}(t_{v}|t)
=Q_{\lambda}(t_{\mu}|t).
$
In \cite{Ik}, it was shown that 
$Q_{\lambda}(t_{\mu}|t)=\sigma_{w_{\mu}}^{(\infty)}|_{w_{\mu}},$
which is equal to $\sigma_{w_{\mu}}^{(\infty)}|_v.$
This completes the proof in this case.
Proofs of the other cases are the same 
with appropriate identification
of the functions and strict partitions.
\end{proof}

\subsection{Equivariant Chevalley formula} 
The {\em 
Chevalley formula} is a rule to multiply a Schubert class with a divisor class.
To state it we need some notation.
For a positive root $\alpha\in R^{+}$ and
a simple reflection $s_{i}$,
set
$$
c_{\alpha,s_{i}}=(\omega_i,\alpha^{\vee}),\quad
\alpha^{\vee}=2\alpha/(\alpha,\alpha),
$$
where $\omega_{i}$ 
is the $i$-th fundamental weight 
of one of the classical types $\mathrm{A_{n}-D_n}$ for sufficiently large $n$.
The number $c_{\alpha,s_{i}}$ - called  {\it Chevalley multiplicity} - does not depend on the choice of $n.$
\begin{prop}[cf. \cite{KK}]\label{prop:codim1}
For any $w\in \W_{\infty}$, the Chevalley multiplicity $\sigma_{s_{i}}^{(\infty)}|_{w}$ is given 
by $\omega_{i}-w(\omega_{i})$,  where
$\omega_{i}$ is the fundamental weight 
for a classical type $\mathrm{A_n-D_n}$ such that $n\geq i.$
\end{prop}

\begin{lem}[Equivariant Chevalley formula]
\label{lem:eqCh}
$$
\sigma_{s_{i}}^{(\infty)}
\sigma_{w}^{(\infty)}
=
\displaystyle\sum_{\alpha\in R^{+},\;\ell(w s_\alpha)=\ell(w)+1} 
c_{\alpha,s_{i}}\,\sigma_{w s_\alpha}^{(\infty)} + \sigma_{s_{i}}^{(\infty)}|_{w} \cdot \sigma_{w}^{(\infty)}. $$
\end{lem}
\begin{proof}
The non-equivariant case is due to Chevalley \cite{C}, but the stable version of this formula was given in  \cite{B93}.
An easy argument using localization shows that the only difference in the equivariant case is the appearance of the 
equivariant term
$\sigma_{s_{i}}^{(\infty)}|_{w} \cdot \sigma_{w}^{(\infty)}.$
\end{proof}
\begin{remark}{\rm
There are only finitely many nonzero terms in the sum
in the right hand side.}
\end{remark}

\begin{lem} \label{lem:z}
The elements $\Phi(z_{i}) \in H_\infty$ 
are expressed in terms of Schubert classes as follows:

\noindent Type $\mathrm{B}$:
$\Phi(z_1)=\sigma_{s_1}^{(\infty)}-2\sigma_{s_{0}}^{(\infty)}+t_{1},\;
\Phi(z_i)=\sigma_{s_i}^{(\infty)}-\sigma_{s_{i-1}}^{(\infty)}+t_{i}
\; (i\geq 2),$

\noindent Type $\mathrm{C}$:
$\Phi({z}_i)=\sigma_{s_i}^{(\infty)}-\sigma_{s_{i-1}}^{(\infty)}+t_{i}
\; (i\geq 1),$

\noindent Type $\mathrm{D}$:
$\Phi({z}_1)=\sigma_{s_{1}}^{(\infty)}-\sigma_{s_{\hat{1}}} ^{(\infty)} +t_{1},\;
\Phi({z}_2)=\sigma_{s_{2}} ^{(\infty)}-\sigma_{s_{1}} ^{(\infty)}
-\sigma_{s_{\hat{1}}} ^{(\infty)} +t_{2}$, \quad and \\ 
$\Phi({z}_i)=\sigma_{s_{i} }^{(\infty)}-\sigma_{s_{i-1}} ^{(\infty)} +t_{i}\;(i\geq 3).
$
\end{lem}
\begin{proof}
This follows by localizing both sides of the formulas, and then using Prop. \ref{prop:codim1}.
\end{proof}

\begin{lem}\label{lem:PhisubH} We have $\mathrm{Im}(\Phi)\subset H_\infty.$
\end{lem}
\begin{proof}
The ring $R_{\infty}$ has a $\Z[t]$-basis
$z^{\alpha}Q_{\lambda}(x|t)$ where $z^{\alpha}$
are monomials in $\Z[z]$ and $\lambda$
are strict partitions.
Since $\Phi$ is $\Z[t]$-linear, it is enough to show that
$\Phi\left(z^{\alpha}Q_{\lambda}(x|t)\right)$
belongs to $H_{\infty}.$
We use induction on degree $d$ of the monomial $z^{\alpha}.$
The case $d=0$ holds by Thm. \ref{PhiFacQ}.
Let $d\geq 1$ and assume that 
$\Phi(z^{\alpha}Q_{\lambda}(x|t))$ lies in $H_{\infty}$
for any monomial $z^{\alpha}$ of degree less than $d.$ 
Note that, by Lem. \ref{lem:z}, we have 
$\Phi(z_{i})\in H_{\infty}.$
Choose any index $i$ such 
that $z^{\alpha}=z_{i}\cdot z^{\beta}.$
By induction hypothesis $\Phi\left(z^{\beta}Q_{\lambda}(x|t)
\right)$ is an element in $H_{\infty},$
i.e., a linear combination
of $\sigma_{w}^{(\infty)}$'s with coefficients
in $\Z[t].$ Lem. \ref{lem:z} together with 
equivariant Chevalley formula imply
that $\Phi(z_{i})\sigma_{w}^{(\infty)}$
belongs to $H_{\infty}.$
It follows that $z^{\alpha}Q_{\lambda}(x|t)$
belongs to $H_{\infty}.$
\end{proof}

\subsection{Transition equations}\label{ssec:trans} To finish the proof of surjectivity of $\Phi$, we need certain recursive relations for the Schubert classes - the {\it transition equations} -  implied by the (equivariant) Chevalley formula. The arguments in this subsection are very similar to those given by S. Billey in \cite{B93}. 
Let $t_{ij}$ denote the reflection with respect to the root $t_{j}-t_{i}$,
$s_{ij}$ the reflection with respect to $t_{i}+t_{j}$ and $s_{ii}$ the reflection
with respect to $t_{i}$ or $2t_{i}$ (depending on type).
From now on we regard $\Z[z]$
as subalgebra of $H_{\infty}$ via $\Phi$
and we identify $z_i$ with its image $\Phi(z_{i})$ in $H_{\infty}$
(cf. Lem \ref{lem:PhisubH}).  

\begin{prop}[Transition equations]
The Schubert classes $\sigma_{w}$ of types
$\mathrm{B,C}$ and $\mathrm{D}$ satisfy the following recursion formula:
\begin{equation}
\sigma_w ^{(\infty)} =({z}_r-v(t_r))\;\sigma_v ^{(\infty)}
+
\sum_{1\leq i<r} \sigma^*_{v t_{ir}}+
\sum_{i\neq r}\sigma^*_{v s_{ir}}+
\chi\sigma^*_{v s_{rr}},\label{Transition}
\end{equation}
where $r$ is the last descent of $w$,
$s$ is the largest index such that $w(s)<w(r)$,
$v=wt_{rs}$,
$\chi=2,1,0$ according to the types $\mathrm{B,C,D}$,
and for each 
$\sigma_{vt}^*=0$ unless $\ell(vt)
=\ell(v)+1=\ell(w)$ for $v,t\in \W_{\infty}$
in which case $\sigma_{vt}^{*}=\sigma_{vt}^{(\infty)}.$
\end{prop}
\begin{proof}
The same as in \cite[Thm.4]{B93}
using the equivariant Chevalley formula (Lemma \ref{lem:eqCh}).
\end{proof}

\begin{remark}{\rm 
The precise recursive nature of the equation (\ref{Transition}) will be explained
in the proof of the next Proposition below.}
\end{remark}
\begin{prop} \label{TransExp}
If $w\in \W_{n}$ then 
the Schubert class $\sigma_{w}^{(\infty)}$
is expressed as a $\Z[z,t]$-linear 
combination of the Schubert classes of maximal Grassmannian type.
More precisely we have
\begin{equation}
\sigma_{w}^{(\infty)}=\sum_{\lambda}
g_{w,\lambda}({z},t)\sigma_{\lambda}^{(\infty)},\label{expansion}
\end{equation}
for some polynomials
$g_{w,\lambda}({z},t)$ in variables $t_i$ and $z_i$, and the sum is over strict partitions $\lambda$ such that $|\lambda|\leq n.$
\end{prop}
\begin{proof}
We will show that the recursion (\ref{Transition})
terminates in a finite number of steps 
to get the desired expression.
Following \cite{B93}, we define a partial ordering on the elements of $\W_{\infty}.$
Given $w$ in $\W_{\infty}$,
let $LD(w)$ be the position of the last descent.
Define a partial ordering on the elements
of $\W_{\infty}$ by $w<_{LD}u$
if $LD(u)<LD(w)$ or if $LD(u)=LD(w)$ and 
$u(LD(u))<w(LD(w))$.
In \cite{B93} it was shown that each element 
appearing on the right hand side of
(\ref{Transition}) is less than $w$ under this ordering.
Moreover it was proved in \cite[Thm.4]{B93} that 
recursive applications of 
(\ref{Transition}) give only terms which correspond to the elements in  
$\W_{n+r}$ where $r$ is the last descent of $w.$
Therefore we obtain the expansion (\ref{expansion}).
\end{proof}

\subsection{Proof of Theorem \ref{thm:isom}}
\begin{proof} By Lem. \ref{lem:PhisubH} we know
$\mathrm{Im}(\Phi)\subset H_{\infty}.$
Clearly $\Phi$ preserves the degree.
So it remains to show
$H_{\infty}\subset\mathrm{Im}(\Phi).$
In order to show this, it suffices to 
$\sigma_{w}^{(\infty)}\in \mathrm{Im}(\Phi).$
In fact we have
\begin{equation}
\Phi\left(\sum_{\lambda}
g_{w,\lambda}(z,t)Q_{\lambda}(x|t)
\right)=\sigma_{w}^{(\infty)}.
\label{eq:DefSch}
\end{equation}
since 
$\Phi$ is $\Z[z,t]$-linear.
\end{proof}

\section{Weyl group actions and divided difference operators on $\R_{\infty}$}\label{sec:WactsR}
We define two commuting actions of $\W_{\infty}$ on the ring $\R_{\infty}.$
It is shown that the Weyl group actions
are compatible with the action on $H_{\infty}$
via $\Phi.$

\subsection{Weyl group actions on $R_{\infty}$}
We start from type $\mathrm{C}$. We make $W_{\infty}$ act as ring automorphisms
on $R_{\infty}$ 
by letting 
$s_{i}^{z}$ interchange $z_{i}$ and $z_{i+1},$
for $i>0,$ and letting
$s_{0}^{z}$  replace $z_{1}$ and $-z_{1},$
and also
$$
s_{0}^{z}Q_{i}(x)=Q_{i}(x)+2\sum_{j=1}^{i}z_{1}^{j}Q_{i-j}(x).
$$
The operator $s_{0}^{z}$ was introduced in 
\cite{BH}. 
Let $\omega: R_{\infty}\rightarrow R_{\infty}$ be an involutive ring automorphism
defined by  
$$
\omega(z_{i})=-t_{i},\quad
\omega(t_{i})=-z_{i},\quad
\omega(Q_{k}(x))=Q_{k}(x).
$$
Define the operators $s_{i}^{t}$ on $R_{\infty}$
by $s_{i}^{t}=\omega s_{i}^{z}\omega$
for $i\in I_{\infty}.$ 
More explicitly, $s_{i}^{t}$ interchange $t_{i}$ and 
$t_{i+1}$, for $i>0$, and $s_{0}^{t}$
replace $t_{1}$ and $-t_{1}$ and also
$$
s_{0}^{t}Q_{i}(x)=Q_{i}(x)+2\sum_{j=1}^{i}(-t_{1})^{j}Q_{i-j}(x).
$$
\begin{lem}\label{lem:s0super} The action of operators 
$s_{0}^{z},s_{0}^{t}$ on any $\varphi(x)\in \Gamma$ 
are written as
$$s_{0}^{z}\varphi(x_{1},x_{2},\ldots)
=\varphi(z_{1},x_{1},x_{2},\ldots),\quad
s_{0}^{t}\varphi(x_{1},x_{2},\ldots)
=\varphi(-t_{1},x_{1},x_{2},\ldots).$$
\end{lem}
Note that the right hand side
of both the formulas above belong to 
$R_{\infty}.$
\begin{proof}
We show this for the generators $\varphi(x)=Q_{k}(x)$
of $\Gamma.$ 
By the definition of $s_{0}^{z}$ we have
$$
\sum_{k=0}^{\infty}
s_{0}^{z}
Q_{k}(x)\cdot u^{k}
=
\left(\prod_{i=1}^{\infty}
\frac{1+x_{i}u}{1-x_{i}u}
\right)\frac{1+z_{1}u}{1-z_{1}u}
=\sum_{k=0}^{\infty} Q_{k}(z_{1},x_{1},x_{2},\ldots)u^{k}.
$$
Thus we have the result for $s_{0}^{z}Q_{k}(x)$
for $k\geq 1.$
The proof for $s_{0}^{t}$ is similar.
\end{proof}

\begin{prop} 
\begin{enumerate}
\item The operators
$s_{i}^{z}\;(i\geq 0)$ give an action of $W_{\infty}$ on $R_{\infty},$
\item The operators 
$s_{i}^{t}\;(i\geq 0)$ give an action of $W_{\infty}$ on $R_{\infty},$
\item The two actions of $W_{\infty}$ commute with
each other.
\end{enumerate}
\end{prop}
\begin{proof}
We show that $s_{i}^{z}$ satisfy the Coxeter relations
for $W_{\infty}.$
The calculation for $s_{i}^{t}$ is the same.
We first show that $(s_{0}^{z})^{2}=1.$
For $f(z)\in \Z[z]$, $(s_{0}^{z})^{2}f(z)=f(z)$ is obvious.
We have for $\varphi(x)\in \Gamma,$
$$
(s_{0}^{z})^{2}
(\varphi(x))=
s_{0}^{z}\varphi(z_{1},x_{1},x_{2},\ldots)
=\varphi(z_{1},-z_{1},x_{1},x_{2},\ldots)
=\varphi(x_{1},x_{2},\ldots),
$$
where we used the super-symmetry
(Lemma \ref{lem:super})
at the last equality.
The verification of the remaining relations and the commutativity are left for the reader.
\end{proof}

In type $\mathrm{B}$, the action of $W_\infty$ on $R_\infty'$ is obtained by extending in the canonical way the action from $R_{\infty}.$ Finally, we consider the type $\mathrm{D}$ case.
In this case, the action is given by restriction the action of $W_{\infty}$ on $R'_\infty$ to the subgroup $W_{\infty}'.$
Namely, if we
set $s_{\hat{1}}^{z}=s_{0}^{z}s_{1}^{z}s_{0}^{z}$
and $s_{\hat{1}}^{t}=s_{0}^{t}s_{1}^{t}s_{0}^{t}$, then
we have the corresponding 
formulas for $s_{\hat{1}}^{t}$ and $s_{\hat{1}}^{t}$ (in type D):
$$
s_{\hat{1}}^{z}\varphi(x_{1},x_{2},\ldots)
=\varphi(z_{1},z_{2},x_{1},x_{2},\ldots),\quad
s_{\hat{1}}^{t}\varphi(x_{1},x_{2},\ldots)
=\varphi(-t_{1},-t_{2},x_{1},x_{2},\ldots).
$$

\subsection{Divided difference operators}\label{ssec:divdiffgeom}
The divided difference operators 
on $\R_{\infty}$ are defined by
$$
\partial_{i}f=\frac{f-s_{i}^{z}f}{\omega(\alpha_{i})},\quad
\delta_{i}f=\frac{f-s_{i}^{t}f}{\alpha_{i}},
$$
where $s_{i}$ and $\alpha_{i}\;(i\in \I_{\infty})$ are the
simple reflections and the corresponding simple roots.
Clearly we have
$\delta_{i}=\omega\partial_{i}\omega \quad (i\in \I_{\infty}).$

\subsection{Weyl group action and 
commutativity with divided difference operators}

\begin{prop}\label{prop:comm}
We have
$
(1)\;s_{i}^{L}\Phi=\Phi s_{i}^{t},\;
(2)\; s_{i}^{R}\Phi=\Phi s_{i}^{z}.
$
\end{prop}
\begin{proof} We will only prove this for type C; the other types can be treated similarly.
We first show $(1).$
This is equivalent to 
$s_{i}\left(\Phi_{s_{i}v}(f)\right)=\Phi_{v}(s_{i}^{t}f)$
for all $f\in R_{\infty}.$
If $f\in \Z[z,t]$ the computation is 
straightforward and we omit the proof.
Suppose $f=\varphi(x)\in \Gamma$.
We will only show 
$s_{0}\left(\Phi_{s_{0}v}(f)\right)=\Phi_{v}(s_{0}^{t}f)$
since the case $i\geq 1$ is straightforward.

By Lem. \ref{lem:s0super}, the right hand side of this
equation is written as
\begin{equation}
\varphi(-t_1,x_1,x_2,\ldots)|_{x_j=t_{v,j}}.
\label{eq:rhs-t}
\end{equation}
Let $k$ be the (unique) index such that $v(k)=
1$ or $\overline{1}.$
Then the string $t_{s_0v}$ differs
from $t_v$ only in $k$-th position.
If $v(k)=\overline{1}$,
then 
$t_{v,k}=t_1,\, t_{s_0v,k}=0$
and $t_{v,j}=t_{s_0v,j}$ for $j\ne k.$
In this case (\ref{eq:rhs-t}) is
$$
\varphi(-t_1,t_{v,1},\ldots,t_{v,k-1},t_1,t_{v,k+1},\ldots).
$$
This polynomial is equal to 
$\varphi(t_{v,1},\ldots,t_{v,k-1},t_{v,k+1},\ldots)$
because $\varphi(x)$ is supersymmetric.
It is straightforward to see that
$s_{0}\Phi_{s_{0}v}(\varphi(x))$
is equal to $\varphi(t_{v,1},\ldots,t_{v,k-1},t_{v,k+1},\ldots).$ The case for $v(k)={1}$ is easier, so 
we left it to the reader.

Next we show (2), i.e. $
\Phi_{vs_{i}}(f)=\Phi_{v}(s_{i}^{z}f)
$ for all $f\in R_{\infty}.$
Again, the case $f\in \Z[z,t]$ is straightforward,
so we we omit the proof of it.
We show $\Phi_{vs_{0}}\left(\varphi(x)\right)
=\Phi_{v}(s_{0}^{z}\varphi(x))$
for $\varphi(x)\in \Gamma.$
The right hand side is 
\begin{equation}
\varphi(z_{1},x_{1},x_{2},\ldots)|_{z_{1}=v(t_{1}),\,
x_{j}=t_{v,j}},
\label{eq:rhs}
\end{equation}
where $t_{v,j}=t_{\overline{v(j)}}$ if 
$v(j)$ is negative and otherwise $t_{v,j}$ is zero.
If $v(1)=-k$ is negative, the above function (\ref{eq:rhs}) is
$$
\varphi(-t_{k},t_{k},t_{v,2},t_{v,3},\ldots).
$$
This is equal to 
$\varphi(0,0,t_{v,2},t_{v,3},\ldots)$
because $\varphi$ is supersymmetric.
Then also this is equal to $\varphi(0,0,t_{v,2},t_{v,3},\ldots)
=\varphi(0,t_{v,2},t_{v,3},\ldots)$ by stability
property.
Now since $\overline{v(1)}$ is positive we have
$t_{vs_{0}}=(0,t_{v,2},t_{v,3},\ldots).$
Therefore the polynomial (\ref{eq:rhs})
coincides with $\Phi_{vs_{0}}(\varphi(x)).$
If $v(1)$ is positive, then 
$t_{v}=(0,t_{v,2},t_{v,3},\ldots)$ and 
$t_{vs_{0}}=(t_{v(1)},t_{v,2},t_{v,3},\ldots).$
Hence the substitution 
$x\mapsto t_{vs_{0}}$ to 
the function $\varphi(x_{1},x_{2},\ldots)$
gives rise to the polynomial (\ref{eq:rhs}).
Next we show $
\Phi_{vs_{i}}(\varphi(x))=\Phi_{v}(s_{i}^{z}\varphi(x))
$
for $i\geq 1.$ 
First recall that $s_{i}^{z}\varphi(x)=\varphi(x).$ 
In this case
$t_{vs_{i}}$ is obtained from
$t_{v}$ by exchanging $t_{v,i}$ and $t_{v,i+1}.$
So $\varphi(t_{vs_{i}})=\varphi(t_{v}).$ 
This completes the proof.
\end{proof}

Using the above proposition, the next 
result follows:
\begin{prop}\label{prop:PhiCommD} The localization map
$\Phi:\R_{\infty}\rightarrow H_{\infty}$
commutes with the divided difference operators 
both on $\R_{\infty}$ and $H_{\infty},$ i.e., 
$$\Phi \,\partial_{i}=\partial_{i}\,\Phi,\quad
\Phi\,\delta_{i}=\delta_{i}\,\Phi$$
\end{prop}
\begin{proof}
Let $f\in R_{\infty}.$
Applying $\Phi$ on the both hand sides of
equation 
$
{\omega(\alpha_{i})}\cdot
\partial_{i}f={f-s_{i}^{z}f}
$ we have 
$\Phi(-\omega(\alpha_{i}))\cdot \Phi(\partial_{i}f)
=\Phi(f)-s_{i}^{R}\Phi(f)$, 
where we used Prop. \ref{prop:comm} and linearity. 
Localizing at $v$ we obtain
$v(\alpha_{i})\cdot \Phi_{v}(\partial_{i}f)
=\Phi_{v}(f)-\Phi_{vs_{i}}(f).$
Note that we used the definition of $s_{i}^{R}$ and
$\Phi_{v}(\omega(\alpha_{i}))=-v(\alpha_{i}).$
The proof for the statement regarding $\delta_{i}$ 
is similar, using $\Phi(\alpha_{i})=\alpha_{i}.$
\end{proof}

\subsection{Proof of the existence and uniqueness Theorem \ref{existC}}
%
\begin{proof} (Uniqueness)
Let  
$\{\mathfrak{S}_{w}\}$ and $\{\mathfrak{S}_{w}'\}$
be two families 
both satisfying the defining conditions 
of the double Schubert polynomials.
By induction on the length of $w$,
we see
$\partial_{i}(\mathfrak{S}_{w}-\mathfrak{S}'_{w})=
\delta_{i}(\mathfrak{S}_{w}-\mathfrak{S}'_{w})=0$
for all $i\in \I_{\infty}.$
This implies that 
the difference $\mathfrak{S}_{w}-\mathfrak{S}'_{w}$ is
invariant for both left and right actions of $\W_{\infty}.$
It is easy to see that the only such invariants 
in $\R_{\infty}$ are the constants.
So $\mathfrak{S}_{w}-\mathfrak{S}'_{w}=0$
by the constant term condition. 

(Existence) Define $\mathfrak{S}_{w}(z,t;x)
=\Phi^{-1}(\sigma_{w}^{(\infty)}).$
By Prop. \ref{prop:PhiCommD}
and Prop. \ref{prop:divdiff},
$\mathfrak{S}_{w}(z,t;x)$
satisfies the defining equations for the double
Schubert polynomials.
The conditions on the constant term are 
satisfied since $\sigma_{w}^{(\infty)}$
is homogeneous of degree $\ell(w)$ (Prop. \ref{StableClass})
and we have $\mathfrak{S}_{e}=1.$
\end{proof}

\begin{remark}\label{rem:TransDSP}
{\rm By construction, $\mathfrak{S}_{w}(z,t;x)$ 
satisfies the transition equation (\ref{Transition}) 
with $\sigma_{w}^{(\infty)}$ replaced by $\mathfrak{S}_{w}(z,t;x)$.
This equation provides an effective way to calculate 
the double Schubert polynomials.}
\end{remark}
%

\begin{remark}\label{rem:A}{\rm The ring
$\Z[z]\otimes_{\Z}\Z[t]$ is
stable under the actions of the divided 
difference operators $\partial_{i},\delta_{i}\,(i\geq 1)$ of type $\mathrm{A}$, and the type $\mathrm{A}$ double Schubert polynomials $\mathfrak{S}_{w}^{A}(z,t)$,
$w\in S_{\infty}$ form the unique family of solutions of the system of equations
involving only $\partial_{i},\delta_{i}$ for $i\geq 1$, and which satisfy the constant term conditions.}\end{remark}

\subsection{Projection to the cohomology of flag manifolds}\label{ssec:projection}
We close this section with a brief discussion 
of the projection from $\R_{\infty}$ onto
$H_{T_{n}}^{*}(\mathcal{F}_{n}).$
For $f\in \Z[t]$, we denote by $f^{(n)}\in \Z[t]^{(n)}$
the polynomial given  
by setting $t_{i}=0$ for $i>n$ in $f.$
Let $\mathrm{pr}_{n}:
H_{\infty}\rightarrow H_{n}$ be the projection
given by 
$(f_{v})_{v\in \W_{\infty}}\mapsto 
(f_{v}^{(n)})_{v\in \W_{n}}.$
Consider the following composition of maps
$$
\pi_{n}: \R_{\infty}\overset{\Phi}{\longrightarrow} H_{\infty}
\overset{\mathrm{pr}_{n}}{\longrightarrow} H_{n}
\cong
H_{T_{n}}^{*}
(\mathcal{F}_{n}).
$$
Explicitly, we have
$
\pi_{n}(f)|_{v}=\Phi_{v}(f)^{(n)}\;
(f\in \R_{\infty},\; v\in \W_{n}).
$
We will give an alternative geometric description
for $\pi_{n}$ in Section \ref{sec:geometry}.

\begin{prop}\label{prop:piCommD}
We have
$\pi_{n}(\mathfrak{S}_{w})=\sigma_{w}^{(n)}$
for $w\in \W_{n}$ and $\pi_{n}(\frak{S}_{w})=0$
for $w\notin\W_{n}.$
Moreover $\pi_{n}$ commutes with divided difference operators
$$
\partial_{i}^{(n)}\circ\pi_{n}=\pi_{n}\circ\partial_{i},\quad
\delta_{i}^{(n)}\circ\pi_{n}=\pi_{n}\circ\delta_{i}
\quad
(i\in \I_{\infty}),
$$ 
where $\partial_{i}^{(n)},\delta_{i}^{(n)}$ are divided difference
operators on $H_{T_{n}}^{*}(\mathcal{F}_{n}).$
\end{prop}
\begin{proof}
The first statement follows from the construction
of $\sigma_{w}^{(\infty)}$ 
and the vanishing property (Prop. \ref{StableClass}).
The second statement follows from Prop. \ref{prop:PhiCommD}
and the commutativity  $\partial_{i}^{(n)}\circ\mathrm{pr}_{n}=
\mathrm{pr}_{n}\circ\partial_{i},\,
\delta_{i}^{(n)}\circ\mathrm{pr}_{n}=
\mathrm{pr}_{n}\circ\delta_{i}$
which is obvious from the construction of 
$\partial_{i},\delta_{i}.$
\end{proof}
\begin{cor} There exists an injective homomorphism of $\Z[t]$-algebras $\pi_\infty: \R_\infty \to \invlim \eqcoh(\F_n)$.\end{cor}
\begin{proof} The proof follows from the above construction and \S \ref{ssec:invlim}.\end{proof}


%


\section{Double Schubert polynomials}\label{sec:DSP}

\subsection{Basic properties} Recall that the {\em double Schubert polynomial} $\mathfrak{S}_w(z,t;x)$ is equal to the inverse image of the stable Schubert class $\sigma_w^{(\infty)}$ under the algebra isomorphism $\Phi:\R_\infty \to H_\infty$. In the next two sections we will study the algebraic properties of these polynomials.
\begin{thm}\label{T:properties}
The double Schubert
polynomials satisfy the following:
\begin{enumerate}
\item (Basis)
The double Schubert polynomials $\{\mathfrak{S}_{w}\}_{w
\in \W_{\infty}}$
form a $\Z[t]$-basis of 
$\R_{\infty}.$
\item (Relation to Billey-Haiman's polynomials) For all $w\in \W_{\infty}$ we have
\begin{equation}
\mathfrak{S}_{w}(z,0;x)=\mathfrak{S}_{w}(z;x),
\label{eq:relBH}
\end{equation}
where $\mathfrak{S}_{w}(z;x)$ denotes
Billey-Haiman's polynomial.
\item (Symmetry) We have
$\mathfrak{S}_{w}(-t,-z;x)=\mathfrak{S}_{w^{-1}}(z,t;x).$
\end{enumerate}
\end{thm}
\begin{proof}
Property (1) holds because the stable Schubert classes $\sigma_w^{(\infty)}$ form a $\Z[t]$-basis for $H_\infty$ (cf. \S\ref{ssec:SSch}).
Property (2) holds because $\mathfrak{S}_{w}(z,0;x)\in \Z[z]\otimes \Gamma'$
satisfies the defining conditions for 
Billey-Haiman's polynomials
involving the right divided difference operators $\partial_{i}.$ 
Then by the uniqueness of Billey-Haiman's polynomials,
we have the results.
For (3), set $\mathfrak{X}_{w}=\omega(\mathfrak{S}_{w^{-1}}).$
Then by the relation $\delta_{i}=\omega\partial_{i}\omega$
we can show that 
$\{\mathfrak{X}_{w}\}$ satisfies the defining 
conditions of the double Schubert polynomials.
So the uniqueness of the double Schubert polynomials
implies $\mathfrak{X}_{w}=\mathfrak{S}_{w}.$
Then we have
$\omega(\mathfrak{S}_{w})
=\omega(\mathfrak{X}_{w})=
\omega(\omega\mathfrak{S}_{w^{-1}})=
\mathfrak{S}_{w^{-1}}.$
\end{proof}

\begin{remark}{\rm
For type $D$ we have
$s_{0}^{z}s_{0}^{t}\mathfrak{D}_{w}=\mathfrak{D}_{\hat{w}}$
where $\hat{w}$ is the image of $w$ under the involution
of $W_{\infty}'$ given by interchanging
$s_{1}$ and $s_{\hat{1}}.$
This is shown by the uniqueness of 
solution as in the proof of the symmetry property.
See \cite[Cor. 4.10]{BH} for the corresponding 
fact for the Billey-Haiman polynomials.}
\end{remark}

\subsection{Relation to type $\mathrm{A}$ double Schubert polynomials}
Let $\mathfrak{S}_{w}^{A}(z,t)$ denote the 
type $\mathrm{A}$ double Schubert polynomials.
Recall that $\W_{\infty}$
has a parabolic subgroup 
generated by $s_{i}\;(i\geq 1)$ which is
isomorphic to
$S_{\infty}.$ 

\begin{lem}\label{lem:typeA}
Let $w\in \W_{\infty}.$
If $w\in S_{\infty}$ then 
$\mathfrak{S}_{w}(z,t;0)=\mathfrak{S}_{w}^{A}(z,t)
$ otherwise we have
$\mathfrak{S}_{w}(z,t;0)=0.$
\end{lem}
\begin{proof} The polynomials
$\{\mathfrak{S}_{w}(z,t;0)\}$, $w\in S_{\infty},$
in $\Z[t]\otimes_{\Z}\Z[z]\subset \R_{\infty}$ satisfy the defining divided difference
equations for the double Schubert polynomials of type $\mathrm{A}$
(see Remark \ref{rem:A}). This proves the first statement.
Suppose $w\not\in S_{\infty}.$
In order to show 
$\mathfrak{S}_{w}(z,t;0)=0$,
we use the universal localization map 
$\Phi^{A}: \Z[z]\otimes \Z[t] \to \prod_{v \in S_\infty}\Z[t]$ of type $\mathrm{A}$, which is defined in the obvious manner.
A similar proof to Lem. \ref{injective} 
shows that the map
$\Phi^{A}$ is injective.
For any $v\in S_{\infty}$ we have 
$\Phi_{v}(\mathfrak{S}_{w}(z,t;0))
=\Phi_{v}(\mathfrak{S}_{w}(z,t;x))
$, which is equal to $\sigma_{w}^{(\infty)}|_{v}$
by construction of $\mathfrak{S}_{w}(z,t;x).$
Since $v\not\geq w,$ we have
$\sigma_{w}^{(\infty)}|_{v}=0.$
This implies 
that the image of $\mathfrak{S}_{w}$
under the universal localization map
$\Phi^{A}$ is zero, thus $\mathfrak{S}_{w}(z,t;0)=0.$
\end{proof}

\subsection{Divided difference operators and the double Schubert polynomials}
We collect here some properties concerning 
actions of the divided difference operators 
on the double Schubert polynomials.
These will be used in the next section.
\begin{prop}\label{prop:divdifSch}
Let $w=s_{i_{1}}\cdots s_{i_{r}}$ be a reduced expression
of $w\in \W_{\infty}.$
Then the operators
$$\partial_{w}=\partial_{i_{1}}\cdots \partial_{i_{r}},
\quad
\delta_{w}=\delta_{i_{1}}\cdots \delta_{i_{r}}$$
do not depend on the reduced expressions and are 
well-defined for $w\in \W_{\infty}.$ 
Moreover we have
\begin{eqnarray}
\partial_{w}\mathfrak{S}_{u}&=&\begin{cases}
\mathfrak{S}_{uw^{-1}} & \mbox{if}\;\ell(uw^{-1})=\ell(u)-\ell(w)\\
0 &\mbox{otherwise}
\end{cases},\label{eq:partialSch}\\
\delta_{w}\mathfrak{S}_{u}&=&\begin{cases}
\mathfrak{S}_{wu}& \mbox{if}\;\ell(wu)=\ell(u)-\ell(w)
\\
0&\mbox{otherwise}
\end{cases}.
\end{eqnarray}
\end{prop}
\begin{proof}
Since $\{\mathfrak{S}_{u}\}$ is a
$\Z[t]$-basis of $\R_{\infty}$,
the equation (\ref{eq:partialSch}) uniquely determine 
a $\Z[t]$-linear operator, which we denote 
by $\varphi_{w}.$
One can prove 
$\partial_{i_{1}}\cdots \partial_{i_{r}}=\varphi_{w}$
by induction on the length of $w.$
The proof for $\delta_{i}$ is similar.
\end{proof}

\begin{remark}{\rm
The argument here is based on the 
existence of $\{\mathfrak{S}_{w}\}$, but one can also prove it in the classical way - using braid relations (cf. e.g. \cite{BGG}) -
by a direct calculation.}
\end{remark}

The next result will be used in the next section
(Prop. \ref{prop:InterP}).

\begin{lem}\label{Phie}
We have $\Phi_{e}(\partial_{u}\mathfrak{S}_{w})=\delta_{u,w}.$
\end{lem}
\begin{proof} 
First note that 
$\Phi_{e}(\mathfrak{S}_{w})=\sigma_{w}^{(\infty)}|_{e}
=\delta_{w,e}.$
If $\ell(wu^{-1})=\ell(w)-\ell(u)$ is satisfied 
then
by Prop. \ref{prop:divdifSch}, we have
$\Phi_{e}(\partial_{u}\mathfrak{S}_{w})
=\Phi_{e}(\mathfrak{S}_{wu^{-1}})=\delta_{w,u}.$
Otherwise we have $\partial_{u}\mathfrak{S}_{w}=0$
again by Prop. \ref{prop:divdifSch}.
\end{proof}

\subsection{Interpolation formulae and their applications} 
In this section we obtain an explicit combinatorial formula for the double Schubert polynomials, based on the explicit formulas for the single Schubert polynomials from \cite{BH}. The main tool for doing this is the interpolation formula, presented next.
\begin{prop}[Interpolation formula]\label{prop:InterP} For any $f\in \R_{\infty}$, we have
$$
f=\sum_{w\in \W_{\infty}}\Phi_{e}(\partial_{w}(f))\mathfrak{S}_{w}(z,t;x).
$$
\end{prop}
\begin{proof}
Since the double Schubert polynomials
$\{\mathfrak{S}_{w}\}$
form a $\Z[t]$-basis 
of the ring $\R_{\infty}$, we write
$
f=\sum_{w\in \W_{\infty}}c_{w}\mathfrak{S}_{w},\;
c_{w}(t)\in \Z[t].
$
As
$\partial_{w}$ is $\Z[t]$-linear 
we obtain by using Lemma \ref{Phie}
$$
\Phi_{e}(\partial_{w}f)=\sum_{u\in \W_{\infty}}c_{u}(t)
\Phi_{e}(\partial_{w}\mathfrak{S}_{u})
=\sum_{u\in \W_{\infty}}c_{u}(t)
\delta_{w,u}=c_{w}(t).
$$
\end{proof}
\begin{remark}\label{rem:y}{\rm
Let $y=(y_{1},y_{2},\ldots)$ be formal parameters.
On the extended ring $\Z[y]\otimes \R_{\infty}$, 
we can introduce the Weyl group actions, divided difference operators, and the localization map
in the trivial way such that they are $\Z[y]$-linear. 
Since the elements $\mathfrak{S}_{w} \,(w\in \W_{\infty})$ clearly 
form a $\Z[y]\otimes \Z[t]$-basis of $\Z[y]\otimes \R_{\infty},$ the interpolation formula holds also
for any $f\in \Z[y]\otimes \R_{\infty}.$}
\end{remark}

\begin{prop}\label{prop:LS-BH}
Let $y=(y_{1},y_{2},\ldots)$ be formal parameters. Then
$$
\mathfrak{S}_{w}(z,t;x)=\sum_{u,v}\mathfrak{S}_{u}^{A}(y,t)
\mathfrak{S}_{v}(z,y;x)
$$
summed over all $
u\in S_{\infty},\,v\in \W_{\infty}$ such that 
$w=uv,\,
\ell(u)+\ell(v)=\ell(w).$
\end{prop}
\begin{proof}
By the interpolation formula (see Remark \ref{rem:y}), 
we have $$\mathfrak{S}_{w}(z,y;x)=\sum_{v}\Phi_{e}(\partial_{v}\mathfrak{S}_{w}(z,y;x))\mathfrak{S}_{v}(z,t;x).$$
By Prop. \ref{prop:divdifSch}, we see 
$\partial_{v}\mathfrak{S}_{w}(z,y;x)
$ is equal to $\mathfrak{S}_{wv^{-1}}(z,y;x)$
if $\ell(wv^{-1})=\ell(w)-\ell(v)$,
and zero otherwise. Suppose $\ell(wv^{-1})=\ell(w)-\ell(v)$, 
then  
$\Phi_{e}\left(\mathfrak{S}_{wv^{-1}}(z,y;x)\right)=\mathfrak{S}_{wv^{-1}}(t,y;0)$ by the definition of $\Phi_{e}.$
By Lemma \ref{lem:typeA} this is 
$\mathfrak{S}_{wv^{-1}}^{A}(t,y)$
if $wv^{-1}=u\in S_{\infty}$ and zero otherwise.
Then interchanging $t$ and $y$ we have the Proposition.
\end{proof}

Making $y=0$ in the previous proposition, and using that $\mathfrak{S}_u^A(y,t) = \mathfrak{S}_{u^{-1}}^A(-t,-y)$ (cf. Theorem \ref{T:properties} (3) above, for type $\mathrm{A}$ double Schubert polynomials) we obtain:
\begin{cor}\label{cor:typeAexpand} Let $\mathfrak{S}_{w}^{A}(z)$ denote the 
(single) Schubert polynomial of type $\mathrm{A}.$
We have
$$
\mathfrak{S}_{w}(z,t;x)=\sum_{u,v}
\mathfrak{S}_{u^{-1}}^{A}(-t)\mathfrak{S}_{v}(z;x)
$$
summed over all $u\in S_{\infty}, v\in \W_{\infty}$
such that $w=uv$ and $\ell(w)=\ell(u)+\ell(v).$
\end{cor}
There is an explicit combinatorial expression
for the Billey-Haiman polynomials $\mathfrak{S}_{w}(z;x)$
in terms of Schur $Q$-functions and type $\mathrm{A}$ (single)
Schubert polynomials (cf. Thms. 3 and 4 in \cite{BH}).
This, together with the above corollary 
implies also an explicit formula in our case. 
Moreover, the formula for 
$\mathfrak{S}_{w}(z;x)$ is {\em positive},
and therefore this yields a positivity property for the double Schubert polynomials (see Thm. \ref{thm:positivity} below).
We will give an alternative proof for this 
positivity result, independent of the results from {\em loc.cit.}


\subsection{Positivity property} To prove the positivity of the double Schubert polynomials, we begin with the following  lemma (compare with Thms. 3 and 4 in \cite{BH}):
\begin{lem}\label{lem:S00}
We have
$
\mathfrak{S}_{w}(z;x)=
\sum_{u,v}\mathfrak{S}_{u}^{A}(z)\mathfrak{S}_{v}(0,0;x)
$
summed over all $u\in \W_{\infty}, v\in S_{\infty}$
such that $w=uv$ and $\ell(w)=\ell(u)+\ell(v).$
\end{lem}

\begin{remark}{\rm
The function $\mathfrak{S}_{v}(0,0;x)$
is the Stanley's symmetric function
involved in the combinatorial expression for $\mathfrak{S}_{w}(z;x)$
from \cite{BH}. This follows from comparing the present
lemma and the Billey-Haiman's formulas 4.6 and 4.8.}
\end{remark}
\begin{proof}
By (\ref{eq:relBH}) and symmetry
property we have 
$
\mathfrak{S}_{w}(z;x)=
\mathfrak{S}_{w}(z,0;x)
=\mathfrak{S}_{w^{-1}}(0,-z;x).$
Applying Prop. \ref{prop:LS-BH} with $y=0$ 
we can rewrite this as follows:
$$
\sum_{w^{-1}=u^{-1}v^{-1}}\mathfrak{S}_{u^{-1}}^{A}
(0,-z)\mathfrak{S}_{v^{-1}}(0,0;x)
=\sum_{w=vu}\mathfrak{S}_{u}^{A}
(z)\mathfrak{S}_{v}(0,0;x),$$
where the sum is over
$v\in \W_{\infty}, u\in S_{\infty}$ such that
$w^{-1}=u^{-1}v^{-1},$ and $\ell(w^{-1})=\ell(u^{-1})+\ell(v^{-1}).$ The last equality 
follows from symmetry property.
\end{proof}

We are finally ready to prove the positivity property of $\mathfrak{S}_{w}(z,t;x)$. Expand $\mathfrak{S}_{w}(z,t;x)$ as $$
\mathfrak{S}_{w}(z,t;x)=\sum_{\lambda
\in \mathcal{SP}}
f_{w,\lambda}(z,t)
F_{\lambda}(x),
$$ where 
$F_{\lambda}(x)=Q_{\lambda}(x)$ for type
$\mathrm{C}$ and $P_{\lambda}(x)$ for type $\mathrm{D}.$
\begin{thm}[Positivity of double Schubert polynomials]\label{thm:positivity} For any $w \in W_n$, the coefficient $f_{w,\lambda}(z,t)$ is a polynomial in 
$ \mathbb{N}[-t_{1},\ldots,-t_{n-1},
z_{1},\ldots,z_{n-1}]$.
\end{thm}
\begin{proof} The proof follows from the expression on Corollary \ref{cor:typeAexpand}, Lemma \ref{lem:Stanley} below, combined with Lemma \ref{lem:S00} and the fact that
$\mathfrak{S}_{u}^A(z)\in \N[z]$.
\end{proof}

\begin{lem}\label{lem:Stanley}
$\mathfrak{S}_{v}(0,0;x)$
is a linear combination of 
Schur's $Q$- (respectively $P$-) Schur functions
with nonnegative integral coefficients. \end{lem}
\begin{proof} This follows from the transition equations
in \S \ref{ssec:trans}
(see also Remark \ref{rem:TransDSP}).
In fact, the functions
$\mathfrak{S}_{w}(0,0;x)$ satisfy the transition equations 
specialized at $z=t=0$ with the Grassmannian
Schubert classes identified with the Schur's $Q$ or $P$-
functions. 
In fact, the 
recursive formula for $F_{w}(x)=\mathfrak{S}_{w}(0,0;x)$ is positive, in the sense that the right hand side of the 
equation is a certain non-negative integral
linear combination of the functions $\{F_{w}(x)\}.$  
This implies that 
$F_{w}(x)=\mathfrak{S}_{w}(0,0;x)$ can be expressed as
a linear combination of Schur's $Q$ (or $P$)
functions with coefficients in non-negative integers.
\end{proof}



\section{Formula for the longest element}\label{sec:Long}
\setcounter{equation}{0}
In this section, we give explicit formula 
for the double Schubert polynomials
associated with the longest element
$w_{0}^{(n)}$ in $W_{n}$ (and $W_{n}'$).
We note that our proof
of Theorem \ref{existC} is independent of this section.

\subsection{Removable boxes}
We start this section with some combinatorial
properties of factorial $Q$ and $P$-Schur functions.
The goal is to prove Prop. \ref{prop:deltaQ}, which shows how the divided difference operators act on the aforementioned functions.
See \S \ref{ssec:FacSchur} to recall the convention
for the shifted Young diagram $Y_{\lambda}.$

\begin{Def}
A box $x\in Y_\lambda$ is removable
if 
$Y_\lambda-\{x\}$ is again a 
shifted Young diagram of a strict partition.
Explicitly,
$x=(i,j)$ is removable if
$j=\lambda_i+i-1$ and $\lambda_{i+1}\leq \lambda_i-2.$
\end{Def}
To each box $x=(i,j)$ in $Y_\lambda$ we
define its {\it content} $c(x)\in I_\infty$,
$c'(x)\in I_\infty'$
by $c(x)=j-i,$
and
$c'(x)=j-i+1$ if $i\ne j,$
$c'(i,i)=\hat{1} $ if $i$ is odd, and $c'(i,i)=1$ if $i$ is even.
Let $i\in I_\infty$ (resp. $i\in I_\infty'$).
We call $\lambda$ $i$-{\it removable\/} 
if there is a removable box $x$ in $Y_\lambda$ 
such that $c(x)=i$ (resp $c'(x)=i)$.
Note that there is at most one
such $x$ for each $i\in I_\infty$
(resp. $i\in I_\infty').$
We say $\lambda$ is $i$-{\it unremovable\/}
if it is not $i$-{\it removable}.

\setlength{\unitlength}{0.5mm}
\begin{center}
\begin{picture}(100,70)
    \put(5,60){Type $\mathrm{C}$}
  \put(5,55){\line(1,0){50}}
  \put(5,45){\line(1,0){50}}
  \put(15,35){\line(1,0){40}}
  \put(25,25){\line(1,0){20}}
  \put(35,15){\line(1,0){10}}
  \put(5,45){\line(0,1){10}}
  \put(15,35){\line(0,1){20}}
  \put(25,25){\line(0,1){30}}
  \put(35,15){\line(0,1){40}}
  \put(45,15){\line(0,1){40}}
  \put(55,35){\line(0,1){20}}
  \put(8.5,47){\small{$0$}}
  \put(18.5,37){\small{$0$}}
  \put(28.5,27){\small{$0$}}
  \put(18.5,47){\small{$1$}}
  \put(28.5,37){\small{$1$}}
  \put(38.5,27){\small{$1$}}
  \put(38.5,17){\small{$0$}}
  \put(38.5,37){\small{$2$}}
  \put(28.5,47){\small{$2$}}
  \put(38.5,47){\small{$3$}}
  \put(48.5,47){\small{$4$}}
  \put(48.5,37){\small{$3$}}
  \put(5,5){$0$ or $3$\,\mbox{-removable}}
  \end{picture}
  \begin{picture}(100,70)
  \put(5,60){Type $D$}
  \put(5,55){\line(1,0){50}}
  \put(5,45){\line(1,0){50}}
  \put(15,35){\line(1,0){40}}
  \put(25,25){\line(1,0){20}}
  \put(35,15){\line(1,0){10}}
  \put(5,45){\line(0,1){10}}
  \put(15,35){\line(0,1){20}}
  \put(25,25){\line(0,1){30}}
  \put(35,15){\line(0,1){40}}
  \put(45,15){\line(0,1){40}}
  \put(55,35){\line(0,1){20}}
  \put(8.5,47){\small{$\hat{1}$}}
  \put(18.5,37){\small{$1$}}
  \put(28.5,27){\small{$\hat{1}$}}
  \put(18.5,47){\small{$2$}}
  \put(28.5,37){\small{$2$}}
  \put(38.5,27){\small{$2$}}
  \put(38.5,17){\small{$1$}}
  \put(38.5,37){\small{$3$}}
  \put(28.5,47){\small{$3$}}
  \put(38.5,47){\small{$4$}}
  \put(48.5,47){\small{$5$}}
  \put(48.5,37){\small{$4$}}
  \put(5,5){$1$ or $4$\,\mbox{-removable}}
  \end{picture}
  \end{center}
The following facts are well-known (see e.g. \S 7 in \cite{IN}).
\begin{lem} \label{lem:s_iGrass}
Let $w_\lambda\in W_\infty^0$ (resp. $w'_\lambda\in W_\infty^{\hat{1}}$) denote the Grassmannian
element corresponding to $\lambda\in \mathcal{SP}.$
For $i\in I_\infty$ (resp. $i\in I_\infty'$),
a strict partition $\lambda$ is $i$-removable if and only if
$\ell(s_iw_\lambda)=\ell(w_\lambda)-1$
(resp. $\ell(s_iw_\lambda')=\ell(w_\lambda')-1$).
If $\lambda$ is $i$-removable 
then $s_iw_\lambda$ (resp. $s_iw_\lambda'$) is also a
Grassmannian
element and the corresponding 
strict partition is the one obtained from $\lambda$
by removing a (unique) box of content $i$.
\end{lem}

\begin{prop} \label{prop:deltaQ} Let $\lambda$ be a strict
partiton and $i\in I_\infty$ (resp. $i\in I_\infty'$).
\begin{enumerate}
\item \label{remove}  If $\lambda$ is $i$-removable,
then $\delta_{i}Q_{\lambda}(x|t)=Q_{\lambda'}(x|t),$
(resp. $\delta_{i}P_{\lambda}(x|t)=P_{\lambda'}(x|t)$)
where $\lambda'$ is the strict partition obtained by removing 
the (unique) box of content $i$ from $\lambda,$
\item  \label{unremov} If $\lambda$ is $i$-unremovable, then
$\delta_{i}Q_{\lambda}(x|t)=0$
(resp. $\delta_{i}P_{\lambda}(x|t)=0$),
that is to say
$s_{i}^{t}Q_{\lambda}(x|t)=Q_{\lambda}(x|t)$
(resp. $s_{i}^{t}P_{\lambda}(x|t)=P_{\lambda}(x|t)$).
\end{enumerate}
\end{prop}
\begin{proof} This follows from Lemma \ref{lem:s_iGrass} and from the fact that $\mathfrak{C}_{w_\lambda}=Q_\lambda(x|t)$ and $\mathfrak{D}_{w_\lambda'}=P_\lambda(x|t)$, hence we can apply the divided difference equations from Theorem \ref{existC}.
\end{proof}

\subsection{Type $\mathrm{C}_{n}$ case}\label{ssec:LongC}
For $\lambda\in \mathcal{SP}$ we define
$$K_{\lambda}=
K_{\lambda}(z,t;x)
=Q_{\lambda}(x|t_{1},-z_{1},t_{2},-z_{2},\ldots,
t_{n},-z_{n},\ldots).
$$
We need the following two lemmata to 
prove Thm. \ref{thm:Top}.

\begin{lem}\label{lem:deltaK} Set $\Lambda_{n}=\rho_{n}+\rho_{n-1}.$
We have
$
\delta_{n-1}\cdots\delta_{1}\delta_{0}
\delta_{1}\cdots\delta_{n-1}
K_{\Lambda_{n}}=K_{\Lambda_{n-1}}.
$
\end{lem}

\begin{lem}\label{lem:piDelta} We have
$
\pi_{n}(K_{\Lambda_{n}})
=\sigma_{w_{0}^{(n)}}^{(n)}$, where $\pi_n:R_\infty \to \eqcoh (\F_n)$ is the projection defined in \S \ref{ssec:projection}.
\end{lem}

\subsubsection{Proof of Theorem \ref{thm:Top} for type $\mathrm{C}$}
\label{ssec:PfLongC}
\begin{proof}
Let $w_{0}^{(n)}$ in $W_{n}$ be the longest 
element in $W_{n}.$
We need to show that
\begin{equation}
\mathfrak{C}_{w_{0}^{(n)}}(z,t;x)=K_{\Lambda_{n}}
(z,t;x).\label{LongElmFormula}
\end{equation}
Let $w\in W_{\infty}.$ Choose any 
$n$ such that $w\in W_{n}$
and set  
 $$
F_{w}:=\delta_{ww_{0}^{(n)}}K_{\Lambda_{n}}.
$$  Since $\ell(w w_0^{(n)}) + 2n-1 = \ell(w w_0^{(n+1)})$ and $ww_0^{(n)} s_n \cdots s_1s_0s_1 \cdots s_n = ww_0^{(n+1)}$, it follows that $\delta_{ww_0^{(n)}} \cdot \delta_n \cdots \delta_{1}\delta_0 \delta_{1}\cdots \delta_n = \delta_{ww_0^{(n+1)}}$. Then
Lem. \ref{lem:deltaK} yields
$\delta_{ww_{0}^{(n+1)}}K_{\Lambda_{n+1}}
=\delta_{ww_{0}^{(n)}}K_{\Lambda_{n}}$
for any $w\in W_{n}$, so $F_{w}$ is
independent of the choice of $n.$ 
In order to prove the theorem
it is enough to prove $F_{w}=\mathcal{C}_{w}$
for all $w\in W_{\infty}.$ 

By definition of $F_{w}$ 
and basic properties of divided differences 
we can show that
\begin{equation}
\delta_{i}F_{w}=\begin{cases}
F_{s_{i}w}& \ell(s_{i}w)=\ell(w)-1\\
0& \mbox{otherwise}
\end{cases}.\label{eq:divF}
\end{equation}
Now we claim that 
$\pi_{n}(F_{w})=\sigma_{w}^{(n)}$
(for any $n$ such that $w\in W_{n}$).
In fact,
by commutativity of $\pi_{n}$ and
divided difference operators (Prop. \ref{prop:piCommD}), we have
$$
\pi_{n}(F_{w})=\delta_{ww_{0}^{(n)}}\pi_{n}(K_{\Lambda_{n}})
=\delta_{ww_{0}^{(n)}}\sigma_{w_{0}^{(n)}}^{(n)}
=\sigma_{w}^{(n)}.
$$
In the second equality we used 
Lem. \ref{lem:piDelta},
and the last equality
is a consequence of (\ref{eq:divF}). Thus the claim is proved.
Since the claim holds for any 
sufficiently large $n$,
we have $\Phi(F_{w})=\sigma_{w}^{(\infty)}$
(cf. Prop. \ref{prop:StabSch}).
\end{proof}

\subsubsection{Proof of Lemma \ref{lem:deltaK}}
\begin{proof} The
Lemma follows from the successive use of the following
equations (see the example below):
\begin{enumerate} 
\item \label{deltaQ1}
$
\delta_{i}K_{\Lambda_{n}-1^{n-i-1}}
=K_{\Lambda_{n}-1^{n-i}}\quad (0\leq i\leq n-1),
$
\item \label{deltaQ2}
$
\delta_{i}K_{\Lambda_{n}-1^{n}-0^{n-i}1^{i-1}}
=K_{\Lambda_{n}-1^{n}-0^{n-i-1}1^{i}}\quad (1\leq i\leq n-1).
$
\end{enumerate}

We first prove (\ref{deltaQ1}).
For the case $i=0$,
we can apply Prop. \ref{prop:deltaQ} (\ref{remove}) 
directly to get the equation.
Suppose $1\leq i\leq n-1.$
Before applying $\delta_{i}$
to $K_{\Lambda_{n}-1^{n-i-1}}$ 
we switch the parameters
at $(2i-1)$-th and $2i$-th positions to get
$$
K_{\Lambda_{n}-1^{n-i-1}}
=Q_{\Lambda_{n}-1^{n-i-1}}
(x|t_{1},-z_{1},\ldots,-z_{i},t_{i},t_{i+1},-z_{i+1},\ldots,t_{n},-z_{n}).
$$
This is valid in view of Prop. \ref{prop:deltaQ} (\ref{unremov}) and
the fact that 
$\Lambda_{n}-1^{n-i-1}$ is $(2i-1)$-unremovable.
In the right hand side,
the parameters $t_{i}$ and $t_{i+1}$ 
are on $2i$-th and $(2i+1)$-th
positions. Thus the operator $\delta_{i}$ on 
this function is equal to the $2i$-th divided difference
operator ``$\delta_{2i}$'' with respect to 
the sequence of the rearranged parameters 
$(t_{1},-z_{1},\ldots,-z_{i},t_{i},t_{i+1},-z_{i+1},\ldots,t_{n},-z_{n})$
(see example below).
Thus we have by Prop. \ref{prop:deltaQ} (\ref{remove})
$$
\delta_{i}K_{\Lambda_{n}-1^{n-i-1}}
=Q_{\Lambda_{n}-1^{n-i}}(x|t_{1},-z_{1},\ldots,-z_{i},t_{i},t_{i+1},-z_{i+1},\ldots,t_{n},-z_{n}),
$$
namely we remove the box of content $2i$
from $\Lambda_{n}-1^{n-i-1}.$
Then again by Prop. \ref{prop:deltaQ} (\ref{unremov}), the last 
function is equal to
$K_{\Lambda_{n}-1^{n-i}}$;
here we notice $\Lambda_{n}-1^{n-i}$ is $(2i-1)$-unremovable.

Next we prove (\ref{deltaQ2}).
In this case, by Prop. \ref{prop:deltaQ} (\ref{unremov}),
we can switch $2i$-th and $(2i+1)$-th
parameters to get
$$
K_{\Lambda_{n}-1^{n}-0^{n-i}1^{i-1}}
=Q_{\Lambda_{n}-1^{n}-0^{n-i}1^{i-1}}(x|t_{1},-z_{1},\ldots,-z_{i-1},t_{i},t_{i+1},-z_{i},\ldots,
t_{n},-z_{n}).
$$
Here we used the fact that 
$\Lambda_{n}-1^{n}-0^{n-i}1^{i-1}$
is $2i$-unremovable.
Now we apply $\delta_i$ to the function. 
The operator $\delta_i$
is now ``$\delta_{2i-1}$'' 
with respect to 
the sequence of the rearranged parameters 
$(t_{1},-z_{1},\ldots,-z_{i-1},t_{i},t_{i+1},-z_{i},\ldots,
t_{n},-z_{n}).
$
By applying Prop.
\ref{prop:deltaQ} (\ref{remove}),
we have
$$
\delta_i
K_{\Lambda_{n}-1^{n}-0^{n-i}1^{i-1}}
=Q_{\Lambda_{n}-1^{n}-0^{n-i-1}1^{i}}
(x|t_{1},-z_{1},\ldots,-z_{i-1},t_{i},t_{i+1},-z_{i},\ldots,
t_{n},-z_{n}).
$$
The last expression is equal to 
$K_{\Lambda_{n}-1^{n}-0^{n-i-1}1^{i}}$ since 
$\Lambda_{n}-1^{n}-0^{n-i-1}1^{i}$ is $2i$-unremovable.
\end{proof}

{\bf Examples.}
Here we illustrate the process to 
show $\delta_{2}\delta_{1}\delta_{0}\delta_{1}\delta_{2}
K_{\Lambda_{3}}=K_{\Lambda_{2}}$
(case $n=3$ in Lem. \ref{lem:deltaK}).
\setlength{\unitlength}{0.4mm}
\begin{center}
  \begin{picture}(800,50)
  \put(5,45){\line(1,0){50}}
  \put(5,35){\line(1,0){50}}
  \put(15,25){\line(1,0){30}}
  \put(25,15){\line(1,0){10}}
  \put(5,35){\line(0,1){10}}
  \put(15,25){\line(0,1){20}}
  \put(25,15){\line(0,1){30}}
  \put(35,15){\line(0,1){30}}
  \put(45,25){\line(0,1){20}}
  \put(55,35){\line(0,1){10}}
  \put(8.5,38){\small{$0$}}
  \put(18.5,28){\small{$0$}}
  \put(28.5,18){\small{$0$}}
  \put(18.5,38){\small{$1$}}
  \put(28.5,28){\small{$1$}}
  \put(38.5,28){\small{$2$}}
  \put(28.5,38){\small{$2$}}
  \put(38.5,38){\small{$3$}}
  \put(48.5,38){\small{$4$}}
  \put(75,45){\line(1,0){40}}
  \put(75,35){\line(1,0){40}}
  \put(85,25){\line(1,0){30}}
  \put(95,15){\line(1,0){10}}
  \put(75,35){\line(0,1){10}}
  \put(85,25){\line(0,1){20}}
  \put(95,15){\line(0,1){30}}
  \put(105,15){\line(0,1){30}}
  \put(115,25){\line(0,1){20}}
  \put(78.5,38){\small{$0$}}
  \put(88.5,28){\small{$0$}}
  \put(98.5,18){\small{$0$}}
  \put(88.5,38){\small{$1$}}
  \put(98.5,28){\small{$1$}}
  \put(108.5,28){\small{$2$}}
  \put(98.5,38){\small{$2$}}
  \put(108.5,38){\small{$3$}}
  \put(130,0){
  \put(5,45){\line(1,0){40}}
  \put(5,35){\line(1,0){40}}
  \put(15,25){\line(1,0){20}}
  \put(25,15){\line(1,0){10}}
  \put(5,35){\line(0,1){10}}
  \put(15,25){\line(0,1){20}}
  \put(25,15){\line(0,1){30}}
  \put(35,15){\line(0,1){30}}
  \put(45,35){\line(0,1){10}}
  \put(8.5,38){\small{$0$}}
  \put(18.5,28){\small{$0$}}
  \put(28.5,18){\small{$0$}}
  \put(18.5,38){\small{$1$}}
  \put(28.5,28){\small{$1$}}
  \put(28.5,38){\small{$2$}}
  \put(38.5,38){\small{$3$}}
  }
  \put(190,0){
  \put(5,45){\line(1,0){40}}
  \put(5,35){\line(1,0){40}}
  \put(15,25){\line(1,0){20}}
  \put(5,35){\line(0,1){10}}
  \put(15,25){\line(0,1){20}}
  \put(25,25){\line(0,1){20}}
  \put(35,25){\line(0,1){20}}
  \put(45,35){\line(0,1){10}}
 \put(8.5,38){\small{$0$}}
  \put(18.5,28){\small{$0$}}
  \put(18.5,38){\small{$1$}}
  \put(28.5,28){\small{$1$}}
  \put(28.5,38){\small{$2$}}
  \put(38.5,38){\small{$3$}}
  }
\put(250,0){
  \put(5,45){\line(1,0){40}}
  \put(5,35){\line(1,0){40}}
  \put(15,25){\line(1,0){10}}
  \put(5,35){\line(0,1){10}}
  \put(15,25){\line(0,1){20}}
  \put(25,25){\line(0,1){20}}
  \put(35,35){\line(0,1){10}}
 \put(45,35){\line(0,1){10}}
  \put(8.5,38){\small{$0$}}
  \put(18.5,28){\small{$0$}}
  \put(18.5,38){\small{$1$}}
  \put(28.5,38){\small{$2$}}
  \put(38.5,38){\small{$3$}}
  }
  \put(310,0){
  \put(5,45){\line(1,0){30}}
  \put(5,35){\line(1,0){30}}
  \put(15,25){\line(1,0){10}}
  \put(5,35){\line(0,1){10}}
  \put(15,25){\line(0,1){20}}
  \put(25,25){\line(0,1){20}}
  \put(35,35){\line(0,1){10}}
  \put(8.5,38){\small{$0$}}
  \put(18.5,28){\small{$0$}}
  \put(18.5,38){\small{$1$}}
 \put(28.5,38){\small{$2$}}}
 \put(56,35){$\longrightarrow$}
 \put(116,35){$\longrightarrow$}
 \put(176,35){$\longrightarrow$}
 \put(236,35){$\longrightarrow$}
 \put(296,35){$\longrightarrow$}
 \put(60,43){\small{$\delta_{2}$}}
 \put(120,43){\small{$\delta_{1}$}}
 \put(180,43){\small{$\delta_{0}$}}
 \put(240,43){\small{$\delta_{1}$}}
 \put(300,43){\small{$\delta_{2}$}}
 \put(19,5){\small{$K_{5,3,1}$}}
  \put(89,5){\small{$K_{4,3,1}$}}
  \put(149,5){\small{$K_{4,2,1}$}}
  \put(209,5){\small{$K_{4,2}$}}
  \put(269,5){\small{$K_{4,1}$}}
  \put(319,5){\small{$K_{3,1}$}}
  \end{picture}
\end{center}
We pick up the first arrow:
$\delta_{2}K_{5,3,1}=K_{4,3,1}$
(equation (\ref{deltaQ1}) in
Lem. \ref{lem:deltaK} for $n=3,\,i=2$).
As is indicated in the proof, we divide
this equality into the following 
four steps: 
\begin{eqnarray*}
K_{5,3,1}=Q_{5,3,1}(x|t_{1},-z_{1},\underline{t_{2},-z_{2}},t_{3},-z_{3})
\underset{(a)}{=}
Q_{5,3,1}(x|t_{1},-z_{1},-z_{2},t_{2},t_{3},-z_{3})\\
\overset{\delta_{2}}{\longrightarrow}
Q_{4,3,1}(x|t_{1},-z_{1},\underline{-z_{2},t_{2}},t_{3},-z_{3})
\underset{(b)}{=}Q_{4,3,1}(x|t_{1},-z_{1},t_{2},-z_{2},t_{3},-z_{3})=K_{4,3,1}.
\end{eqnarray*}
In the equality $(a)$ we used the fact that 
$\Lambda_{3}=(5,3,1)$ is $3$-unremovable,
so the underlined pair of variables can be exchanged
(by Prop. \ref{prop:deltaQ}, (2)).
Then we apply $\delta_{2}$ to this function.
Note that the variables $t_{2},t_{3}$ are in the
$4$-th and $5$-th positions in the parameters of 
the function. So if we rename the parameters as 
$
f=Q_{5,3,1}(x|t_{1},-z_{1},-z_{2},t_{2},t_{3},-z_{3})
=Q_{5,3,1}(x|u_{1},u_{2},u_{3},u_{4},u_{5},u_{6}),$
then $\delta_{2}$ is ``$\delta_{4}$''
with respect to the parameter sequence $(u_{i})_{i}.$
Namely we have
$$
\delta_{2}f=\frac{f-s_{2}^{t}f}{t_{3}-t_{2}}
=\frac{f-s_{4}^{u}f}{u_{5}-u_{4}},
$$
where $s_{4}^{u}$ exchanges $u_{4}$ and $u_{5}.$
Since $\Lambda_{3}=(5,3,1)$ is $4$-removable,
we see from Prop. \ref{prop:deltaQ}, (1) that 
$\delta_{2}=$``$\delta_{4}$'' removes the box of content $4$
from $(5,3,1)$ to obtain the shape $(4,3,1).$
Then finally, in the equality $(b)$, we exchange the variables
$-z_{2},t_{2}$ again using Prop. \ref{prop:deltaQ}, (2).
This is valid since $(4,3,1)$ is $3$-unremovable.
Thus we obtained $K_{4,3,1}.$

\subsubsection{Proof of Lemma \ref{lem:piDelta}}
\begin{proof}
We calculate $\Phi_{v}(K_{\Lambda_{n}})$ for $v\in W_{n}.$
Recall that the map $\Phi_{v}: R_{\infty}\rightarrow
\Z[t]$ is the $\Z[t]$-algebra homomorphism 
given by $x_{i}\mapsto t_{v,i}$
and $z_{i}\mapsto t_{v(i)}.$
So we have
$$
\Phi_{v}(K_{\Lambda_{n}})
=Q_{\Lambda_{n}}
(t_{v,1},\ldots,t_{v,n}|t_{1},-{t_{v(1)}},
\ldots, t_{n},-{t_{v(n)}}).
$$
Note that $t_{v,i}=0$ for $i>n$ since $v$ is an element
in $W_{n}.$
From the factorization formula \ref{lem:factorization},
this is equal to
$$\prod_{1\leq i\leq n}2t_{v,i}
\prod_{1\leq i<j\leq n}(t_{v,i}+t_{v,j})
\times s_{\rho_{n-1}}
(t_{v,1},\ldots,t_{v,n}|t_{1},-{t_{v(1)}},
\ldots, t_{n},-{t_{v(n)}}).
$$
The presence of the  
factor $\prod_{i}2t_{v,i}$
implies that
$\Phi_{v}(K_{\Lambda_{n}})$ vanishes unless
$v(1),\ldots,v(n)$ are all negative.
So from now on we assume
$v=(\overline{\sigma(1)},\ldots,\overline{\sigma(n)})$
for some permutation $\sigma \in S_{n}.$
Then we have $t_{v,i}=t_{\sigma(i)}$
and $t_{v(i)}=-t_{\sigma(i)}$ so the last factor
of factorial Schur polynomial becomes
$$s_{\rho_{n-1}}
(t_{\sigma(1)},\ldots,t_{\sigma(n)}|t_{1},{t_{\sigma(1)}},
\ldots, t_{n},{t_{\sigma(n)}}).$$
This is equal to
$s_{\rho_{n-1}}
(t_{1},\ldots,t_{n}|t_{1},{t_{\sigma(1)}},
\ldots, t_{n},{t_{\sigma(n)}})$ because 
$s_{\rho_{n-1}}$ is symmetric in ther first set of variables.
From Lem. \ref{lem:A-long} we know that this polynomial factors into  
$\prod_{1\leq i<j\leq n}
(t_{j}-t_{\sigma(i)}).$
This is zero except for the case $\sigma=\mathrm{id},$
namely $v=w_{0}^{(n)}.$
If $\sigma=\mathrm{id}$ then
$\Phi_{v}(K_{\Lambda_{n}})$ becomes 
$\prod_{1\leq i\leq n}2t_{i}
\prod_{1\leq i<j\leq n}(t_{i}+t_{j})
\prod_{1\leq i<j\leq n}
(t_{j}-t_{i})=\sigma_{w_{0}^{(n)}}^{(n)}|_{w_{0}^{(n)}}.$
\end{proof}

\subsection{Type $\mathrm{D}_{n}$ case}
Set 
$K'_{\lambda}(z,-t;x)=
P_{\lambda}
(x|t_{1},-z_{1},\ldots,t_{n-1},-z_{n-1},\ldots).
$
Our goal in this section is
$$
\mathfrak{D}_{w_{0}^{(n)}}
=K'_{2\rho_{n-1}}
(z,t;x).
$$
We use the same strategy 
as in \S \ref{ssec:LongC} to
prove this.
Actually the proof in \S \ref{ssec:PfLongC}
works also in this case
using the following two lemmata,
which will be proved below.

\begin{lem}\label{lem:DeltaKD} We have
$\delta_{n-1}\cdots\delta_{2}\delta_{\hat{1}}\delta_{1}\delta_{2}
\cdots\delta_{n-1}K_{2\rho_{n-1}}'
=K_{2\rho_{n-2}}'.$
\end{lem}

\begin{lem}\label{lem:piDeltaD} We have
$\pi_{n}(K'_{2\rho_{n-1}})=\sigma_{w^{(n)}_{0}}^{(n)}.$
\end{lem}

\subsubsection{A technical lemma}

We need the following technical lemma which is used
in the proof of Lem. \ref{lem:DeltaKD}.
Throughout the section, 
$(u_{1},u_{2},u_{3},\ldots)$
denote any sequence of variables independent of $t_{1},t_{2}.$
%
\begin{lem}\label{lem:hat1}
Let $\lambda=(\lambda_{1},\ldots,\lambda_{r})$
be a strict partition such that $r$ is odd and $\lambda_{r}\geq 3.$
Set $\tilde{t}=(u_{1},t_{1},t_{2},u_{2},u_{3},\ldots).$
Then
$
\delta_{\hat{1}}P_{\lambda_{1},\ldots,\lambda_{r},1}(x|\tilde{t})
=P_{\lambda_{1},\ldots,\lambda_{r}}(x|\tilde{t}).
$
\end{lem}
\begin{sublem}\label{lem:deltahat1} Suppose 
$\lambda$ is $1, 2$ and $\hat{1}$-unremovable.
Then we have
$
\delta_{\hat{1}}
P_{\lambda}(x|\tilde{t})=0.\label{eq:deltahat1}
$
\end{sublem}
\begin{proof}
Since 
$\lambda$ 
is $1,2$-unremovable,
we can rearrange the first three parameters
by using Prop. \ref{prop:deltaQ}, so we have
$P_{\lambda}(x|\tilde{t})=P_{\lambda}(x|t_{1},t_{2},u_{1},u_{2},u_{3},\ldots).$
Because $\lambda$ is also $\hat{1}$-unremovable,
it follows that $\delta_{\hat{1}}P_{\lambda}(t_{1},t_{2},u_{1},u_{4},\ldots)=0$ from Prop. \ref{prop:deltaQ}.
\end{proof}

\begin{sublem}[Special case of Lem. \ref{lem:hat1} for $r=1$]\label{lem:Pk1}
We have $
\delta_{\hat{1}}P_{k,1}(x|\tilde{t})=P_{k}(x|\tilde{t})$
for $k\geq 3.$
\end{sublem}
\begin{proof}
Substituting $\tilde{t}$ for $t$ into (\ref{eq:P2row})
we have 
$$P_{k,1}(x|\tilde{t})=
P_{k}(x|\tilde{t})P_{1}(x|\tilde{t})
-P_{k+1}(x|\tilde{t})-(u_{k-1}+u_{1})P_{k}(x|\tilde{t}).$$
By the explicit formula $P_{1}(x|\tilde{t})=P_{1}(x),$
we have 
$\delta_{\hat{1}}P_{1}(x|\tilde{t})=1.$
We also have $\delta_{\hat{1}}P_{k}(x|\tilde{t})=\delta_{\hat{1}}P_{k+1}(x|\tilde{t})=0$
by Sublemma \ref{lem:deltahat1}.
Then we use the Leibnitz rule
$\delta_{\hat{1}}(fg)=\delta_{\hat{1}}(f)g+(s_{\hat{1}}f)\delta_{\hat{1}}(g)$
to get $\delta_{\hat{1}}P_{k,1}(x|\tilde{t})=P_{k}(x|\tilde{t}).$
\end{proof}

\begin{proof}[Proof of Lem. \ref{lem:hat1}.]
From the definition of the Pfaffian
it follows that 
$$
P_{\lambda_{1},\ldots,\lambda_{r},1}(x|\tilde{t})
=\sum_{j=1}^{r}(-1)^{r-j}P_{\lambda_{j},1}(x|\tilde{t})P_{\lambda_{1},\ldots,
\widehat{\lambda_{j}},\ldots,\lambda_{r}}(x|\tilde{t}).
$$
Then the Leibnitz rule
combined with Sublemma
\ref{eq:deltahat1} and Sublemma \ref{lem:Pk1}
implies
$$
\delta_{\hat{1}}
P_{\lambda_{1},\ldots,\lambda_{r},1}(x|\tilde{t})
=\sum_{j=1}^{r}(-1)^{r-j}P_{\lambda_{j}}(x|\tilde{t})P_{\lambda_{1},\ldots,
\widehat{\lambda_{j}},\ldots,\lambda_{r}}(x|\tilde{t})
=P_{\lambda_{1},\ldots,\lambda_{r}}(x|\tilde{t}),
$$
where in the last equality we used the expansion
formula of Pfaffian again. \end{proof}

\subsubsection{Proof of Lem. \ref{lem:DeltaKD}}

\begin{proof}
Consider the case when $n$ is even.
By applying the same method of calculation as in
type $\mathrm{C}$ case, we have 
$$
\delta_{1}\delta_{2}\cdots\delta_{n-1}K'_{2\rho_{n-1}}
=P_{\rho_{n-1}+\rho_{n-2}}(x|-z_{1},t_{1},t_{2},-z_{2},t_{3},-z_{3},\ldots).
$$
The problem here is that $\rho_{n-1}+\rho_{n-2}$ is not
$1$-unremovable when $n$ is even.
So we can not rewrite the function
as $K_{\rho_{n-1}+\rho_{n-2}}'.$
Nevertheless, by using Lem. \ref{lem:hat1}, 
we can show 
$$
\delta_{\hat{1}}P_{\rho_{n-1}+\rho_{n-2}}(x|-z_{1},t_{1},t_{2},-z_{2},t_{3},-z_{3},\ldots)
=K'_{\rho_{n-1}+\rho_{n-2}-0^{n-2}1}.
$$
The rest of calculation is similar to type $\mathrm{C}$ case.
If $n$ is odd, we can show this equation
using only Prop. \ref{prop:deltaQ}
as in type $\mathrm{C}$ case.
\end{proof}

\subsubsection{Proof of Lem. \ref{lem:piDeltaD}}
\begin{proof}
Similar to the proof of Lem. \ref{lem:piDelta} using
Lem. \ref{lem:A-long}, \ref{lem:A-longOdd},
\ref{lem:factorizationD}, and 
\ref{lem:factorizationDodd}.

We calculate $\Phi_{v}(K'_{2\rho_{n-1}})$ for $v\in W'_{n}.$
We have
$$
\Phi_{v}(K'_{2\rho_{n-1}})
=P_{2\rho_{n-1}}
(t_{v,1},\ldots,t_{v,n}|t_{1},-{t_{v(1)}},
\ldots, t_{n},-{t_{v(n)}}).
$$
Note that $t_{v,i}=0$ for $i>n$ since $v$ is an element
in $W'_{n}.$

Assume now that $n$ is even.
From the factorization formula (Lem. \ref{lem:factorizationD}),
this is equal to
$$
\prod_{1\leq i<j\leq n}(t_{v,i}+t_{v,j})
\times s_{\rho_{n-1}}
(t_{v,1},\ldots,t_{v,n}|t_{1},-{t_{v(1)}},
\ldots, t_{n},-{t_{v(n)}}).
$$
The factorial Schur polynomial 
factorizes further into linear terms by Lem. \ref{lem:A-long},
and we finally obtain $$
\Phi_{v}(K'_{2\rho_{n-1}})=
\prod_{1\leq i<j\leq n}(t_{v,i}+t_{v,j})
\prod_{1\leq i<j\leq n}
(t_{j}+t_{v(i)}).
$$
We set
$v=(\overline{\sigma(1)},\ldots,\overline{\sigma(n)})$
for some $\sigma \in S_{n}$ since otherwise 
$t_{v,i}=t_{v,j}=0$ for some $i,j$ with $i\ne j$
and then $\prod_{1\leq i<j\leq n}
(t_{v,i}+t_{v,j})$ vanishes.
Then
the factor $\prod_{1\leq i<j\leq n}
(t_{j}+t_{v(i)})$ is $\prod_{1\leq i<j\leq n}
(t_{j}-t_{\sigma(i)}).$
This is zero except for the case $\sigma=\mathrm{id},$
namely $v=w_{0}^{(n)}.$
If $w=w_{0}^{(n)}$, we have
$\Phi_{v}(K'_{\rho_{n-1}})=
\prod_{1\leq i<j\leq n}(t_{i}+t_{j})
\prod_{1\leq i<j\leq n}
(t_{j}-t_{i})=\sigma_{w_{0}^{(n)}}^{(n)}|_{w_{0}^{(n)}}.$

Next we consider the case when $n$ is odd.
Note that the longest element $w_{0}^{(n)}$ in this case
is $1\bar{2}\bar{3}\cdots\bar{n}.$
Let $s(v)$ denote the number of
nonzero entries in $t_{v,1},\ldots,t_{v,n}.$ 
Then we have 
$s(v)\leq n-1$ since $v\in W_{n}'$ and $n$ is odd. 
 We use the following identity:
\begin{eqnarray}
&&P_{2\rho_{n-1}}(x_{1},\ldots,x_{n-1}|t_{1},-z_{1},\ldots,t_{n-1},-z_{n-1})\nonumber\\
&=&\prod_{1\leq i<j\leq n-1}(x_{i}+x_{j})
\times s_{\rho_{n-1}+1^{n-1}}(x_{1},\ldots,x_{n-1}|t_{1},-z_{1},\ldots,t_{n-1},-z_{n-1}).\label{eq:facD}
\end{eqnarray}
If $s(v)<n-1$,
then $s(v)\leq n-3$ because $v\in W_{n}'.$
This means that there are at least $3$
zeros in $t_{v,1},\ldots,t_{v,n}.$
Because there is the factor
$\prod_{1\leq i<j\leq n-1}(x_{i}+x_{j})$
in (\ref{eq:facD}) we have $\Phi_{v}(K'_{2\rho_{n-1}})=0.$
So we suppose $s(v)=n-1.$
By a calculation using the definition of the factorial
Schur polynomial, we see that 
$ s_{\rho_{n-1}+1^{n-1}}(x_{1},\ldots,x_{n-1}|t_{1},-z_{1},\ldots,
t_{n-1},-z_{n-1})$
is divisible by the factor $\prod_{i=1}^{n-1}(t_{1}-x_{i}).$
By this fact we may assume
$t_{v,1},\ldots,t_{v,n}$ is a permutation of
$0,t_{2},t_{3},\ldots,t_{n}$ since
otherwise $\Phi_{v}(K'_{2\rho_{n-1}})$ is zero.
Thus under the assumption, we have
$$
\Phi_{v}(K'_{2\rho_{n-1}})
=\prod_{2\leq i<j\leq n}(t_{i}+t_{j})
s_{\rho_{n-1}+1^{n-1}}(t_{2},\ldots,t_{n}|t_{1},-t_{v(1)},t_{2},-t_{v(2)},\ldots,
t_{n-1},-t_{v(n-1)}).
$$
By Lem. \ref{lem:A-longOdd}
this factorizes into 
$
\prod_{2\leq i<j\leq n}(t_{i}+t_{j})
\prod_{j=2}^{n}(t_{j}-t_{1})
\prod_{1\leq i<j\leq n}(t_{j}+t_{v(i)}).
$
Now our assumption is that
the negative elements in 
$\{v(1),\ldots,v(n)\}$ are exactly $\{2,3,\ldots,n\}.$
Among these elements, only 
$w_{0}^{(n)}$ gives a non-zero
polynomial, which 
is shown to be 
$\prod_{1\leq i<j\leq n}(t_{i}+t_{j})(t_{j}-t_{i})
=\sigma_{w_{0}^{(n)}}^{(n)}|_{w_{0}^{(n)}}.$
\end{proof}


\section{Geometric construction 
of the universal localization map}\label{sec:geometry}
\setcounter{equation}{0}
In this section, we construct the morphism of $\Z[t]$-algebras $\tilde{\pi}_{\infty}:\R_\infty \longrightarrow \invlim\, \eqcoh(\F_n)$
from a geometric point of view.
We start this section by describing 
the embedding $\F_{n}\hookrightarrow \F_{n+1}$
explicitly, and calculate the localization 
of the Chern roots of tautological bundles.
Then we introduce some particular cohomology classes
$\beta_{i}$ in $H_{T_{n}}^{*}(\F_{n}) $, by using the geometry of isotropic flag varieties. These classes satisfy the relations of the $Q$-Schur functions $Q_{i}(x)$, and mapping $Q_i(x)$ to $\beta_i$ ultimately leads to the homomorphism 
$\tilde{\pi}_{\infty}.$ In particular, this provides an explanation on why
the Schur $Q$-functions enter into our 
theory (cf. Prop. \ref{prop:pi}).
The final goal is to establish the connection
of $\tilde{\pi}_{\infty}$ and the universal 
localization map $\Phi$ (Thm. \ref{thm:piQpiQ}).

The arguments in the preceding sections are 
logically independent from this section.
However, we believe that the results in this section provide
the reader with some 
insight into the underlying geometric idea
of the algebraic construction.

\subsection{Flag varieties of isotropic flags}\label{ssection:flagv} The groups $G_n$ are the group of automorphisms preserving a non-degenerate, bilinear form $\langle \cdot, \cdot \rangle$ on a complex vector space $V_n$. The pair $(V_n, \langle \cdot , \cdot \rangle)$ is the following:
\begin{enumerate}
\item In type $\mathrm{C}_n$, $V_n = \C^{2n}$; fix $\e_n^{*}, \dots, \e_1^{*}, \e_1, \dots ,\e_n$ a basis for $V_n$. Then $\langle\cdot , \cdot \rangle$ is the skew-symmetric form given by $\langle\e_i, \e_j^*\rangle = \delta_{i,j}$ (the ordering of the basis elements will be important later, when we will embed $G_n$ into $G_{n+1}$).
\item In types $\mathrm{B}_n$ and $\mathrm{D}_n$, $V_n$ is an odd, respectively even-dimensional complex vector space. Let $\e_n^{*}, \dots , \e_1^{*}, \e_0, \e_1, \dots , \e_n$ respectively $\e_n^{*}, \dots , \e_1^{*}, \e_1, \dots , \e_n$ be a basis of $V_n$. Then $\langle \cdot , \cdot \rangle$ is the symmetric form such that $\langle\e_i,\e_j^*\rangle=\delta_{i,j}$.
\end{enumerate}

A subspace $V$ of $V_n$ will be called {\em isotropic} if $\langle \vecu,\vecv\rangle  = 0 $ for any $\vecu,\vecv \in V$. Then $\F_n$ is the variety consisting of complete isotropic flags with respect to the appropriate bilinear form. For example, in type $\mathrm{C}_n$, $\F_n$ consists of nested sequence of vector spaces
\[ F_1 \subset F_2 \subset \dots \subset F_n \subset V_n = \C^{2n} \/, \] such that each $F_i$ is isotropic and $\dim F_i =i$. Note that the maximal dimension of an isotropic subspace of $\C^{2n}$ is $n$; but the flag above can be completed to a full flag of $\C^{2n}$ by taking $V_{n+i}= V_{n-i}^\perp$, using the non-degeneracy of the form $\langle \cdot , \cdot \rangle$. A similar description can be given in types $\mathrm{B}_n$ and $\mathrm{D}_n$, with the added condition that, in type $\mathrm{D}_n$, \[\dim F_n \cap \langle \e_n^{*}, \dots , \e_1^*\rangle \equiv 0 \mod 2 \/; \] in this case we say that all $F_n$ are {\em in the same family} (cf. \cite[pag.68]{FP}).

The flag variety $\F_{n}$ carries 
a transitive left action of the group $G_n,$
and can be identified with the homogeneous space $G_n/B_n$, where $B_{n}$ is the Borel subgroup consisting of
upper triangular matrices in $G_{n}$.
Let $T_{n}$ be the maximal torus 
in $G_{n}$ consisting of
diagonal matrices in $G_{n}.$
Let $t = \diag(\xi_n^{-1}, \dots ,\xi_1^{-1},\xi_1, \dots , \xi_n)$ be a torus element in types $\mathrm{C_n,D_n}$, and $t= \diag(\xi_n^{-1}, \dots ,\xi_1^{-1},1,\xi_1, \dots , \xi_n)$ in type $\mathrm{B}_n$. We denote by $t_{i}$ the character 
of $T_{n}$ defined by $t\mapsto \xi_{i}^{-1}\;(t\in T_{n}).$
Then the weight of $\C \,\e_i$ is $-t_i$ and that of $\C\, \e_i^*$ is  $t_i$. 
We identify $t_i \in H^2_{T_n}(pt)$ with $c_1^T(\C\e_i^*)$, where $\C\e_i^*$ is the (trivial, but not equivariantly trivial) line bundle over $pt$
with fibre $\C\e_i^*$. For $v \in W_n$, the corresponding $T_{n}$-fixed point ${e}_v$ is \[{e}_v: \langle  \e_{v(n)}^{*} \rangle \subset \langle   \e_{v(n)}^{*}, \e_{v(n-1)}^{*}\rangle \subset \dots \subset \langle \e_{v(n)}^{*}, \e_{v(n-1)}^{*}, \dots , \e_{v(1)}^{*}\rangle \subset V_n \/. \]

\subsection{Equivariant embeddings of flag varieties}\label{ssec:embeddings} 
There is a natural embedding $G_n \hookrightarrow G_{n+1}$, given explicitly by
$$
g\rightarrow
\left(\begin{array}{c|c|c}
1 & & \\
\hline
 & g & \\
 \hline
 & & 1
\end{array}\right).
$$ 
This corresponds to the embedding of Dynkin diagrams in each type.
This also induces embeddings $B_{n}\hookrightarrow B_{n+1},$ $T_{n}\hookrightarrow T_{n+1},$ and ultimately $\varphi_n
:\F_{n}\hookrightarrow \F_{n+1}.$ 
The embedding $\varphi_n$ sends 
the complete isotropic flag $F_1 \subset \cdots \subset F_n$ of $V_n$ to the complete isotropic flag of $V_{n+1}= 
\C\, \e_{n+1}^{*} \oplus V_{n}\oplus \C\, \e_{n+1}$: 
\[ \C\, \e_{n+1}^* \subset \C\, \e_{n+1}^* \oplus F_1 \subset \cdots \subset \C\, \e_{n+1}^* \oplus F_n \/. \]
Cleary $\varphi_n$ is 
equivariant with respect to the embedding 
$T_{n}\hookrightarrow T_{n+1}.$

\subsection{Localization of Chern classes of tautological bundles} 
Consider the flag of tautological (isotropic) vector bundles
$$
0=\mathcal{V}_{n+1}\subset \mathcal{V}_{n}\subset
\cdots
\subset \mathcal{V}_{1}\subset \mathcal{E},
\quad \mathrm{rank}\,\mathcal{V}_{i}=n-i+1,
$$
where $\mathcal{E}$ is the trivial bundle
with fiber $V_{n}$ and
$\mathcal{V}_{i}$ is defined to be the vector 
subbundle of $\mathcal{E}$ whose 
fiber over the point $F_{\bullet}= F_1 \subset \dots \subset F_n$ in $\mathcal{F}_{n}$ 
is $F_{n-i+1}.$
Let $z_{i}=c_{1}^{T}(\mathcal{V}_{i}/\mathcal{V}_{i+1})$
\begin{footnote}{The bundle $\mathcal{V}_i/\mathcal{V}_{i+1}$ is in fact negative; for example, in type $\mathrm{C}$, if $n=1$, $\F_1=\mathbb{P}^1$ and $\mathcal{V}_1 = \mathcal{O}(-1)$. The reason for choosing positive sign for $z_i$ is to be consistent with the conventions used by Billey-Haiman in \cite{BH}.}\end{footnote}denote the equivariant Chern class of the line bundle
$\mathcal{V}_{i}/\mathcal{V}_{i+1}$.
\begin{prop}\label{prop:loc} Let $v \in W_n$. Then the localization map $\iota_v^*: \eqcoh (\F_n) \to \eqcoh ({e}_v)$ satisfies $\iota_v^*(z_i) = t_{v(i)}$. \end{prop}
\begin{proof} The pull-back of the line bundle $\mathcal{V}_{i}/\mathcal{V}_{i+1}$ via $\iota_v^*$ is the line bundle over ${e}_v$ with fibre $\C\e_{v(i)}^{*}$, which has (equivariant) first Chern class $t_{v(i)}$. \end{proof}

\subsection{The cohomology class $\beta_{i}$} 
In this section we introduce the cohomology classes
$\beta_{i}$, which will later be identified to $Q$-Schur functions $Q_i(x)$.


The torus action on $V_{n}$ induces a $T_{n}$-equivariant splitting $\mathcal{E}=\oplus_{i=1}^{n}\mathcal{L}_{i}\oplus\mathcal{L}_{i}^{*}$ ($\mathcal{E}=\oplus_{i=1}^{n}\mathcal{L}_{i}\oplus\mathcal{L}_{i}^{*}\oplus\mathcal{L}_{0}$ for type $\mathrm{B}_{n}$)
where $\mathcal{L}_{i}$ (resp. $\mathcal{L}_{i}^{*}$) is the trivial line bundle over $\F_n$
with fiber $\C\e_{i}$ (resp. $\C\e_{i}^{*}$). Recall from \S \ref{ssection:flagv} that $T_n$ acts on $\mathcal{L}_i^*$ by weight $t_i$ and that $t_i = c_1^T(\mathcal{L}_i^*)$. 

Let $\F_n$ be the flag variety of type $\mathrm{C}_{n}$ or $\mathrm{D}_{n}$ and
set $\mathcal{V}=\mathcal{V}_{1}$.
We have the following 
exact sequence of $T_{n}$-equivariant 
vector bundles:
\begin{equation}
0\longrightarrow \mathcal{V}\longrightarrow
\mathcal{E}
\longrightarrow
\mathcal{V}^{*}\longrightarrow 0.\label{tauto}
\end{equation}
where $\mathcal{V}^{*}$ denotes the dual bundle
of $\mathcal{V}$ in $\mathcal{E}$ with respect to 
the bilinear form. 
Let $\mathcal{L}=\oplus_{i=1}^{n}\mathcal{L}_{i}$
and $\mathcal{L}^{*}=\oplus_{i=1}^{n}\mathcal{L}_{i}^{*}.$
Since $\mathcal{E}=\mathcal{L}\oplus
\mathcal{L}^{*},$ 
we have $c^{T}(\mathcal{E})=c^{T}(\mathcal{L})c^{T}(\mathcal{L}^{*}).$ Define
the class $\beta_{i}\in H_{T_{n}}^{*}(\F_{n})$ by 
$$
\beta_{i}=c_{i}^{T}(\mathcal{V}^{*}-\mathcal{L}),
$$
where $c_{i}^{T}(\mathcal{A}-\mathcal{B})$
is the term of degree $i$ in the formal expansion
of $c^{T}(\mathcal{A})/c^{T}(\mathcal{B}).$
Using the relation $c^{T}(\mathcal{L})c^{T}(\mathcal{L}^{*})
=c^{T}(\mathcal{V})c^{T}(\mathcal{V}^{*})$,
we also have the expression:
$$
\beta_{i}=c_{i}^{T}(\mathcal{L}^{*}-\mathcal{V}).
$$
In terms of the Chern classes $z_{i},t_{i}$, 
the class $\beta_{i}$ has
the following two equivalent expressions:
\begin{equation}
\sum_{i=0}^{\infty}\beta_{i}u^{i}=\prod_{i=1}^{n}\frac{1-z_{i}u}{1-t_{i}u}=
\prod_{i=1}^{n}\frac{1+t_{i}u}{1+z_{i}u}.\label{eq:twoexpr}
\end{equation}
\begin{lem}\label{lem:beta}
The classes $\beta_{i}$ satisfy
the same relations as the $Q$-Schur functions
of $Q_{i}(x)$, i.e.
$$
\beta_{i}^{2}+2\sum_{j=1}^{i}(-1)^{j}\beta_{i+j}\beta_{i-j}=0\quad
\mbox{for}\; i\geq 1. 
$$ 
\end{lem}
\begin{proof}
We have the following two expressions:
$$
\sum_{i=0}^{\infty}\beta_{i}u^{i}
=
\prod_{i=1}^{n}\frac{1-z_{i}u}{1-t_{i}u},\quad
\sum_{j=0}^{\infty}(-1)^{j}\beta_{j}u^{j}=
\prod_{i=1}^{n}\frac{1-t_{i}u}{1-z_{i}u}.
$$
The lemma follows from multiplying both sides, and then extracting 
the degree $2i$ parts.
\end{proof}

Minor modifications need to be done if $\F_n$ is the flag variety of type $\mathrm{B}_{n}$. In this case the tautological sequence of isotropic flag subbundles consists of $
0=\mathcal{V}_{n+1}\subset \mathcal{V}_{n}\subset
\cdots
\subset \mathcal{V}_{1}\subset  \mathcal{E}= \C^{2n+1} \times \F_n$, but the dual bundle $\mathcal{V}_1^*$ of $\mathcal{V}_1$ is not isomorphic to $\mathcal{E}/\mathcal{V}_1$, which has rank $n+1$. 
However, the line bundle $\mathcal{V}_1^\perp/\mathcal{V}_1$ is equivariantly isomorphic to $\bigwedge^{2n+1} \mathcal{E}$ - cf. \cite[pag.75]{FP} - so $c_1^T(\mathcal{V}_1^\perp/\mathcal{V}_1)=0$; here $\mathcal{V}_1^\perp$ denotes the bundle whose fibre over $V_1 \subset \cdots \subset V_n$ is the subspace of vectors in $\C^{2n+1}$ perpendicular to those in $V_n$ with respect to the non-degenerate form $\langle \cdot , \cdot \rangle$.
It follows that the bundle $\mathcal{E}/\mathcal{V}_{1}$ has (equivariant) total Chern class $(1-z_1u) \cdots (1-z_nu)$, which is the same as the total Chern class of $\mathcal{V}_{1}^*.$ Similarly, the total Chern class of $\mathcal{E}/\mathcal{L}$
with $\mathcal{L}=\oplus_{i=1}^{n}\mathcal{L}_{i}$ is $(1+t_1u)\cdots (1+t_nu)$ and equals $c^T(\mathcal{L}^*)$. So the definition of $\beta_i$ and the proofs of its properties remain unchanged.


Recall that in \S \ref{ssec:projection} we introduced $\pi_{n}:\R_{\infty}
\longrightarrow H_{T_{n}}^{*}(\mathcal{F}_{n})$
by using the universal localization map $\Phi$. The following is the key fact used in the proof of the main 
result of this section.
\begin{lem}\label{lem:key} We have
$
\pi_{n}(Q_{i}(x))=\beta_{i}.
$
\end{lem}
\begin{proof}
It is enough to show that
$\iota_{v}^{*}(\beta_{i})=Q_{i}(t_{v})
$ for $v\in \W_{n}.$
By Prop. \ref{prop:loc} and the definition of 
$\beta_{i}$, we have
$$
\iota_{v}^{*}\left(
\sum_{i=0}^{\infty}\beta_{i}u^{i}
\right)
=\iota_{v}^{*}\left(\prod_{i=1}^{n}\frac{1-z_{i}u}{1-t_{i}u}
\right)
=\prod_{i=1}^{n}\frac{1-t_{v(i)}u}{1-t_{i}u}.$$
If $v(i)$ is positive, the factors $1-t_{v(i)}u$ cancel out and 
the last expression 
becomes
$$
\prod_{v(i)\;{\small \mbox{negative}}}\frac{1-t_{v(i)}u}{1+t_{v(i)}u}
=\sum_{i=0}^{\infty}
Q_{i}(t_{v})u^{i}$$
where the last equality follows 
from the definition of $Q_{i}(x)$ and
that of $t_{v}.$
\end{proof}

\subsection{Homomorphism $\tilde{\pi}_{n}$}
We consider $\F_{n}$ of one of the types
$\mathrm{B}_{n},\mathrm{C}_{n},$ and $\mathrm{D}_{n}.$
We will define next the projection homomorphism
from $\R_{\infty}$ to $H_{T_{n}}^{*}(\F_{n})$, which will be used to construct the geometric analogue $\tilde{\pi}_n$ of $\pi_n$.
Note that $R_{\infty}$ is a proper subalgebra of $\R_{\infty}$ in types $\mathrm{B}$ and $\mathrm{D}.$
We regard $H_{T_{n}}^{*}(\F_{n})$
as $\Z[t]$-module via the natural projection
$\Z[t]\rightarrow \Z[t_{1},\ldots,t_{n}].$
\begin{prop}\label{prop:pi}
There exists a homomorphism of graded
$\Z[t]$-algebras
$\tilde{\pi}_{n}:R_{\infty}\rightarrow H_{T_{n}}^{*}(\F_{n})$
such that
$$
\tilde{\pi}_{n}(Q_{i}(x))=\beta_{i}\quad
(i\geq 1)
\quad\mbox{,}\quad
\tilde{\pi}_{n}(z_{i})=z_{i}\quad
(1\leq i\leq n)\quad\mbox{and}\quad
\tilde{\pi}_{n}(z_{i})=0\quad(i>n).$$
\end{prop}
\begin{proof} This follows from the fact that $R_\infty$ is generated as a $\Z[t]$-algebra by $Q_{i}(x),z_{i}\;(i\geq 1)$, and that the ideal of relations among $Q_i(x)$ is generated by those in (\ref{eq:quadraticQ}) (see \cite{Mac}, III \S 8). Since the elements $\beta_i$ satisfy also those relations by Lemma \ref{lem:beta}, the result follows.\end{proof}

\subsection{Types $\mathrm{B}$ and $\mathrm{D}$} In this section, we extend $\tilde{\pi}_{n}$
from $R_{\infty}'$ to $H_{T_{n}}^{*}(\F_{n})$. The key to that is the identity $P_i(x) = \frac{1}{2} Q_i(x)$.


%

%
\begin{prop} Let 
$\F_{n}$ be the flag variety of type $\mathrm{B}_{n}$ or $\mathrm{D}_{n}.$
Then there is an (integral) cohomology class
$\gamma_{i}$ such that 
$2\gamma_{i}=\beta_{i}.$
Moreover, the classes $\gamma_i$ satisfy the following quadratic relations:
$$
\gamma_{i}^{2}+2\sum_{j=1}^{i-1}(-1)^{j}
\gamma_{i+j}\gamma_{i-j}
+(-1)^{i}\gamma_{2i}=0 \quad(i>0).
$$
\end{prop}
\begin{proof} Define $\gamma_i=\frac{1}{2} \beta_i$. Then, as in the proof of Lemma \ref{lem:key}, the localization $\iota^*_v(\gamma_i) = \frac{1}{2} Q_i(t_v) = P_i(t_v)$ which is a polynomial with integer coefficients.
The quadratic relations follow immediately from Lem. \ref{lem:beta}. 
\end{proof}
The proposition implies immediately the following: 
\begin{prop}\label{prop:piBD}
Let $\F_{n}$ be the flag variety of type $\mathrm{B}_{n}$ or $\mathrm{D}_{n}.$
There exists a homomorphism of graded
$\Z[t]$-algebras
$\tilde{\pi}_{n}:R_{\infty}'\rightarrow H_{T_{n}}^{*}(\F_{n})$
such that
$$
\tilde{\pi}_{n}(P_{i}(x))=\gamma_{i}\quad
(i\geq 1)
\quad\mbox{and}\quad
\tilde{\pi}_{n}(z_{i})=z_{i}\quad
(1\leq i\leq n)\quad\mbox{and}\quad
\tilde{\pi}_{n}(z_{i})=0\quad(i>n).$$
\end{prop}
\begin{remark} {\rm It is easy to see (cf, \cite[\S 6.2]{FP}) that the morphism $\tilde{\pi}_n:R_\infty \to \eqcoh(\F_n)$ is surjective in type $\mathrm{C}$, and also in types $\mathrm{B,D}$, but with coefficients over $\Z[1/2]$. But in fact, using that $\Phi:R_\infty' \to H_\infty$ is an isomorphism, one can show that surjectivity holds over $\Z$ as well.}\end{remark}
\subsection{The geometric interpretation of the universal localization map $\Phi$}
From Prop. \ref{prop:pi} and Prop. \ref{prop:piBD},
we have $\Z[t]$-algebra homomorphism 
$
\tilde{\pi}_{n}: \R_{\infty}\longrightarrow
H_{T_{n}}^{*}(\F_{n})
$
for all types $\mathrm{B,C,D}.$
Since $\tilde{\pi}_n$ is compatible with maps $\varphi_{n}^{*}:\eqcoh(\F_{n+1}) \to \eqcoh(\F_n)$ induced by embeddings $\F_n \to \F_{n+1}$
there is an induced homomorphism $$\tilde{\pi}_{\infty}:\R_\infty \longrightarrow \invlim \eqcoh (\F_n).$$
Recall from \S \ref{ssec:projection} that we have the natural embedding
$\pi_\infty: R_{\infty}\hookrightarrow 
\invlim \eqcoh (\F_n)$, defined via the localization map $\Phi$. Then:
\begin{thm}\label{thm:piQpiQ}
We have that $\tilde{\pi}_{\infty}=\pi_\infty$.
\end{thm}
\begin{proof}
It is enough to show that $\tilde{\pi}_{n}=\pi_{n}$. To do that, we
compare both maps on the generators
of $\R_{\infty}.$
We know that $\tilde{\pi}_{n}(Q_{i}(x))={\pi}_{n}(Q_{i}(x))
=\beta_{i}$
by Lem. \ref{lem:key} and this implies that
$\tilde{\pi}_{n}(P_{i}(x))={\pi}_{n}(P_{i}(x))$
for types $\mathrm{B}_{n}$ and $\mathrm{D}_{n}.$
It remains to show
$\pi_{n}(z_{i})=\tilde{\pi}_{n}(z_{i}).$
In this case, for $v\in \W_{n}$, 
\[
\iota_{v}^{*}\pi_{n}(z_{i})
=\Phi_{v}(z_{i})^{(n)}=t_{v(i)}^{(n)}
=\iota_{v}^{*}\tilde{\pi}_{n}(z_{i})\/.
\]
This completes the proof.
\end{proof}

\subsection{Integrality of Fulton's classes $c_{i}$}
We take to opportunity to briefly discuss, in the present setting, an integrality property of some
cohomology classes considered by Fulton in \cite{F}, in relation to degeneracy loci in classical types. This property was proved before, in a more general setting, by Edidin and Graham \cite{EG}, by using the geometry of quadric bundles.

Let $\mathcal{F}_{n}$ be the flag variety of type
$\mathrm{B}_{n}$ or $\mathrm{D}_{n}.$ Recall that $c_i$ is the equivariant cohomology class in $H_{T_{n}}^{*}(\F_{n},\Z[{\textstyle\frac{1}{2}}])$ defined by:
$$
c_{i}={\textstyle\frac{1}{2}}
\left(
e_{i}(-z_{1},\ldots,-z_{n})+e_{i}(t_{1},\ldots,t_{n})
\right)\quad(0\leq i\leq n).
$$
\begin{prop}
We have
$
c_{i}=
\sum_{j=1}^{i-\epsilon}(-1)^{j}e_{j}(t_{1},\ldots,t_{n})\gamma_{i-j}\;(0\leq i\leq n),
$
where $\epsilon=0$ if $i$ is even and $\epsilon=1$
if $i$ is odd. 
In particular, $c_{i}$ are classes defined over $\Z$.
\end{prop}
\begin{proof}
Using the definition of $\beta_{i}$
and $
\sum_{i=0}^{n}2c_{i}u^{i}
=\prod_{i=1}^{n}(1-z_{i}u)
+\prod_{i=1}^{n}(1+t_{i}u)
$
we have
$$\sum_{i=0}^{n}2c_{i}u^{i}
=\prod_{i=1}^{n}(1+t_{i}u)+
\sum_{i=0}^{\infty}\beta_{i}u^{i}\prod_{i=1}^{n}(1-t_{i}u).
$$
Comparing both hand side of degree $i$ and using $2\gamma_{i}=\beta_{i}$
we have the equation.
\end{proof}

\section{Kazarian's formula for Lagrangian Schubert classes}
\label{sec:Kazarian}
\setcounter{equation}{0}
In this section, we give a brief discussion 
of a ``multi-Schur Pfaffian'' expression for the Schubert classes
of the Lagrangian Grassmannian. 
This formula appeared in a preprint of Kazarian, regarding a
degeneracy loci formula for the Lagrangian 
vector bundles \cite{Ka}.
\subsection{Multi-Schur Pfaffian} 
We recall the definition of the multi-Schur Pfaffian from \cite{Ka}.  
Let $\lambda=(\lambda_{1}>\cdots>\lambda_{r}\geq 0)$ 
be any strict partition with $r$ even.
Consider an $r$-tuple of infinite sequences
$c^{(i)}=\{c^{(i)}_k\}_{k=0}^\infty \;(i=1,\ldots,r)$,
where each $c^{(i)}_k$ is an element in
a commutative ring with unit.
For $a\geq b\geq 0$, we set
$$
c_{a,b}^{(i),(j)}
:=c_{a}^{(i)}
c_{b}^{(j)}
+2\sum_{k=1}^{b}(-1)^{k}
c_{a+k}^{(i)}
c_{b-k}^{(j)}.
$$
Assume that the matrix $(c_{\lambda_i,\lambda_j}^{(i),(j)})_{i,j}$
is skew-symmetric, i.e. 
$c_{\lambda_i,\lambda_j}^{(i),(j)}
=-c_{\lambda_j,\lambda_i}^{(j),(i)}$ for
$1\leq i,j\leq r.$
Then we consider its Pfaffian
$$\mathrm{Pf}_\lambda(c^{(1)},\ldots,c^{(r)})
=
\mathrm{Pf}\left(
c_{\lambda_i,\lambda_j}^{(i),(j)}
\right)_{1\leq i<j\leq r},
$$
called {\it multi-Schur Pfaffian}.

\subsection{Factorial Schur functions 
as a multi-Schur Pfaffian}
We introduce the following 
versions of factorial $Q$-Schur functions $Q_k(x|t)$:
$$
\sum_{k=0}^\infty Q_k^{(l)}(x|t)u^k=
\sum_{i=1}^\infty 
\frac{1+x_iu}{1-x_iu}
\prod_{j=1}^{l-1}(1-t_ju).
$$
Note that, by definition, 
$Q_k^{(k)}(x|t)=Q_k(x|t)$ and
$Q_k^{(1)}(x|t)=Q_k(x).$

\begin{prop} \label{prop:multiPf} Let $\lambda=(\lambda_{1}>\cdots>\lambda_{r}\geq 0)$ 
be any strict partition with $r$ even.
Set $c^{(i)}_k=Q_{k}^{(\lambda_i)}(x|t)$
for $i=1,\ldots,r.$
Then the matrix $(c_{\lambda_i,\lambda_j}^{(i),(j)})_{i,j}$
is skew-symmetric
and 
we have 
$$
\mathrm{Pf}_\lambda(c^{(1)},\ldots,c^{(r)})
=Q_\lambda(x|t).
$$
\end{prop}
\begin{proof}
In view of the Pfaffian formula for $Q_{\lambda}(x|t)$ (Prop. \ref{prop:Pf}),
it suffices to show the following identity:
\begin{equation}
Q_{k,l}(x|t)
=Q_{k}^{(k)}(x|t)Q_{l}^{(l)}(x|t)
+2\sum_{i=1}^{l}(-1)^{i}Q_{k+i}^{(k)}(x|t)Q_{l-i}^{(l)}(x|t).\label{Fkl}
\end{equation}
By induction we can show that for $k\geq 0$
\begin{eqnarray*}
Q_{j}^{(j+k)}(x|t)&=&\sum_{i=0}^{j}(-1)^{i}
e_{i}(t_{j+k-1},t_{j+k-2},\ldots,t_{j-i+1})Q_{j-i}(x|t),\\
Q_{j}^{(j-k)}(x|t)&=&\sum_{i=0}^{k}
h_{i}(t_{j-k},t_{j-k+1},\ldots,t_{j-i})Q_{j-i}(x|t).
\end{eqnarray*}
Substituting these 
expressions into (\ref{Fkl}), we get
a quadratic expression 
in $Q_i(x|t)$'s.
The obtained expression
coincides with a formula for 
$Q_{k,l}(x|t)$ proved in \cite[Prop.7.1]{Ik}.
\end{proof}

\subsection{Schubert classes in the Lagrangian Grassmannian
as multi-Pfaffians} We use the notations from \S \ref{sec:geometry}. 
The next formula expresses the equivariant Schubert class $\sigma_{w_{\lambda}}^{(n)}$ in a flag variety of type $\mathrm{C}$ in terms of a multi-Pfaffian. Recall that this is also the equivariant Schubert class for the Schubert variety indexed by $\lambda$ in the Lagrangian Grassmannian, so this is a "Giambelli formula" in this case. Another such expression, in terms of ordinary Pfaffians, was proved by the first author in \cite{Ik}.
\begin{prop}[cf. \cite{Ka}, Thm. 1.1] Set $\mathcal{U}_{k}
=\oplus_{j=k}^n \mathcal{L}_i.$
Then $
\sigma_{w_{\lambda}}^{(n)}
=\mathrm{Pf}_{\lambda}(c^{T}(\mathcal{E}-\mathcal{V}-\mathcal{U}_{\lambda_{1}}),\ldots,
c^{T}(\mathcal{E}-\mathcal{V}-\mathcal{U}_{\lambda_{r}})).$
\end{prop}
\begin{proof} By Thm. \ref{PhiFacQ}, we know
$\pi_{n}(Q_{\lambda}(x|t))=\sigma_{w_{\lambda}}^{(n)}.$
On the other hand the formula of Prop. \ref{prop:multiPf}
writes $Q_{\lambda}(x|t)$ as a multi-Pfaffian. 
So it is enough to show that:
$$
c^{T}_{i}(\mathcal{E}-\mathcal{V}-\mathcal{U}_{k})=
\pi_{n}(Q_{i}^{(k)}(x|t)).
$$
We have
$$
c^{T}(\mathcal{E}-\mathcal{V}-\mathcal{U}_{k})
=\frac{\prod_{i=1}^n(1-t_i^2u^{2})}{
\prod_{j=1}^{n}(1+z_{i}u)\prod_{j=k}^n(1-t_ju)}
=\prod_{i=1}^n
\frac{1+t_iu}{1+z_{i}u}
\prod_{j=1}^{k-1}(1-t_ju).
$$
The first factor of the 
right hand side is the generating 
function for $\beta_{i}=\pi_{n}(Q_{i}(x))
(i\geq 0).$
So the last expression is 
$$
\sum_{i=0}^{\infty}\pi_{n}(Q_{i}(x))u^{i}
\prod_{j=1}^{k-1}(1-t_ju)
=\sum_{i=0}^{\infty}\pi_{n}(Q_{i}^{(k)}(x|t))u^{i}.
$$
Hence the proposition is proved.
\end{proof}



\section{Type $\mathrm{C}$ double Schubert polynomials for $w \in W_3$}

\tiny
\renewcommand{\arraystretch}{1.2}

$\begin{array}{|c|l|}
\hline
123 & 1 \\
\hline
\123 & Q_1 \\
\hline
213 & Q_1+(z_1-t_1) \\
\hline
\213 & Q_2+Q_1(-t_1) \\
\hline
2\13 & Q_2+Q_1z_1 \\
\hline
\2\13 & Q_{21} \\
\hline
1\23 & Q_3+Q_2(z_1-t_1)+Q_1(-z_1t_1) \\
\hline
\1\23 & Q_{31}+Q_{21}(z_1-t_1) \\
\hline
132 & Q_{1}+(z_1+z_2-t_1-t_2) \\
\hline
\132 & 2Q_{2}+Q_{1}(z_1+z_2-t_1-t_2) \\
\hline
312 & Q_2+Q_1(z_1-t_1-t_2)+(z_1-t_1)(z_1-t_2) \\
\hline
\312 & Q_3+Q_2(-t_1-t_2)+Q_1 t_1t_2 \\
\hline
3\12 & Q_{3}+Q_{21}+Q_{2}(2z_1-t_1-t_2)+Q_{1}(z_1)(z_1-t_1-t_2) \\
\hline
\3\12 & Q_{31}+Q_{21}(-t_1-t_2) \\
\hline
1\32 & Q_{4}+Q_{3}(z_1-t_1-t_2)+Q_{2}(t_1t_2-z_1(t_1+t_2))+Q_1z_1t_1t_2 \\
\hline
\1\32 & Q_{41}+Q_{31}(z_1-t_1-t_2)+Q_{21}(t_1t_2-z_1(t_1+t_2)) \\
\hline
231 & Q_2+Q_{1}(z_1+z_2-t_1)+(z_1-t_1)(z_2-t_1) \\
\hline
\231 & Q_{3}+Q_{21}+Q_{2}(z_1+z_2-2t_1)+Q_{1}(-t_1)(z_1+z_2-t_1) \\
\hline
321 & Q_{3}+Q_{21}+Q_{2}(2z_1+z_2-2t_1-t_2)+Q_{1}(z_1+z_2-t_1)(z_1-t_1-t_2)+(z_1-t_1)(z_1-t_2)(z_2-t_1) \\
\hline
\321 & Q_{4}+Q_{31}+Q_{3}(z_1+z_2-2t_1-t_2)+Q_{21}(-t_1-t_2)+\\ & Q_{2}(t_1t_2-(t_1+t_2)(z_1+z_2-t_1)+Q_{1}t_1t_2(z_1+z_2-t_1) \\
\hline
3\21 & Q_{31}+Q_{3}(z_1-t_1)+Q_{21}(z_1-t_1)+Q_{2}(z_1-t_1)^2+Q_1 z_1(-t_1)(z_1-t_1) \\
\hline
\3\21 & Q_{32}+Q_{31}(-t_1)+Q_{21} t_1^2 \\
\hline
2\31 & Q_{41}+Q_{4}(z_1-t_1)+Q_{31}(z_1-t_1-t_2)+Q_{3}(z_1-t_1)(z_1-t_1-t_2)+Q_{21}(t_1t_2-z_1(t_1+t_2))+\\ & 
Q_{2}(z_1-t_1)(t_1t_2-z_1(t_1+t_2))+Q_{1}(z_1-t_1)z_1t_1t_2  \\
\hline
\2\31 & Q_{42}+Q_{32}(z_1-t_1-t_2)+Q_{41}(-t_1)+Q_{31}(-t_1)(z_1-t_1-t_2)+Q_{21}t_1^2(z_1-t_2) \\
\hline
23\1 & Q_{3}+Q_{2}(z_1+z_2)+Q_{1} z_1z_2 \\
\hline
\23\1 & Q_{31}+Q_{21}(z_1+z_2) \\
\hline
32\1 & Q_{4}+Q_{31}+Q_{3}(2z_1+z_2-t_1-t_2)+Q_{21}(z_1+z_2)+\\ & Q_{2}((z_1+z_2)(z_1-t_1-t_2)+z_1z_2)+Q_{1}z_1z_2(z_1-t_1-t_2) \\
\hline
\32\1 & Q_{32}+Q_{41}+Q_{31}(z_1+z_2-t_1-t_2)+Q_{21}(z_1+z_2)(-t_1-t_2) \\
\hline
3\2\1 & Q_{32}+Q_{31}z_1+Q_{21}z_1^2 \\
\hline
\3\2\1 & Q_{321} \\
\hline
2\3\1 & Q_{42}+Q_{32}(z_1-t_1-t_2)+Q_{41}z_1+Q_{31} z_1(z_1-t_1-t_2)+Q_{21}z_1^2(-t_1-t_2) \\
\hline
\2\3\1 & Q_{421}+Q_{321}(z_1-t_1-t_2) \\
\hline
13\2 & Q_{4}+Q_{3}(z_1+z_2-t_1)+Q_{2}(z_1z_2-t_1(z_1+z_2))+Q_1(-t_1)z_1z_2  \\
\hline
\13\2 & Q_{41}+Q_{31}(z_1+z_2-t_1)+Q_{21}(z_1z_2-t_1(z_1+z_2)) \\
\hline
31\2 & Q_{41}+Q_{4}(z_1-t_1)+Q_{31}(z_1+z_2-t_1)+Q_3(z_1-t_1)(z_1+z_2-t_1)+Q_{21}(z_1z_2-t_1(z_1+z_2))+\\ & Q_{2}(z_1-t_1)(z_1z_2-t_1(z_1+z_2))+Q_{1}(z_1-t_1)z_1z_2(-t_1)  \\
\hline
\31\2 & Q_{42}+Q_{32}(z_1+z_2-t_1)+Q_{41}(-t_1)+Q_{31}(z_1+z_2-t_1)(-t_1)+Q_{21}t_1^2(z_1+z_2) \\
\hline
3\1\2 & Q_{42}+Q_{41}z_1+Q_{32}(z_1+z_2-t_1)+Q_{31}z_1(z_1+z_2-t_1)+Q_{21}z_1^2(z_2-t_1)  \\
\hline
\3\1\2 & Q_{421}+Q_{321}(z_1+z_2-t_1) \\
\hline

\end{array}$

$\begin{array}{|c|l|}
\hline
1\3\2 & Q_{43}+Q_{42}(z_1-t_1)+Q_{32}(z_1^2+t_1^2-z_1t_1)+Q_{41}(-z_1t_1)+Q_{31}z_1(-t_1)(z_1-t_1)+Q_{21}(z_1^2t_1^2) \\
\hline
\1\3\2 & Q_{431}+Q_{421}(z_1-t_1)+Q_{321}(z_1^2-z_1t_1+t_1^2) \\
\hline
12\3 & Q_{5}+Q_4(z_1+z_2-t_1-t_2)+Q_3(z_1z_2+t_1t_2-(z_1+z_2)(t_1+t_2))+\\ & Q_2(z_1z_2(-t_1-t_2)+t_1t_2(z_1+z_2) )+Q_1 z_1z_2t_1t_2 \\
\hline
\12\3 & Q_{51}+Q_{41}(z_1+z_2-t_1-t_2)+Q_{31}(z_1z_2+t_1t_2-(z_1+z_2)(t_1+t_2))+Q_{21}(z_1z_2(-t_1-t_2)+t_1t_2(z_1+z_2)) \\
\hline



21\3 & Q_{51}+Q_{5}(z_1-t_1)+Q_{41}(z_1+z_2-t_1-t_2)+Q_{4}(z_1-t_1)(z_1+z_2-t_1-t_2)+ \\ & Q_{31}(z_1z_2+t_1t_2-(z_1+z_2)(t_1+t_2)) +Q_{3}(z_1-t_1)(z_1z_2+t_1t_2-(z_1+z_2)(t_1+t_2))+ \\ & Q_{21}(z_1z_2(-t_1-t_2)+t_1t_2(z_1+z_2))  +Q_2(z_1-t_1)(z_1z_2(-t_1-t_2)+t_1t_2(z_1+z_2))+Q_1 z_1z_2t_1t_2(z_1-t_1) \\
\hline
\21\3 & Q_{52}+Q_{42}(z_1+z_2-t_1-t_2)+Q_{32}(z_1z_2+t_1t_2-(z_1+z_2)(t_1+t_2))+Q_{51}(-t_1)+\\ & Q_{41}(-t_1)(z_1+z_2-t_1-t_2)+Q_{31}(-t_1)(z_1z_2+t_1t_2-(z_1+z_2)(t_1+t_2)) \\
 & +Q_{21}(t_1)^2(z_1z_2-(z_1+z_2)t_2)) \\
\hline

2\1\3 & Q_{52}+Q_{42}(z_1+z_2-t_1-t_2)+Q_{32}(z_1z_2+t_1t_2-(z_1+z_2)(t_1+t_2))+Q_{51}z_1+Q_{41}z_1(z_1+z_2-t_1-t_2)+\\ & Q_{31}z_1(z_1z_2+t_1t_2-(z_1+z_2)(t_1+t_2))  +Q_{21} z_1^2(t_1t_2-z_2(t_1+t_2)) \\
\hline
\2\1\3 & Q_{521}+Q_{421}(z_1+z_2-t_1-t_2)+Q_{321}(z_1z_2+t_1t_2-(z_1+z_2)(t_1+t_2)) \\
\hline
1\2\3 & Q_{53}+Q_{52}(z_1-t_1)+Q_{51}(-z_1t_1)+Q_{43}(z_1+z_2-t_1-t_2)+Q_{42}(z_1-t_1)(z_1+z_2-t_1-t_2)+\\ & Q_{41}(-z_1t_1)(z_1+z_2-t_1-t_2) +Q_{32}(z_1(z_1-t_1)(z_2-t_1-t_2)+t_1^2(z_2-t_2))+\\ & Q_{31}(-z_1t_1)(z_1(z_2-t_1-t_2)-t_1(z_2-t_2))+Q_{21}z_1^2t_1^2(z_2-t_2) \\
\hline
\1\2\3 & Q_{531}+Q_{431}(z_1+z_2-t_1-t_2)+Q_{521}(z_1-t_1)+
Q_{421}(z_1-t_1)(z_1+z_2-t_1-t_2)+\\ & Q_{321}((z_1^2-z_1t_1+t_1^2)(z_2-t_2)+z_1t_1(t_1-z_1)) \\
\hline
\end{array}
$

\normalsize

\section{Double Schubert polynomials in type $\mathrm{D}$ for $w \in W'_3$}

\tiny
\renewcommand{\arraystretch}{1.2}
$\begin{array}{|c|l|}
\hline
123 & 1 \\
\hline
213 & P_1+(z_1-t_1) \\
\hline
\2\13 & P_1 \\
\hline
\1\23 & P_2+P_1(z_1-t_1) \\
\hline
132 & 2P_1+(z_1+z_2-t_1-t_2) \\
\hline
312 & P_2+P_1(2z_1-t_1-t_2)+(z_1-t_1)(z_1-t_2) \\
\hline
\3\12 & P_2+P_1(-t_1-t_2) \\
\hline
\1\32 & P_3+P_2(z_1-t_1-t_2)+P_1(t_1t_2-z_1t_1-z_1t_2) \\
\hline
231 & P_2+P_1(z_1+z_2-2t_1)+(z_1-t_1)(z_2-t_1) \\
\hline
321 & P_{3}+P_{21}+P_2(2z_1+z_2-2t_1-t_2)+P_1(z_1^2+2z_1z_2+t_1^2+2t_1t_2-3t_1z_1-t_1z_2-t_2z_1-t_2z_2)+\\ & (z_1-t_1)(z_1-t_2)(z_2-t_1) \\
\hline
\3\21 & P_{21}+P_2(-t_1)+P_1t_1^2 \\
\hline
\2\31 & P_{31}+P_{21}(z_1-t_1-t_2)+P_3(-t_1)+P_2(-t_1)(z_1-t_1-t_2)+P_1 t_1^2(z_1-t_2) \\
\hline
\23\1 & P_2+P_1(z_1+z_2) \\
\hline
3\2\1 & P_{21}+P_{2}z_1+P_{1}z_1^2 \\
\hline
\32\1 & P_3+P_{21}+P_2(z_1+z_2-t_1-t_2)+P_1(z_1+z_2)(-t_1-t_2) \\
\hline
2\3\1 & P_{31}+P_{3}z_1+P_{21}(z_1-t_1-t_2)+P_2 z_1(z_1-t_1-t_2)+P_1 z_1^2(-t_1-t_2) \\
\hline
\13\2 & P_3+P_2(z_1+z_2-t_1)+P_1(z_1z_2-t_1(z_1+z_2)) \\
\hline
3\1\2 & P_{31}+P_{21}(z_1+z_2-t_1)+P_3 z_1+P_2z_1(z_1+z_2-t_1)+P_1 z_1^2(z_2-t_1) \\
\hline
\31\2 & P_{31}+P_{21}(z_1+z_2-t_1)+P_3(-t_1)+P_2(-t_1)(z_1+z_2-t_1)+P_1(z_1+z_2)t_1^2 \\
\hline
1\3\2 & P_{32}+P_{31}(z_1-t_1)+P_3(-z_1t_1)+P_{21}(z_1^2-z_1t_1+t_1^2)+P_2(-z_1t_1)(z_1-t_1)+P_1z_1^2t_1^2 \\
\hline
\12\3 & P_4+P_3(z_1+z_2-t_1-t_2)+P_2(z_1z_2+t_1t_2-(z_1+z_2)(t_1+t_2))+P_1(z_1z_2(-t_1-t_2)+t_1t_2(z_1+z_2)) \\
\hline
2\1\3 & P_{41}+P_{4}z_1+P_{31}(z_1+z_2-t_1-t_2)+P_3(z_1)(z_1+z_2-t_1-t_2)+P_{21}(z_1z_2+t_1t_2-(z_1+z_2)(t_1+t_2)) \\
 & +P_2z_1(z_1z_2+t_1t_2-(z_1+z_2)(t_1+t_2))+P_1z_1^2(t_1t_2-z_2t_1-z_2t_2) \\
\hline
\21\3 & P_{41}+P_{4}(-t_1)+P_{31}(z_1+z_2-t_1-t_2)+P_{3}(-t_1)(z_1+z_2-t_1-t_2)+P_{21}(z_1z_2+t_1t_2-(z_1+z_2)(t_1+t_2)) \\
 & +P_2(-t_1)(z_1z_2+t_1t_2-(z_1+z_2)(t_1+t_2))+P_1t_1^2(z_1z_2-z_1t_2-z_2t_2) \\
\hline
1\2\3 & P_{42}+P_{32}(z_1+z_2-t_1-t_2)+P_{41}(z_1-t_1)+P_{31}(z_1-t_1)(z_1+z_2-t_1-t_2)+ \\ & P_{21}(z_1^2z_2-t_1^2t_2+z_1t_1t_2-z_1z_2t_1+z_1^2(-t_1-t_2)+t_1^2(z_1+z_2)) +P_4(-z_1t_1)+\\ & P_3(-z_1t_1)(z_1+z_2-t_1-t_2)+P_2(-z_1t_1)(-z_1t_1-z_2t_1-z_1t_2+z_1z_2+t_1t_2)+P_1(z_1^2t_1^2)(z_2-t_2) \\
\hline
\end{array}
$

\end{document}